\documentclass[a4paper]{article}
\usepackage[english]{babel}
\usepackage[utf8]{inputenc}
\usepackage{amsmath,amsthm}
\usepackage{amsfonts,amssymb}

\usepackage{tikz,tikz-cd} 

\usepackage{xifthen} 

\newcommand{\toposdefaut}{0}
\newcommand{\topos}[1][\toposdefaut]{ 
\ifthenelse{\equal{#1}{0}}{ \mathcal{T} }
{
\ifthenelse{\equal{#1}{1}}{ \mathcal{E} }{ #1 }
}
}

\newcommand{\sh}{\textsf{Sh}}

\newcommand{\spec}[1][]{\text{Spec}_{#1}\text{ }}

\newcommand{\scal}[2]{ \left\langle #1 , #2 \right\rangle }

\newcommand{\R}{\mathbb{R}}
\newcommand{\C}{\mathbb{C}}
\newcommand{\N}{\mathbb{N}}

\newcommand{\U}{\mathbb{U}}

\newcommand{\Acal}{\mathcal{A}} 
 
\newcommand{\Ecal}{\mathcal{E}} 
 
\newcommand{\Tcal}{\mathcal{T}}

\newcommand{\Scal}{\mathcal{S}} 
\newcommand{\Dcal}{\mathcal{D}} 
\newcommand{\Fcal}{\mathcal{F}} 
 
\newcommand{\Hcal}{\mathcal{H}} 
 
\newcommand{\Kcal}{\mathcal{K}} 
\newcommand{\Lcal}{\mathcal{L}}

\newcommand{\Ccal}{\mathcal{C}} 
 
\newcommand{\Bcal}{\mathcal{B}}

\newcommand{\Ago}{\mathfrak{A}}

\newcommand{\Ygo}{\mathfrak{Y}} 
 
\newcommand{\Igo}{\mathfrak{I}}

\newcommand{\Dgo}{\mathfrak{D}} 
\newcommand{\Fgo}{\mathfrak{F}}

\newcommand{\Cgo}{\mathfrak{C}}

\newcommand{\Yon}{\mathbb{Y}}

\usepackage{titlesec}
\titleformat{\subsubsection}[runin]{\normalfont}{\thesubsubsection}{0pt}{}[.]

\renewcommand{\thesubsubsection}{\arabic{section}.\arabic{subsection}.\arabic{subsubsection}}
\csname @addtoreset\endcsname{subsubsection}{section}

\newcommand{\block}[1]
{

\par \subsubsection{} #1

\bigskip}

\newcommand{\blockn}[1]{\par #1 \bigskip}

\newcommand{\Th}[1]
	{
	\bigskip	
	\textbf{Theorem : }{\itshape #1}
		
	\bigskip
	}

\newcommand{\Prop}[1]
	{

	\bigskip
	
	\textbf{Proposition : }{\itshape #1}
		
	\bigskip
	
	}

\newcommand{\Cor}[1]
	{

	\bigskip
	
	\textbf{Corollary : }{\itshape #1}	
		
	\bigskip

	}

\newcommand{\Lem}[1]
	{

	\bigskip
	
	\textbf{Lemma : }{\itshape #1}
		
	\bigskip
	
	}

\newcommand{\Def}[1]
	{
	
	\bigskip
	
	\textbf{Definition : }{\itshape #1}
	
	\bigskip
	
	}

\newcommand{\Dem}[1]{
	
	\smallskip
	
	\textbf{Proof : } \par
	 {#1} $\square$
	 
	 \bigskip
}

\begin{document}

\pagestyle{plain}
\title{Complete $C^{*}$-categories and a topos theoretic Green-Julg theorem}
\author{Simon Henry}

\maketitle

\begin{abstract}

We investigate what would be a correct definition of categorical completeness for $C^{*}$-categories and propose several variants of such a definition that make the category of Hilbert modules over a $C^{*}$-algebra a free (co)completion. We extend results about generators and comparison theory known for $W^{*}$-categories with direct sums and splitting of symmetric projections to our ``complete $C^{*}$-categories'' and we give an abstract characterization of categories of Hilbert modules over a $C^{*}$-algebra or a $C^{*}$-category as ``complete $C^{*}$-category having enough absolutely compact morphisms''.

We then apply this to study the category of Hilbert spaces over a topos showing that this is an example of a complete $C^{*}$-category. We prove a topos theoretic Green-Julg theorem: The category of Hilbert spaces over a topos which is locally decidable, separated and whose localic reflection is locally compact and completely regular is a category of Hilbert modules over a $C^{*}$-algebras attached to the topos. 

All the results in this paper are proved constructively and hence can be applied themselves internally to a topos. Moreover we give constructive proof of some known classical results about $C^{*}$-algebras and Hilbert modules.

\end{abstract}

\renewcommand{\thefootnote}{\fnsymbol{footnote}} 
\footnotetext{\emph{Keywords.} $C^{*}$-categories, topos, separated topos, complete $C^{*}$-categories}
\footnotetext{\emph{2010 Mathematics Subject Classification.} 18B25, 03G30, 46L05}
\footnotetext{\emph{email:} simon.henry@college-de-france.fr}
\renewcommand{\thefootnote}{\arabic{footnote}} 


\tableofcontents

\section{Introduction}

\blockn{Topos theory and non-commutative geometry are two topics that offer generalizations of ordinary topology which can deal with certain spaces that are too singular from the point of view of ordinary topology. It appears that there is a lot of examples of such ``singular spaces'' that can be studied both from the point of view of topos theory and from the point of view of non-commutative geometry: the space of leaves of a foliation, space of orbits of a dynamical system or of a topological groupoid, the space of ``asymptotic'' paths of a graph or of generalization of graphs, etc.

This paper is part of a project initiated in the author's thesis to try to reconciliate these two points of view on these singular spaces by finding relations between toposes and operator algebras. The introduction of the author's thesis contains a more extended discussion of the motivations for this project that would be out of the scope of the present paper.

Because of the close relation between toposes and topological groupoids (initiated in \cite{joyal1984extension}), and the well known relation between topological groupoids and $C^{*}$-algebras (see for example \cite{renault1980groupoid} ) it is natural to expect that it should be possible to associate a $C^{*}$-algebra to a ``reasonable'' topos, while it should be difficult, if not impossible, to construct interesting toposes out of a $C^{*}$-algebra, unless the $C^{*}$-algebra comes with some additional structures. The $C^{*}$-algebra attached to a topos should then contains a lot of informations about the geometry of the topos, hence such constructions should allow to transport results and techniques from non-commutative geometry to topos theory, and topos theory should produce examples of $C^{*}$-algebras which can be studied geometrically.

Our key technique to attach $C^{*}$-algebras to a topos consist in a very nice notion of continuous fields of Hilbert spaces over a topos (also called Hilbert bundles over a topos). It has been introduced and studied in the early years of categorical logic by Mulvey, Burden, Banachewski and others (see for example \cite{burden1979banach},\cite{mulvey1980banach}). It generalizes the notion of continuous fields of Hilbert spaces over a topological space, and are defined as Hilbert spaces in the internal logic. The category of continuous fields of Hilbert spaces over a topos appears to be a $C^{*}$-category (see section \ref{subsection_HB1} for a precise statement and a proof), providing this way a lots of $C^{*}$-algebras attached to a topos (the algebras of endomorphism of a continuous fields of Hilbert space). Unfortunately these $C^{*}$-algebras are a little too large to be interesting: the $C^{*}$-algebras that we want to consider appears as sub-algebras of these algebras and we need to find a process to select nice sub-algebras in general.
}

\blockn{In previous works (\cite{henry2014measure}, \cite{henry2015toward}) we restricted our attention to boolean toposes for which the $C^{*}$-categories obtained is monotone complete, and in fact often a $W^{*}$-category. In this situation we have been able to obtain some promising results on the relation between boolean toposes and monotone complete $C^{*}$-algebras (for example, a reconstruction theorem in \cite{henry2015toward}). It appears in this previous work that the good point of view is to see the $C^{*}$-category of Hilbert bundle over such a topos as a ``complete\footnote{It appears in the results of \cite{wstarcat} that the correct notion of completeness for $W^{*}$-category is the existence of orthogonal sums and splitting of symmetric idempotent.}'' (in the categorical sense) $W^{*}$-categories. Unfortunately, completeness for $C^{*}$-categories is a subtle notion that is different from completeness for ordinary category and was only well understood for $W^{*}$-categories, hence it seemed like a very natural thing to do to develop a notion of complete $C^{*}$-category that will encompass the examples coming from topos theory, this is a key step if one wants for example to extend the reconstruction theorem of \cite{henry2015toward} to non-boolean situations.

One of the main goal of this paper is to develop such a notion. The resulting formalism allows to deal with the example we had in mind and prove all the results we were interested in, but despite that it does not seem completely satisfying: for example one is not able to prove that for a $C^{*}$-category being complete is an intrinsic property, and not an additional structure that we need to choose, for this reason we will only talk of ``pre-complete'' $C^{*}$-category and we will never define what it means for a $C^{*}$-category to be complete.

We show that we can recover a large part of the properties of ``complete'' $W^{*}$-categories in this framework (see section \ref{subsection_Generators_comparison}), and we give an abstract characterization of categories of Hilbert modules over a $C^{*}$-algebra in terms of existence of enough ``compact'' operators (see section \ref{subsection_compact}).
}

\blockn{We then use this concept of (pre)complete $C^{*}$-categories to deal with the case of a relatively restricted class of topos which are ``almost'' non-singular (namely, the locally decidable, separated with a locally compact localic reflection) and for which it is reasonable to have an analogue of the Green-Julg theorem\footnote{stating that for a proper locally compact groupoid there is a natrual equivalence of categories between Hilbert module over the $C^{*}$-algebra of the groupoid and equivariant continuous fields of Hilbert space over the groupoid.}, and the main result for the paper is indeed a proof of such a theorem, which allow to associate to toposes in this restricted class a well defined $C^{*}$-algebra with a nice universal properties (in the category of $C^{*}$-algebras and bi-modules as morphisms), which by the ordinary Green-Julg theorem is the universal property of the reduced/maximal $C^{*}$-algebra that we want to associate to these examples in the case of proper groupoids. The proof of the theorem is based on the abstract characterization of (complete) $C^{*}$-categories which are category of Hilbert module over a $C^{*}$-algebra (or over a $C^{*}$-category) mentioned above.}

\blockn{One could think that because this class is made of topos that are ``almost non-singular'' it avoids all the interesting phenomenons of non-commutative geometry, and this is not far from being true. The reason why this class is interesting is because all the examples mentioned above are ``locally'' in this restricted class and in a future work we will use this to associate a reduced and a maximal $C^{*}$-algebra to any topos which is only locally of the form which is treated in this paper.}

\blockn{This paper is written in a completely constructive formalism (more precisely within the internal logic of a topos), which allows to have a ``familial'' version of all our results. For example if one has a geometric morphism $f :\Tcal \rightarrow \Ecal$ between two toposes which satisfies a relative version of the hypothesis of our theorem then the category of continuous fields of Hilbert spaces over $\Tcal$ is equivalent to the category of continuous fields of Hilbert modules over a field of $C^{*}$-algebras over $\Ecal$. It appears that this paper is (to our knowledge) the first reference to deal with non-commutatve $C^{*}$-algebra within a topos theoretic constructive framework, for this reason we needed to prove a lot ``classical'' and ``basic'' results about $C^{*}$-algebras and $C^{*}$-categories within constructive mathematics. Those results have been concentrated in the first section and the first appendix. Even if we were not interested in the constructive/relative version of the main theorem we need to have a large part of the material on $C^{*}$-categories constructively valid in order to prove the classical version of the main theorem as the proof involves working with categories of Hilbert modules in the internal logic of the topos.}

\blockn{Hence this paper deal with three inter-connected topics divided in three sections:

\begin{itemize}

\item Extend the constructive theory of $C^{*}$-algebras (and $C^{*}$-categories) outside of the commutative case: the commutative case is well understood by the proof of the Gelfand duality first by C.J.Mulvey and B.Banachewski in Gorthendieck topos in \cite{banaschewski2006globalisation} and then by T.Coquand and B.Spitters \cite{coquand2009constructive} in a fully constructive framework. But some aspect of the non-commutative case where still a little out of reach, for example the notion of spectrum, spectral radius and positivity were not completely obvious, and are developed in this paper. We also re-develop in this framework a large part of the basic theory of $C^{*}$-categories, including Hilbert modules, compact operators, Morita equivalences. This mostly done in sections \ref{section_prelims} and appendix \ref{appendix_positivity}.

\item To develop a notion of ``completeness'' for $C^{*}$-categories. Because $C^{*}$-categories are dagger categories (or $*$-categories) the classical notion of limits and co-limits are not very well suited for them. On the other hand, it has been showed in \cite{wstarcat} that there is a very good notion of completeness for $W^{*}$-categories, that it is equivalent to the existence of arbitrary orthogonal sums and splitting of symmetric projections and that functors are ``limit preserving'' if and only if they are normal. It seems natural that there should also be such a notion of completeness and limit preserving functor for general $C^{*}$-categories which includes both $W^{*}$-categories with normal functors between them, the category of Hilbert modules over a $C^{*}$-algebra and functor coming from tensorization by bi-module, as well as all $C^{*}$-categories coming from geometry or topos theory (category of Hilbert bundle with structure over some ``space'') that can be in between the two other examples. We propose here several definitions that might not be completely satisfying but includes all the examples we had in mind and proves all the properties we wanted to prove. This is done in section \ref{section_completeness}.

\item To prove a topos theoretic Green-Julg theorem: for a reasonable \emph{separated} toposes $\Tcal$, there should be a $C^{*}$-algebra attached (which is both the reduced and the maximal $C^{*}$-algebra) such that continuous fields of Hilbert spaces over $\Tcal$ are the same as Hilbert $\Ccal$-modules. And this last result should be proved constructively. This is done in section \ref{section_Hilbert_bundle}.

\end{itemize}

}

\blockn{Those are interconnected in the following way: section \ref{section_Hilbert_bundle} is entirely based on section \ref{section_completeness} and need it to be constructive. Section \ref{section_completeness} needs section \ref{section_prelims} and appendix \ref{appendix_positivity}, but only if we are interested in the constructive aspect, it can be considered as self-contained if we are only interested in the classical results. Section \ref{section_prelims} is mostly self-contained (it use at some very basic result of classical analysis whose constructive proof is in the appendix) while appendix \ref{appendix_positivity} require a bit of familiarity with locale theory and relies to some extent on the theory of localic $C^{*}$-algebra developed by the author in \cite{henry2014localic}.

Moreover, sections \ref{section_prelims}, \ref{section_completeness} and appendix \ref{appendix_positivity} can be read without any knowledge or familiarity with topos theory, and should be very accessible\footnote{A reader unfamiliar with constructive mathematics might found some part of the paper more complicated than necessary , but should be able to ignore those aspect.}, while the final section \ref{section_Hilbert_bundle} require a lot of familiarity with topos theory, and in particular with their internal logic and the interplay between internal and external logic. }

\subsection{General preliminaries}

\blockn{The general framework we are working in is the internal logic of a Heyting pre-topos with a family of small maps as in \cite{joyal1995algebraic} which satisfies all the additional axioms presented in \cite{joyal1995algebraic}. Appendix \ref{appendiw_foundation} at the end of the paper contains the details about this, including the precise list of axioms. Roughly, this framework is very similar to a ``topos theoretic'' version of Neuman-Bernays-Godel set theory: we have a notion of class and a notion of sets and a class of all sets. We can talk about equality of set (because they are element of a class) although it is not really interesting. A class if said to be small if it is isomorphic to a set (strictly speaking, and contrary to NBG set theory, it would not make sense to say that a small class ``is'' a set, but we might use this as an abuse of language occasionally). The category of sets form an elementary topos with a natural number object which we call the ``base topos'', but the category of all class is just a Heyting pre-topos (so no class of sub-class of a given class, or class of function between two classes). 

Any model of Neuman-Bernays-Godel set theory, or any model of Zermelo-Frankel set theory with either a Grothendieck universe or an inaccessible cardinal gives us such Heyting pre-topos with a class of small map, any Grothendieck topos over such a base universe also give us such a category (by looking at ``class valued'' sheaves). }

\blockn{In this paper, categories are always assumed to be ``locally small'' in the sense that the morphisms between two object form a sets, but we do not assume the class of objects to be a set or small. A small category is category whose class of objects is a set.
By ``functors'' we actually mean ``anafunctors'', or more precisely ``saturated anafunctor'' (hence in particular locally small anafunctors as all our categories are locally small) in the sense of \cite{makkai1996avoiding}. Informally an anafunctor is the same thing as a functor with the exception that the value of the functor on an object is only well defined ``up to unique isomorphism''. We refer the reader to \cite{makkai1996avoiding} for more details about this notion, their usage will be completely transparent. }

\blockn{When we talk about a topos, we generally mean a topos which is a Grothendieck topos over our base topos, i.e. a category equivalent to a category of sheaf over a small site, it can be extended into a Heyting pre-topos with a class of small map satisfying the axiom given in appendix \ref{appendiw_foundation} by looking at class valued sheaves over the site (this is sketched in appendix \ref{appendiw_foundation}). }

\section{Preliminaries on $C^{*}$-categories}
\label{section_prelims}

\blockn{This section contains basic preliminaries on $C^{*}$-algebras, $C^{*}$-categories and Hilbert modules over them. The only new contributions are to formulate well known results of the theory of $C^{*}$-algebras in the language of $C^{*}$-categories and/or to prove them constructively. We included some of the proofs just in order to show the constructiveness of the results.

A reader already familiar with the theory of $C^{*}$-algebras, $C^{*}$-categories and Hilbert modules, and either not interested in the constructive aspect or willing to thrust the author on the constructiveness of these results can skip this section almost entirely. There is only two concepts that can be considered as new to some extent and are worth looking at: the notion of small Hilbert modules over a (large) $C^{*}$-categories between \ref{genHilbmod} and the end of subsection \ref{subsection_CcatHM} and to some extend the subsection \ref{Restricton_of_Hilbert_modules} about restriction of Hilbert modules. }

\subsection{General analytic preliminaries}
\label{subsection_genAnaprelim}
\blockn{This subsection is only here to make a bit more precise how we do analysis in constructive mathematics: what do we call a real number, a metric space and how completeness is defined.}

\blockn{One of the main subtleties of constructive analysis is the definition of real numbers: all the classical definitions of real numbers tends to become non-equivalent without the axiom of choice and the law of excluded middle and all gives different objects of real numbers with very different properties. We will only recall the property of the sets of real numbers we are interested in, without proof or to much details. The details and proof can be found in Section D4.7 of \cite{sketches} which is a very good places to find several such definitions with their relations and properties from a topos theoretic perspective.
In the topos theoretic framework the more natural choice for being ``the object of real number'' when we do analysis is the object of Dedekind real numbers, also called the object continuous real numbers. A Dedekind (or continuous) real number is a two sided Dedekind cut $(L,U)$ of rational number, with the relation between $L$ and $U$ being expressed as $\forall l \in L, \forall u \in U, l <u$ and for all $q < q'$ either $q \in L$ or $q' \in U$. Dedekind real numbers have relatively good algebraic properties, they of course do not form a field in the geometric sense as equality to zero is non-decidable, but they form a local ring. They have relatively well behaved order relation $<$ and $\leqslant$ and one can show that $\leqslant$ is the negation of $>$. Moreover, a continuous real number can be approximated arbitrary well by rational numbers in the sense that for every integer $n$, one has $\exists q, q \leqslant x \leqslant q+1/n$, and this property is characteristic of continuous real number in the sense that it implies the defining property that for all $q < q'$ either $q <x $ or $x < q'$. The choice of Dedekind real numbers make the topos theoretic version of constructive analysis rather different from older approaches, like for example the Brouwer-Bishop constructive analysis, which are in general based on sequential definitions of real numbers.

\bigskip

The object of continous real number will be simply denoted by $\mathbb{R}$, and one defines the complex number $\mathbb{C}$ as $\mathbb{R} \times \mathbb{R} = \mathbb{R}[X]/(X^{2} + 1)$, it is not exactly true in constructive mathematics that complex number are algebraically closed but this is mostly because they are not a field, one can for example prove\footnote{We will not do it and we will not use this result, but at least in Grothendieck toposes it follows from Barr's covering theorem.} in our framework that every separable polynomial can be splited. There is at least three reasons why this notion is the correct one: the first is that it is a geometric notion (one has a classifying space for the theory of continuous complex numbers) hence one can expect that at least because of Barr's covering theorem all classical ``geometric'' results can be extended to them, the second is that their classifying space is exactly the space of complex numbers\footnote{Only assuming the law of excluded middle, in general it is the formal locale of complex number, see below.}, hence sections of the object of continuous complex numbers corresponds exactly to continuous complex valued functions, the last reason is that it is for this notion that Mulvey's comparison result with the notion of continuous field of Banach space works (see \cite{mulvey1980banach}).}

\block{The main defect of this object of continuous real numbers is the absence of supremum and infimum of inhabited bounded subset (the existence of such supremum is equivalent to the law of excluded middle) which makes them not well suited to be the target of norm and distances function which are generally defined as supremums or infimums. For this reason we will use a different set of real numbers as natural target for our distances and norms, (following previous works on the subject as \cite{burden1979banach} and \cite{mulvey1980banach}). Norms and distances will take value in the set of positive (bounded) upper semi-continuous real numbers, they are upper Dedekind cut: so sub-sets $U \subset \mathbb{Q}$ such that $q \in U \Leftrightarrow \exists q' <q, q' \in U$, $q \in U \Rightarrow q>0$ and $\exists q \in U$, the last condition corresponds to the ``bounded'' property, and could be remove when we talk about distance (but is preferable for norms). Positive continuous real numbers can be embedded into postive upper semi-continuous real numbers by $(L,U) \mapsto U$, which is an injection. Semi-continuous real numbers can be added and multiplied but even if we removed the positivity condition, there is no way to substract them or to extend the multiplication to negative elements: indeed, section of the object of positive bounded upper-continuous real numbers corresponds to upper-continuous functions with values in positive real numbers which are locally bounded, hence the opposite of such a function is no longer a upper semi-continuous function, but rather a lower semi-continuous function. }

\block{A distance on a set $X$ is just a function from $X \times X$ to the set of positive upper semi-continuous real numbers that satisfies the usual three axioms of metric set. A set endowed with a distance is said to be complete if every Cauchy filter is convergent. A constructive account of the notion of Cauchy filter, completeness and completion can be found in \cite{henry2014localic}, it is written in the more general framework of locales endowed with a (symmetric) distance functions but it can be specialized to sets with a distance functions. The short introduction that follows should suffice for the present paper.}

\block{We just recall basic fact about Cauchy filters: In constructive mathematics Cauchy filters are defined such that for every element $f$ of the filter $\exists x \in f$, and not just $f \neq \emptyset$, aside from that the theory is similar to the classical notion of Cauchy filter (as in \cite{bourbaki1966elements} ). Every Cauchy filter contains a unique regular/minimal Cauchy filter, i.e. a Cauchy filter generated by subset which are already $\epsilon$ thickening of an other inhabited subset. In particular two Cauchy filters are contained in a same filter if and only if their intersection is a Cauchy filter if and only they contains the same regular Cauchy filter. Hence a Cauchy filter is convergent if and only the regular Cauchy filter it contains is the filter of neighbourhood of a point and the completion can be constructed as the set of regular Cauchy filter. In \cite{henry2014localic} we only defined the localic completion, which is the classifying space for the theory of regular Cauchy filter, but the ordinary completion is the set of point of the localic completion.}

\blockn{In constructive mathematics, it is helpful to have an alternative to Cauchy filter where the convergence speed is controlled, this is especially useful when we try to relate internal results to external results, typically in lemma \ref{Lem_sectionarecomplete}. This is achieved by the following notion:}

\block{\Def{Let $X$ be a set with a distance $d$. A Cauchy approximation on $X$ is a collection of subset $A_n \subset X$ such that:

\begin{itemize}

\item $\forall n \ \exists x \in A_n,$

\item $\forall n \ A_{n+1} \subset A_n,$

\item $\displaystyle \forall n \ \forall x,y \in A_n \ d(x,y)<\frac{1}{n}.$

\end{itemize}

A Cauchy approximation is said to converge to $x \in X$ if and only if:

\[ \forall n \ \forall y \in A_n \ d(x,y) \leqslant \frac{1}{n} \]

}
This relates to Cauchy filter in the following way: each Cauchy approximation generates a Cauchy filter $\Fcal = \{ F \subset X | \exists n \ A_n \subset F \}$. Every minimal Cauchy filter $\Fcal$ is generated by the Cauchy approximation defined by: 

\[ A_n = \bigcup_{f \in \Fcal \atop \delta(f) < \frac{1}{2n}} f. \]

And finally, a Cauchy approximation converge to $x$ if and only if the Cauchy filter it generates converge to $x$. All these fact are easy, $\epsilon,\delta$ manipulation. and this implies that:

\Cor{A metric set $X$ is complete if and only if every Cauchy approximation in $X$ converge. }

In \cite{burden1979banach}, \cite{mulvey1980banach} or other references by the same authors, a different notion of Cauchy approximation is used where there is no control on the convergence speed like for Cauchy filter. Our notion is more restrictive, the two notions of convergence can be shown to be equivalent for our definition of Cauchy approximation, and the two notion define the same notion of completeness.

}

\block{Following \cite{mulvey1980banach}, a Banach space is a $\C$-module (where $\C$ is the object of ``continuous complex numbers'') endowed with a norm function which takes value in the set of upper semi-continuous real number, such that all the usual algebraic axiom of Banach space are satisfied and which is complete in the sense of Cauchy filter (or Cauchy approximation as above). By a ``Banach norm'' on a vector space, we mean a norm that makes it into a Banach space.}

\block{Banach spaces cannot be directly pulled back along geometric morphism, if $B$ is a banach space then $f^{*}(B)$ comes with a norm that can fail to be satisfies the axiom of separation and completeness, and it has an action of $f^{*}(\mathbb{C})$ that can be different from the $\mathbb{C}$ of the source topos of $f$, but if we quotient by the subspace of norm zero elements and take the completion one obtain $f^{\sharp}(B)$ which is a Banach space. This is this operation that we call the pullback of a Banach space, it also corresponds to the pullback of localic Banach spaces as defined in \cite{henry2014localic}.
}

\subsection{$C^{*}$-algebras, $C^{*}$-categories and Hilbert modules}
\label{subsection_CcatHM}

\block{\Def{A $C^{*}$-algebra is a $\C$-algebra $\Ccal$ endowed with a Banach norm $\Vert \_ \Vert$ and an isometric anti-linear involution $* : \Ccal \rightarrow \Ccal$ such that:

\begin{itemize}

\item $\Vert a b \Vert \leqslant \Vert a \Vert \Vert b \Vert$.

\item $(ab)^{*} = b^{*}a^{*}$.

\item $\Vert x \Vert^{2} = \Vert x^{*} x \Vert$.

\end{itemize}

}
}

\blockn{A lot of classical results about $C^{*}$-algebras can also be proved in constructive mathematics. One of the biggest achievement in this direction is the proof\footnote{see \cite{banaschewski2006globalisation} for the proof in Grothendieck toposes and \cite{coquand2009constructive} for a fully constructive proof, and \cite{henry2014nonunital} for the case of non-unital algebras} of a Gelfand duality between commutative $C^{*}$-algebra and (locally) compact completely regular locales. A locales is a variant on the idea of topological spaces which does not necessary have points (see \cite{picado2012frames}) assuming the axiom of choice any locally compact locales has enough point and is the same as a locally compact topological space, so in classical mathematics this result recovers the ordinary Gelfand duality. The fact that we need locales instead of topological spaces in constructive mathematics corresponds to the fact that without the Hahn-Banach theorem it is not possible to construct character of an arbitrary commutative $C^{*}$-algebra and those are exactly the points of the spectrum. This has one consequence for the theory of $C^{*}$-algebras: in order to use continuous functional calculus one need a function defined on the ``formal locale of complex numbers'' and not just on the topological space of complex numbers (which might be different if we do not assume the law of excluded middle), but this distinction is not going to be a problem: a function on this formal locale is exactly a function on the topological space of complex numbers which is uniformly continuous on every bounded subset, so all classical functions can be used in functional calculous, indeed, by the theory of \cite{henry2014localic} any uniformly continuous function on the rational numbers can be extended to their localic completion which is the formal locale of real numbers, and conversely any function defined on the formal locale is uniformly continuous on bounded subset by compactness of closed bounded interval.}

\blockn{The theory of general $C^{*}$-algebras works relatively well in constructive mathematics. There is of course some difference with the classical theory: for example the existence of a faithful representation of an arbitrary $C^{*}$-algebra on a Hilbert space which is a consequence of the Hahn-Banach theorem in classical mathematics and only holds after pull-back along a proper geometric morphism in constructive mathematics.}

\blockn{Appendix \ref{appendix_positivity} contains statement and constructive proofs of some basic results about $C^{*}$-algebras that are well known in classical mathematics. We will use those results freely here. The most important is proposition \ref{Prop_def_positivity} about the notion of positivity in non-commutative $C^{*}$-algebras.
}

\block{\Def{\begin{itemize} 
\item A Banach category is a category $\Cgo$ such that for each pair of objects $A,B \in |\Cgo|$ the set of morphisms $\Cgo(A,B)$ is endowed with a Banach space structure such that the composition is a bilinear map of norm smaller than one i.e. $\Vert f \circ g \Vert \leqslant \Vert f \Vert \Vert g \Vert$.

\item A $*$-category\footnote{The terminology ``Dagger-category'' is more common for these categories, but in the context of $C^{*}$-category this terminology seems more appropriate, by analogy with $*$-algebras.} is a category $\Cgo$ endowed with an involutive contravariant functor $* : \Cgo \rightarrow \Cgo^{op}$ which is the identity on object. A Banach $*$-category is a category with both a Banach category structure and a $*$-category structure such that $*$ is isometric and anti-linear.
\end{itemize}
}
}

\block{\label{Def_Cstarcar}\Def{A $C^{*}$-category is a Banach $*$-category $\Cgo$ in which the following holds:

\begin{itemize}

\item For any arrow $h$ in $\Cgo$ one has $\Vert h \Vert^{2} = \Vert h^{*}h \Vert $.

\item For all $h : X \rightarrow Y \in \Cgo$, $h^{*}h \in \Cgo(X)^{+}$.

\end{itemize}
}

The second condition cannot be removed: Consider the Banach $*$-category $\Fgo$ with two objects $a,b$ such that $\Fgo(a)=\Fgo(b)=\Fgo(a,b)=\Fgo(b,a)=\C$ as Banach space, composition is complex multiplication and the $*$ operation is complex conjugation on endomorphisms and $z \mapsto - \overline{z}$ on morphisms between $a$ and $b$. This satisfies all the conditions for being a $C^{*}$-category except the positivity condition.}

\block{In the case where the category has bi-product (i.e. binary co-product, which are also binary product as the category is enriched over abelian group), and that they are ``orthogonal'' in the sense that the inclusion map of $A$ and $B$ into $A \oplus B$ are adjoint to the projection maps then it is enough to check that for each object $A$ the algebra of endomorphism is a $C^{*}$-algebra, indeed the conditions on morphism $f :X \rightarrow Y$ can be recovered by seeing them as endomorphisms of $X \oplus Y$.}

\block{\label{Prop_alternatedDefCstarcat}One can give a more concise and often easier to check definition of $C^{*}$-categories: 
\Prop{A Banach $*$-category $\Cgo$  is a $C^{*}$-category if and only if it satisfies the following strong $C^{*}$-inequality: If $A,B,C$ are three objects of $\Cgo$, $f \in \Ccal(A,B), g \in \Cgo(A,C) $ then:

\[ \Vert f \Vert^{2} \leqslant \Vert f^{*}f + g^{*}g \Vert \]

}

In fact, it appears that putting the strong $C^{*}$-inequality in the definition of $C^{*}$-algebra makes the basic theory of $C^{*}$-algebra considerably simpler (it makes all the material presented in appendix \ref{appendix_positivity} either useless or trivial). In some sense the fact that the ordinary $C^{*}$-equality implies the strong $C^{*}$-inequality for complex Banach algebra can be considered as an accident (and the main result of appendix \ref{appendix_positivity}). It is for example not the case for real Banach algebras either: One also need to consider the strong $C^{*}$-inequality in order to obtain a correct notion of real $C^{*}$-algebra.

\Dem{Any $C^{*}$-category satisfies this strong $C^{*}$-inequality: $\Vert f \Vert^{2} = \Vert f ^{*} f \Vert$ by the $C^{*}$-equality and $\Vert f^{*}f \Vert \leqslant \Vert f^{*}f +g^{*}g \Vert$ because $g^{*}g$ is positive. Conversely, if we assume the strong $C^{*}$-inequality, then taking $g=0$ gives $\Vert f \Vert ^{2} \leqslant \Vert f ^{*} f \Vert$, and $\Vert f^{*}f \Vert \leqslant \Vert f \Vert ^{2}$ follows from the axioms of Banach categories and the fact that $*$ is isometric. Finally the positivity axioms follow from the strong inequality because if $f \in \Cgo(X,Y)$ is any arrow and $K$ is bigger than $\Vert f \Vert ^{2}$ then $K - f^{*}f$ is a positive element of $\Cgo(X)$, hence there exists a $g \in \Cgo(X)$ such that $g^{*}g = K- f^{*}f$ and one has by the strong $C^{*}$-inequality:

\[ \Vert K- f^{*}f \Vert = \Vert g^{*}g \Vert \leqslant \Vert g^{*}g +f^{*}f \Vert = \Vert K \Vert =K \] 

Hence $f^{*}f$ is positive.
}

}

\block{We will also talk about \emph{non-unital} $C^{*}$-categories, i.e. structure which satisfy all the properties of a $C^{*}$-category but for which we do not require the existence of ``identity" $1_A \in \Cgo(A)$ in the underlying category. In fact, unless explicitly stated otherwise \emph{(small) $C^{*}$-categories will not be assumed to be unital}. This being said, the reader should also note that the definition of pre-complete $C^{*}$-categories given in definition \ref{Defpre-comp} assume that they are unital.}

\block{The first example of $C^{*}$-category is the following: If $\Ccal$ is a $C^{*}$-algebra then $\Ccal$ can be seen as a $C^{*}$-category with only one object and $\Ccal$ as its algebra of endomorphisms. The second example will be the category of Hilbert modules over a $C^{*}$-category defined below.}

\block{\Def{Let $\Cgo$ be a $C^{*}$-category (possibly non-unital), A (right) pre-Hilbert module $H$ over $\Cgo$, or a (right) pre-Hilbert $\Cgo$-module $H$, is the data of:

\begin{itemize}

\item A contravariant $\C$-linear functor $H$ from $\Cgo$ to the category of $\C$-vector spaces. If $a \in H(A)$ and $f:B \rightarrow A$ is an arrow in $\Cgo$ we denote by $af$ or $a \circ f$ the corresponding element of $H(B)$.

\item For all pair of objects $A,B \in \Cgo$, for each $a \in H(A)$ and $b \in H(B)$ one has an element:

\[ \scal{b}{a} \in \Cgo(A,B) \]

\end{itemize}

Which satisfies the following axioms, for $a \in H(A)$, $b\in H(B)$ and $f:C\rightarrow A$ in $\Cgo$:

\begin{itemize}

\item $\scal{b}{a} = \scal{a}{b}^{*}$.

\item $(x,y) \mapsto \scal{x}{y}$ is linear in the second variable (and hence anti-linear in the first).

\item $\scal{b}{a \circ f}=\scal{a}{b} \circ f$.

\item $\scal{a}{a}$ is a positive (self-adjoint) element of $\Cgo(A,A)$.

\end{itemize}

}

One can of course also defines a notion of \emph{left} pre-Hilbert $\Cgo$-modules (for which the scalar product would be taken linear in the first variable and anti-linear in the second) but they will play almost no role in the present paper and we will says (pre-)Hilbert module instead of right (pre-)Hilbert module.
}

\block{\Prop{Let $H$ be a pre-Hilbert $\Cgo$-module, then:
\begin{itemize}

\item One has a ``Cauchy-Schwarz'' inequality: for of any $a \in H(A), b \in H(B)$ and $K$ a real number such that $\Vert \scal{a}{a} \Vert \leqslant K$ one has\footnote{One cannot state the inequality as $ \scal{a}{b}^{*}\scal{a}{b} \leqslant \Vert \scal{a}{a} \Vert \scal{b}{b} $, because $\Vert \scal{a}{a} \Vert$ is in general only a semi-continuous number and hence it does not make sense to multiply an element of the algebra by such a number, this is why we need to introduce this number $K$.}:

\[ \scal{a}{b}^{*}\scal{a}{b} \leqslant K \scal{b}{b} \]

In particular:

\[ \Vert \scal{a}{b} \Vert \leqslant \Vert \scal{a}{a} \Vert^{\frac{1}{2}} \Vert \scal{b}{b} \Vert^{\frac{1}{2}} \]

\item For each object $A$ of $\Cgo$, \[ \Vert a \Vert := \Vert \scal{a}{a} \Vert ^{\frac{1}{2}}\] defines a semi-norm on $H(A)$.

\end{itemize}

}

The proof given in \cite[Prop. 1.1]{lance1995hilbert} can easily be adapted from $C^{*}$-algebras to $C^{*}$-categories and made constructive.
}

\block{\Def{A Hilbert $\Cgo$-module is a pre-Hilbert $\Cgo$-module such that for each object $A \in |\Cgo|$, $H(A)$ is complete and separated for the norm defined in the previous proposition.}

As always, completeness is taken in the sense of Cauchy filters or equivalently Cauchy approximations (see section \ref{subsection_genAnaprelim}). Any pre-Hilbert module can be completed into a Hilbert module by taking the separated completion of $H(A)$ for each object $A$ of $\Cgo$.

}

\block{Here is the key example of Hilbert module: let $\Cgo$ be any $C^{*}$-category, $A \in |\Cgo|$ any object, then there is a Hilbert $\Cgo$-module $\Yon_A$, or $\Yon(A)$, defined by: 

\[\Yon_A(X):=\Cgo(X,A)\]
\[ \scal{u}{v}:=u^{*}v \]
}

\block{If $\Ccal$ is a $C^{*}$-algebra a Hilbert $\Ccal$-module is a Hilbert module over the $C^{*}$-category with one object corresponding to $\Ccal$. If $I$ is a right (closed) ideal of $\Ccal$ then $I$ is a Hilbert $\Ccal$-module for the product $\scal{a}{b} = a^{*}b$. If $\Cgo$ is a $C^{*}$-category, $A$ and $B$ are two objects of $\Cgo$ then $\Cgo(A,B)$ is a right Hilbert $\Cgo(A)$-module for the scalar product $\scal{f}{g}=f^{*}g$, and also a left Hilbert $\Cgo(B)$-module for the scalar product $\scal{g}{f}=gf^{*}$.}

\blockn{For example, in the case where $\Ccal$ is just the algebra $\C$ of (continuous) complex numbers then a Hilbert $\C$-modules is just a Hilbert space.}

\block{\Def{If $H$ and $G$ are two Hilbert $\Cgo$-modules, a bounded linear map from $H$ to $G$ is a $\C$-linear natural transformation from $H$ to $G$ which is bounded in the sense that:

 \[ \Vert f \Vert = \sup_{A \in |\Cgo|} \Vert f_A \Vert < \infty \]
 
An \emph{operator} from $H$ to $G$ is a pair of bounded linear maps $f : H \rightarrow G$, $f^{*}:G \rightarrow H$ such that for all $h \in H(A)$, $g \in G(B)$:

\[ \scal{h}{f^{*}(g)} = \scal{f(h)}{g} \]
}

$f^{*}$ is called the adjoint of $f$ and is entirely determined by $f$, but its existence is not automatic when $f$ is an arbitrary bounded linear map.

In the case of Hilbert space (i.e. $\Cgo$ is just $\C$), assuming classical logic any bounded linear map is an operator, but this can fail in intuitionist logic or over an arbitrary $C^{*}$-algebras or $C^{*}$-category.
}

\block{\Lem{For any (small) $C^{*}$-category, the space of all bounded linear maps between two Hilbert modules is complete. The space of operators is closed in the space of bounded linear map. In particular the space of operators is also complete.}

Note that smallness is not used explicitly in the proof , it is only here to ensure that there is only a set of operators between two Hilbert modules.

\Dem{Let $A_n$ be a Cauchy approximation in the space of bounded linear maps from $E$ to $F$ two Hilbert $\Cgo$-modules. For each $c \in |\Cgo|$ and $e \in E_c$ the set $(A_n)(c) =\{a(c), a \in A_n \}$ is easily checked to be a Cauchy approximation on $F_c$, hence it converge to a unique element denoted $v(c)$ in $F_c$, one easily check that $c \mapsto v(c)$ is a linear map and is bounded (for example by $1+B$ where $B$ is a bound for an element in $A_1$) and that is is the limit of the Cauchy approximation $A_n$. 

Let now $f$ be a bounded linear map between two Hilbert $\Cgo$-modules $f :E \rightarrow F$ such that for every $\epsilon >0$ there is an operator $g : E \rightarrow F$ such that $\Vert f - g \Vert < \epsilon$, then the sets $A_n = \{ g : F \rightarrow E | \Vert f- g^{*} \Vert < 1/n \}$ form a Cauchy approximation hence it converges to a bounded linear map $g: F \rightarrow E$ which is easily seen to be an adjoint for $f$ hence $f$ is an operator.
}

}

\blockn{\Prop{Let $\Cgo$ be a small $C^{*}$-category. Then, endowed with the above defined norm and the adjunction $ f \mapsto f^{*}$, the category of Hilbert $\Cgo$-modules and operators between them is a $C^{*}$-category denoted $\Hcal\Cgo$.}

If $\Cgo$ is non-small, then the only things that goes wrong in this proposition is that operators between two arbitrary Hilbert modules might not be a set (they can form a proper class). This problem will be dealt with in \ref{DefnonsmallHC}.

\Dem{The category of Hilbert modules has bi-product, so it is enough to check that for any Hilbert $\Cgo$-module $H$ the algebra of endomorphism if a $C^{*}$-algebra. One clearly has $\Vert f g \Vert \leqslant \Vert f \Vert \Vert g \Vert$ as for any $x \in H_A$:

 \[ \Vert f(g(a)) \Vert \leqslant \Vert f \Vert \Vert g(a) \Vert \leqslant \Vert f \Vert \Vert g \Vert \Vert a \Vert,\]
 
one also clearly have $(fg)^{*} = g^{*} f^{*}$ and completeness follow from the above lemma. For the $C^{*}$-equality:

\[ \Vert f(x) \Vert^{2} \leqslant \Vert \scal{f(x)}{f(x)} \Vert = \Vert \scal{f^{*} f(x)}{x} \Vert \leqslant \Vert f^{*}f \Vert \Vert x \Vert^{2}   \] 

Hence $\Vert f \Vert ^{2} \leqslant \Vert f^{*} f \Vert$, while in the other direction $\Vert f^{*} f \Vert \leqslant \Vert f^{*} \Vert \Vert f \Vert = \Vert f \Vert ^{2}$.

}

}

\block{For example, if $\Cgo$ is a unital $C^{*}$-category and $A,B \in \Cgo$ are two objects then the operators between $\Yon_A$ and $\Yon_B$ correspond (functorially) to morphism in $\Cgo$ between $A$ and $B$ (this is just the Yoneda lemma). More generally, operators from $\Yon_A$ to an arbitrary Hilbert $\Cgo$-modules $H$ are in natural bijection with elements of $H(A)$.

If $\Cgo$ is non-unital then operators between $\Yon_A$ and $\Yon_B$ are exactly the morphisms between $A$ and $B$ in the category of multipliers of $\Cgo$ as defined in \cite{vasselli2007bundles}. There is still one operator from $\Yon_A$ to $H$ for each element of $H(A)$, but there is in general more such operators (operators of this kind are sometimes called the multipliers of $H$).}

\block{\Prop{Let $A_1,\dots, A_n$ be a finite sequence of objects of a $C^{*}$-category $\Cgo$. Then there exists a unique structure of $C^{*}$-algebra on the set $M_{\Cgo}(A_1,\dots,A_n)$ or simply $M(A_1,\dots,A_n)$ whose element are matrix:
\[ \begin{pmatrix} 
f_{1,1} & f_{1,2} & \dots & f_{1,n} \\
f_{2,1} & f_{2,2} & \dots & f_{2,n} \\
\vdots & \vdots & \ddots & \vdots  \\
f_{n,1} & f_{n,2} & \dots & f_{n,n} \\
\end{pmatrix}
\]

with $f_{i,j} \in \Cgo(A_j,A_i)$ such that the addition and multiplication are the usual matrix operation, the $*$ operation is $(f^{*})_{i,j} = (f_{j,i})^{*}$, and the topology is the ordinary product topology.
}

\Dem{For the existence, simply take $M(A_1,\dots , A_n)$ to be the algebra of endomorphism of the Hilbert $\Cgo$-module $\Yon(A_1) \oplus  \dots \oplus \Yon(A_n)$ (if $\Cgo$ is a non-unital $C^{*}$-category then take the sub-algebra of those endomorphisms whose matrix elements are in the image of the Yoneda embeddings instead).

Classically, the uniqueness follows from a well known fact that the norm of a $C^{*}$-algebra is entirely determined by its algebraic structure. Indeed, the norm can be computed from the spectral radius (proposition \ref{Prop_specctral_radius}) and the spectrum only depends on the algebraic structures. Constructively, the spectrum does not only depends on the algebraic structures: it also depends on the localic completion (or to put it an other way, one need to know the algebraic structure and the algebraic structures of all the pullback along geometric morphisms). This is why we need this additional assumption that the topology is the usual product topology: same topology and same additive structure implies same uniform structure hence same localic completion and this concludes the proof.

}
}

\block{\label{Pos_op_on_module}\Lem{Let $\Cgo$ be a $C^{*}$-category, $H$ be a Hilbert $\Cgo$-module and $F$ a self-adjoint endomorphism of $H$. Then:

\begin{itemize}

\item If $F$ is positive then $\Vert F \Vert \leqslant K$ if and only if for each $X \in |\Cgo|$ for each $h \in H(X)$ one has $ \Vert \scal{h}{F h } \Vert \leqslant K \Vert \scal{h}{h} \Vert$ in $\Cgo(X)$.

\item  $F$ is positive if and only if for each $X \in |\Cgo|$ for each $h \in H(X)$, $\scal{h}{Fh}$ is a positive element of $\Cgo(X)$.

\end{itemize}

}

\Dem{As $F$ is positive one has $F=S^{2}$ with $S$ self-adjoint, hence $\scal{h}{Fh} = \Vert S h \Vert^{2}$. Hence the inequality $\forall h, \Vert \scal{h}{Fh} \Vert \leqslant K \Vert \scal{h}{h} \Vert$ is equivalent to $\Vert S \Vert ^{2} \leqslant K$ as one has $\Vert S \Vert ^{2} = \Vert S^{*}S\Vert = \Vert F \Vert $ concludes the proof.

For the second point: If $F$ is positive then $\scal{h}{F(h)}$ can be written as $h^{*} F h$ in the $C^{*}$-category of Hilbert modules and hence is positive. Conversely, assume that for every $h$ $\scal{h}{F h}$ is positive. Fix $K$ a real number bigger than the norm of $F$. Then $F \leqslant K$,  and hence, as an application of the first implication to the positive operator $K-F$, for any $h \in H(X)$, $\scal{h}{F h} \leqslant K \scal{h}{h}$, one hence has:

 \[ 0 \leqslant \scal{h}{F(h)} \leqslant K \scal{h}{h},\]
 which gives:
\[ 0 \leqslant \scal{h}{K h - F(h)} \leqslant K \scal{h}{h} \]

As $K-F$ is positive, the first point implies that $\Vert K-F \Vert \leqslant K$ and hence that $F$ is positive.
}

}

\block{\label{pos_matrix}\Cor{A self-adjoint matrix $ M \in M(A_1,\dots,A_n)$ is positive if and only if for all $X \in |\Cgo|$, for all $n$-tuple $t =(t_i)_{i=1\dots n}$ with $t_i \in \Cgo(X,A_i)$ seen as a column vector one has:

\[t^{*} M t = \sum_{i,j} t^{*}_i M_{i,j} t_j \]

is a positive element of $\Cgo(X)$.
}

Note that by applying this lemma to the full subcategory of $\Cgo$ on the object $A_1,\dots, A_n$ one can restrict to the case where $X$ is one of the $A_i$.
\Dem{This follows directly from lemma \ref{Pos_op_on_module} and the construction of the algebra of matrices as a sub-algebra of endomorphisms of $\Yon(A_1) \oplus \dots \oplus \Yon(A_n)$.}
}

\block{\label{genHilbmod}Let $\Cgo$ be a $C^{*}$-category (non necessarily small) and let $X$ be a set together with:

\begin{itemize}
\item for each $x \in X$, an object $s(x) \in |\Cgo|$
\item for each pair $x,y \in X$, an arrow $\scal{x}{y} \in \Cgo(s(y),s(x))$
\end{itemize}

One then has:

\Prop{The following conditions are equivalent:

\begin{itemize}
\item For any finite family $(x_1,\dots, x_n)$ of elements of $X$, the matrix whose coefficients are given by:

\[ A_{i,j} = \scal{x_i}{x_j} \in \Cgo(s(x_j),s(x_i)) \]

is a positive element of the $C^{*}$-algebra $M_{\Cgo}(s(x_1),\dots,s(x_n))$.

\item There exists a Hilbert module $T$, and for each $x\in X$ an element $t_x \in T(s(x))$ such that $\scal{x}{y}=\scal{t_x}{t_y}$.

\item There exists a unique (up to unique unitary isomorphism) Hilbert $\Cgo$-module $\langle X \rangle$ with elements $e_x \in \langle X \rangle(s(x))$ such that $\scal{e_x}{e_y} = \scal{x}{y}$ for each $x,y \in X$ and such that for each $A \in |\Cgo|$ the finite sum of elements of the form $e_xf$ for $f : A \rightarrow s(x)$ are dense in $\langle X \rangle(A)$.

\end{itemize}

}

A set $X$ with such map $s$ and $\scal{ \_ }{ \_ }$ satisfying the conditions of the proposition will be called a set of generators for a Hilbert $\Cgo$-module.

\Dem{The key observation is that by corollary \ref{pos_matrix}, the first condition is equivalent to the fact that any formal linear combination of the form $c = \sum_i x_i \circ f_i $ for $x_i \in X$ and $f_i : B_0 \rightarrow s(x_i)$ will gives:

 \[ \scal{c}{c} = \sum_{i,j} f_i^{*}\scal{x_i}{x_j} f_j \geqslant 0 \]

Hence the second condition implies the first because such a $c$ (replacing the $x_i$ by $t_{x_i}$) would be an element of $T$ and hence would indeed satisfy $\scal{c}{c} \geqslant 0$, the third condition implies the second tautologically and the first condition implies the third by constructing $\langle X \rangle $ first as a pre-Hilbert module whose elements over an object $B_0$ are the linear combinations of the same shape as $c$ and then completing it into a Hilbert module. The uniqueness of $\langle X \rangle$ up to unique isomorphism follows from the fact that uniformly continuous maps (for example, bounded linear map) extend uniquely to maps between the completions.
}

}

\block{\Def{Let $\Cgo$ be a $C^{*}$-category, a Hilbert $\Cgo$-module is said to be \emph{small} if it can be ``generated by a set of elements'' as in the previous proposition.}

If $\Cgo$ is a small category, then every module over $\Cgo$ is small: $H$ can be generated  by $\cup_{A \in |\Cgo|} H(A)$, but for a non-small $C^{*}$-category this is a non trivial condition. In fact one has:

\Lem{If $G$ and $H$ are two Hilbert modules over a $C^{*}$-category, and if either $G$ or $H$ is small, then operators from $G$ to $H$ form a set. }

\Dem{Because the operation $*$ induces a bijection between operators from $G$ to $H$ and operators from $H$ to $G$ one can freely assume that it is $G$ which is small. But a linear map from $G$ to $H$ is then entirely determined by the choice of the image of each generator and hence the set of operators can be identified as a certain subset of this set of choice of image of each generator and hence is itself a set.}
}

\block{\label{DefnonsmallHC}\Def{If $\Cgo$ is a $C^{*}$-category, we denote by $\Hcal\Cgo$ the $C^{*}$-category of small Hilbert $\Cgo$-modules, that is those Hilbert modules that admit a set of generators.

We denote by $\Hcal' \Cgo$ the category of small Hilbert $\Cgo$-modules and bounded linear maps between them. It is \emph{not} a $C^{*}$-category.
}

This definition is in fact a bit too informal to make sense in our framework, indeed the ``class of all Hilbert modules'' actually does not exists. Hence we will fix a more precise model for the category $\Hcal \Cgo$:

Objects of $\Hcal \Cgo$ will be sets of generators for Hilbert $\Cgo$-modules as in \ref{genHilbmod}. Morphisms between two such sets $X$ and $Y$ will be functions $\lambda$ which to every pair $x \in X, y \in Y$ associate a $\lambda(y,x) \in \Cgo(s(x),s(y))$ such that there exists an operator $f_{\lambda}:\langle X \rangle \rightarrow \langle Y \rangle$ satisfying $\lambda(y,x) = \scal{e_y}{f(e_x)}$ (the map $f_{\lambda}$ being unique). One do the same for $\Hcal' \Cgo$.

The Yoneda embeddings take values in this $C^{*}$-category $\Hcal\Cgo$: if $\Cgo$ is unital then $\Yon_A$ is generated by a single element corresponding to the identity element of $A$, if $\Cgo$ is non-unital then $\Yon_A$ still admit a set of generators, for example one can take one generator in $\Yon_A(A)$ for each element of $\Cgo(A)$, or for each element in an approximate unit of $\Cgo(A)$.

}

\block{There is a small and uninteresting framework related difficulty that we did not mentioned in the above two paragraph and that we will discus now for the sake of completeness: if $f :A \rightarrow B$ is a bounded linear map between Hilbert module over a large category it is not clear that the proposition ``$f$ is an operator'' can be defined internally in our framework as it involves existential quantification on bounded linear map from $B$ to $A$ and that such map does not even form a class as we did not assume that there is a class of functions between two class. Fortunately, one can define that $f$ is an operator using only quantification over $A$ and $B$ (which is allowed as we assumed that the category of class if a \emph{Heyting} pre-topos). Indeed, $f$ is an operator if for any $X$ in the $C^{*}$-category, $\forall b \in B(X)$, $\exists a \in A(X)$, $\forall x \in A(x)$, $\scal{b}{f(x)} = \scal{a}{x}$. Indeed, assuming this condition one can then prove that for any $X$ and any $b$ in  $B(X)$, a $a$ satisfying this condition is unique, define $f^{*}(b)=a$ and show that this is a bounded linear map adjoint to $f$.

}

\block{If $f : \Cgo \rightarrow \Dgo$ is a $*$-functor then there is a canonical extension $\Hcal f : \Hcal \Cgo \rightarrow \Hcal \Dgo$ and a further extension $\Hcal' f: \Hcal' \Cgo \rightarrow \Hcal' \Dgo$. It can be seen on the presentation given above: If $(X,s,\scal{\_}{\_})$ is a set of generators for a Hilbert $\Cgo$-module then $(X,f(s),f(\scal{\_}{\_}))$ can be checked to be a set of generators for a Hilbert $\Dgo$-module. And if one has an operator $h$ from $\langle X \rangle$ to $\langle Y \rangle$ described by a function $\lambda$ as above then one can easily see that there exists an operator $\Hcal f (h)$ described by the function $f(\lambda)$. This clearly turn $\Hcal$ into a functor from $C^{*}$-category to $C^{*}$-category.

Moreover, we also note that if $f:\Cgo \rightarrow \Dgo$ is a $C^{*}$-functor, then the functor $(\Hcal f) : \Hcal \Cgo \rightarrow \Hcal \Dgo$ also acts on linear maps (even if they have no adjoint) in a way which is functorial, preserves the adjunction relation and send isometric inclusion to isometric inclusion.
}

\block{\label{Rq_boundedlinemapongen}Note that if $H$ is a Hilbert $\Cgo$-module given by a set of generators $H = \lbrace X \rbrace$ and $E$ is an arbitrary Hilbert $\Cgo$-module a bounded linear map from $H$ to $E$ is obviously determined by the image of the generators. Moreover, if one choose $e_x \in E(s(x))$ for each $x \in X$ this defines a bounded linear map $f$ from $H$ to $E$ such that $e_x = f(t_x)$ is and only if the following condition holds: There exists a constant $K$ such that for each $n$-tuple $(x_1,\dots,x_n) \in X^{n}$ one has the inequality in $M(s(e_{x_1}),\dots, s(e_{x_n}))^{+}$:  

\[[\scal{e_{x_i}}{e_{x_j}}]_{i,j} \leqslant K^{2} . [\scal{t_{x_i}}{t_{x_j}}]_{i,j}.  \]

Indeed, if $f$ exists then this equality follow from lemma \ref{Lem_boundedlinearmapIneq} below, and conversely if the inequality is satisfied for some constant $K$, then $f$ is defined on an element of the form $y=\sum_{i=1}^{n} t_{x_i} \lambda_i \in H(Y)$ with $\lambda_i \in \Cgo(Y,s(x_i))$ as $f(y)=\sum_{i=1}^{n} e_{x_i} \lambda_i \in E(Y)$ and the matrix inequality implies that:

\[ \scal{f(y)}{f(y)}=[\lambda_i^{*} \scal{e_{x_i}}{e_{x_j}} \lambda_j]_{i,j} \leqslant K^{2} . [\lambda_i^{*}\scal{t_{x_i}}{t_{x_j}}\lambda_j^{*}]_{i,j} = K^{2}\scal{y}{y} \]

Hence $f$ is bounded linear and hence extend into a bounded linear map from $\lbrace X \rbrace$ to $E$..
}

\block{\label{warningNonunital}It is important to note the following: if one works with non-unital categories and if $f: \Cgo \rightarrow \Dgo$ is a $C^{*}$-functor then the following square is in general \emph{not} commutative:
\[
\begin{tikzcd}[ampersand replacement=\&]
\Cgo \arrow{r}{f} \arrow{d}{\Yon} \& \Dgo \arrow{d}{\Yon} \\
\Hcal\Cgo \arrow{r}{\Hcal f} \& \Hcal\Dgo \\
\end{tikzcd}
\]

i.e. the Yoneda embeddings $\Yon$ might not be a natural transformation from the identity functor to the functor $\Hcal$. For example, if $f$ send every arrows to $0$ then $\Hcal f$ will send any objects to the zero module. 

The $C^{*}$-functors $f$ for which this square commute are called ``non-degenerate'' and this is equivalent to the fact that for each object $A$ of $\Cgo$ the corresponding morphism $f : \Cgo(A) \rightarrow \Dgo(f(A))$ is a non-degenerate morphism of $C^{*}$-algebras. It include in particular every unital $C^{*}$-functor between unital $C^{*}$-categories.
}

\subsection{Cones, ideals and hereditary subcategories}

\block{
\Def{
\begin{itemize}

\item If $A$ is a $C^{*}$-algebra, we denote by $A^{s}$ the set of self-adjoint elements of $A$ and by $A^{+}$ the subset of $A$ of positive self-adjoint elements of $A$. Self-adjoint elements are ordered by the relation $a \leqslant b$ if $b-a \in A^{+}$.

\item A hereditary cone of $A^{+}$ is a closed subset $X \subset A^{+}$ stable under linear combinations with real non-negative coefficients and such that if $a \in X$ and $b \leqslant a$ then $b \in X$.

\item A hereditary sub-algebra $B \subset A$ is a sub-$C^{*}$-algebra $B \subset A$ such that $BAB \subset B$.

\item If $\Cgo$ is a $C^{*}$-category a hereditary cone of $\Cgo$ is the data of a hereditary cone of $\Cgo(X)^{+}$ for each object $X \in |\Cgo|$.

\item A hereditary subcategory of $\Cgo$ is a sub-$C^{*}$-category $\Ago \subset \Cgo$ containing all the objects and such that $\Ago.\Cgo.\Ago \subset \Cgo$, i.e. for each $W,X,Y,Z \in |\Cgo|$ and for each $f \in \Ago(W,X), g \in \Cgo(X,Y), h \in \Ago(Y,Z)$ one has $h\circ g \circ f \in \Ago(X,Z)$. 

\end{itemize}
}

Also, a hereditary subcategory of a full subcategory of $\Cgo$ will be called a partial hereditary subcategory. Any partial hereditary category $\Ago$ can be considered as a hereditary category $\Ago'$ by adjoining all the missing objects and by defining $\Ago'(X,Y)= \Ago(X,Y)$ if both $X$ and $Y$ are objects of $\Ago$ and $\{0\}$ otherwise, or to say that without using the law of excluded middle:

\[ \Ago'(X,Y) = \{0\} \cup \bigcup_{X \in |\Ago|, Y \in |\Ago|} \Ago(X,Y) \]

Where the union is indexed by a subset of a singleton.
One can check that the category $\Hcal \Ago$ and $\Hcal \Ago'$ are equivalents: it follows for example from proposition \ref{PropHiwheniishereditary} below because the two sided ideal generated by $\Ago$ in $\Ago'$ is $\Ago'$ itself. For this reason $\Ago$ and $\Ago'$ will often be identified and we will for example say that $\Cgo(X)$ is a hereditary subcategory of $\Cgo$.

}

\block{\label{Lem_hered}Here are a few lemmas that allow to prove most of the result relating these notions:

\Lem{Let $\Cgo$ be a $C^{*}$ category, and $X,Y$ two objects of $\Cgo$. Let also $\Acal$ be a $C^{*}$-algebra.

\begin{enumerate}

\item \label{Lem_hered1} For any $a,b \in \Cgo(X,Y)$ one has:

\[ (a+b)^{*}(a+b) \leqslant 2 a^{*}a + 2b^{*}b \]
\[ (a+b)(a+b)^{*} \leqslant 2 aa^{*} + 2bb^{*} \]

\item \label{Lem_hered2} If $C$ is a hereditary cone of $\Acal^{+}$, $x \in C$ and $f: [ 0 , \Vert x \Vert ] \rightarrow \R^{+}$ is a morphism of locale such that $f(0)=0$ then $f(x) \in C$.

\item \label{Lem_hered3} Let $f \in \Cgo(X,Y), a \in \Cgo(X,X)^{+}, b \in \Cgo(Y,Y)^{+}$ such that $f^{*}f \leqslant a$ and $ff^{*}\leqslant b$ then:

\[ \lim_{n \rightarrow \infty} f a^{1/n} = \lim_{n \rightarrow \infty} b^{1/n} f = \lim_{n \rightarrow \infty} b^{1/n} f a^{1/n}  =f \]

\end{enumerate}

}

By $[ 0 , \Vert x \Vert ]$ we mean the the sub-locale of $\R$ defined as $\bigcap_{\Vert x \Vert < q} [0,q]$, i.e. the locale which classify the continuous real number $y$ such that $0 \leqslant y \leqslant x$.

\Dem{

\begin{enumerate}

\item One has:

\[(a+b)^{*}(a+b) + (a-b)^{*}(a-b)= 2a^{*}a +2b^{*}b \]

which directly proves the first equality and the second is obtained by replacing $a$ and $b$ with $a^{*}$ and $b^{*}$.

\item Let $f_n(t) = \min(f(t),n t)$, It is also a morphism of locales (because $t \mapsto nt$ and $\min$ are), and $(f_n)$ converge to $f$ uniformly (on $[0,\Vert x \Vert ]$) when $n$ goes to infinity. Hence $f_n(x)$ converge to $f(x)$, but $0 \leqslant f_n(x) \leqslant n x$ hence $f_n(x) \in C$ which concludes the proof.

\item We will first show that $f a^{1/n} \rightarrow f$. Indeed:

\[ \Vert f a^{1/n} - f \Vert  = \Vert (a^{1/n}-1) f^{*}f(a^{1/n}-1) \Vert^{1/2} \]

as $f^{*}f \leqslant a$ one has $(a^{1/n}-1) f^{*}f(a^{1/n}-1) \leqslant (a^{1/n}-1)a(a^{1/n}-1)$ but because the function $x(x^{1/n}-1)^{2}$ can be checked to converge to $0$ uniformly on $[0, \Vert a \Vert]$ this proves that $\Vert fa^{1/n}-f \Vert \rightarrow 0$. The case of $b^{1/n}f$ is dual, and for $b^{1/n} f a^{1/n}$ one can simply write that:

\[ b^{1/n} f a^{1/n} - f = ( b^{1/n} f a^{1/n} - b^{1/n} f ) + (b^{1/n} f - f)  \]

The second term is already known to tend to $0$, and for the first, the inequality:

\[ \Vert  b^{1/n} f a^{1/n} - b^{1/n} f \Vert \leqslant \Vert b \Vert ^{1/n} \Vert fa^{1/n} - f \Vert \]

allows to conclude.

\end{enumerate}

}

}

\block{\label{hered_Imp_cone}\Prop{Let $\Ago$ be a hereditary subcategory of $\Cgo$, then $\Ago^{+}$ is a hereditary cone of $\Cgo$ and an element $a \in \Cgo(X,Y)$ is in $\Ago$ if and only if $(a^{*}a) \in \Ago(X)^{+}$ and $(aa^{*}) \in \Ago(Y)^{+}$. }

\Dem{ $\Ago^{+}(X)$ is clearly a cone, we need to prove that it is hereditary, but if $b \leqslant a$ with $b \in \Cgo(X)$ and $a \in \Ago(X)^{+}$ then by lemma \ref{Lem_hered}(\ref{Lem_hered3}) one has $a^{1/n} b^{1/2}a^{1/n} \rightarrow b^{1/2}$ but $a^{1/n} \in \Ago$ hence this proves that $b^{1/2} \in \Ago$ and hence that $b \in  \Ago$.

If $a \in \Ago$ then $a^{*}a$ and $aa^{*}$ are in $\Ago^{+}$. Conversely, if $a^{*}a$ and $aa^{*}$ are in $\Ago$, then by \ref{Lem_hered}(\ref{Lem_hered3}) one has $a = \lim (aa^{*})^{1/n} a (a^{*}a)^{1/n}$ and $(aa^{*})^{1/n} a (a^{*}a)^{1/n} \in \Ago$, which concludes the proof.

}

}

\block{\Def{Let $\Cgo$ be a $C^{*}$-category.

\begin{itemize}

\item A left (resp. right) ideal of $\Cgo$ is an additive (generally non-unital) sub-category $I \subset \Cgo$ (containing all objects) such that for all objects $X, Y \in |\Cgo|$, $I(X,Y)$ is a closed subspace of $\Cgo(X,Y)$ and $\Cgo.I \subset I$ (resp. $I.\Cgo \subset I$) i.e. $I$ is stable under post-composition (resp. pre-composition) by arbitrary morphism.

\item A two sided ideal is an additive subcategory which is both a left and right ideal.

\item If $X$ is an object of $\Cgo$, a left ideal on $X$ is a left ideal $I \subset \Cgo$ which only contains arrows from $X$ and zero morphisms. Similarly a right ideal on $X$ is a right ideal which only contains arrows into $X$ and zero morphisms.

\end{itemize}
}

Note that a left ideal $I$ on $\Cgo$ is the same as a collection of a left ideals $I_X$ on $X$ for each object object $X \in |\Cgo|$ and dually a right ideal is the same as a collection of right ideals on every object.

Also a right ideal on $X$ for an object $X \in \Cgo$ is exactly the same as a Hilbert sub-module of $\Yon(X)$.

}

\block{\label{IdealGeneratedByCone}\Prop{Let $C$ be a hereditary cone of a $C^{*}$-category $\Cgo$. Then:

\begin{itemize}

\item The set $LC$ of $f \in \Cgo$ such that $ f^{*}f \in C$ is a left ideal of $\Cgo$. Moreover it is the left ideal generated by $C$ in the sense that for each $X,Y \in |\Cgo|$, $LC(X,Y)$ is the closure of the set of $f \circ c$ for $c\in C(X)$ and $f \in \Cgo(X,Y)$.

\item The set $RC$ of $f \in \Cgo$ such that $ ff^{*} \in C$ is a right ideal of $\Cgo$. Moreover it is the left ideal generated by $C$ in the sense that for each $X,Y \in |\Cgo|$, $RC(X,Y)$ is the closure of the set of $c \circ f$ for $c\in C(X)$ and $f \in \Cgo(X,Y)$.

\item The set $HC = RC \cap LC$ of $f$ such that both $ff^{*}$ and $f^{*}f$ are in $C$ is a hereditary sub-$C^{*}$-category and $HC(X,Y)$ is the closure of the set of element of the form $c \circ f \circ c'$ with $c \in C(Y), c' \in C(X)$ and $f \in  \Cgo(X,Y)$.

\end{itemize}

}

Note that $LC$ and $RC$ are ideals on an object $X \in \Cgo$ if and only if $C$ is included in $\Cgo(X)^{+}$, seen as a hereditary cone of $\Cgo$, i.e. if $C \subset \Cgo(X)^{+} \bigcup \{ 0_Y, Y \in |\Cgo| \}$.

\Dem{For the first point, the stability by addition follow directly from \ref{Lem_hered}(\ref{Lem_hered1}), the stability under post-composition comes from $(fc)^{*}(fc) = c^{*}f^{*}fc \leqslant K c^{*}c \in C$ and the closedness is obvious. $LC$ contains $C$ because if $c \in C$ then $c^{*}c = c^{2} \in C$ by \ref{Lem_hered}(\ref{Lem_hered2}) hence $c \in LC$, hence elements of the form $(fc)$ with $c \in C$ are in $LC$ because it is a left ideal. Conversely if $f^{*}f \in C$ then by \ref{Lem_hered}(\ref{Lem_hered2}) $(f^{*}f)^{1/n} \in C$ and hence \ref{Lem_hered}(\ref{Lem_hered3}) exhibit $f$ as a limit of a sequence of elements of the form $f \circ c$ with $c \in C$.
The second point is just the dual of the first.

For the third point, the stability under $*$ of $HC$ is obvious from the definition, the fact that it is a hereditary subcategory follows directly from the previous two points, it contains the elements of the form $c f c'$ because it is a hereditary subcategory and we already noticed that $C \subset LC \cap RC$. Conversely, if $f \in HC$ then by the same argument as for $LC$, lemma \ref{Lem_hered}(\ref{Lem_hered3}) exhibit $f$ as a limit of a sequence of element of the form  $c f  c'$.

}

}

\block{\label{CorBijHered}\Cor{Let $\Cgo$ be a $C^{*}$-category, the construction $\Ago \mapsto \Ago^{+}$ and $C \mapsto HC$ induces a bijection between the class of hereditary subcategories of $\Cgo$ and the class of hereditary cone of $\Cgo$.}

This result automatically extends to partial hereditary subcategories and ``partial cones'' by applying it to each full subcategory.

\Dem{The third point of proposition \ref{IdealGeneratedByCone} prove that $HC$ is a hereditary subcategory, proposition \ref{hered_Imp_cone} show that $\Ago^{+}$ is hereditary and that $H\Ago^{+} = \Ago$, and finally $(HC)^{+} = C$ because lemma \ref{Lem_hered}(\ref{Lem_hered2}) imply that for a positive element $c$ if $c^{*}c = c^{2} \in C$ then $c = \sqrt{c^{2}} \in C$ and conversely that if $c \in C$ then $c^{*}c = cc^{*} = c^{2} \in C$. }

}

\block{\label{IeqL(Iplus)}

\Prop{If $I$ is a left ideal of $\Cgo$ then $I^{+} = I \cap \Cgo^{+}$ is a hereditary cone of $\Cgo$ and $I = L(I^{+})$. This induces a bijection between the hereditary cones of $\Cgo$ and its left ideals and a bijection between hereditary cones on $\Cgo(X)$ and left ideals on $X$ for any object $X \in |\Cgo|$.}

Putting this together with the previous results one has compatible bijections between the right ideals, left ideals, hereditary subcategory and hereditary cone. One pass from left ideal to right ideal by $I \mapsto I^{*}$ and from left or right ideal to hereditary sub-algebra by $I \mapsto I \cap I^{*}$.

\Dem{Let $I$ be a left ideal then $I \cap I^{*}$ is easily checked to be a hereditary sub-algebra, hence $I^{+} = (I \cap I^{*})^{+}$ is a hereditary cone by proposition \ref{hered_Imp_cone}. If $f \in I$ then $f^{*}f \in I^{+}$ hence $f \in L(I^{+})$ conversely if $f \in L(I^{+})$ then by proposition \ref{IdealGeneratedByCone} $f$ is in the closure of elements of the form $f \circ c$ for $c \in I^{+}$ hence is in $I$. Finally those bijections reduces to a bijection between ideals on an object $X \in |\Cgo|$ and hereditary cones of $\Cgo(X)$.
} }

\block{\Cor{A two sided ideal of a $C^{*}$-category is stable under conjugation.}
In particular it is a sub-$C^{*}$-category, even an hereditary one.
\Dem{If $I$ is a two-sided ideal then it is both a left and right ideal hence it is both equal to $L(I^{+})$ and $R(I^{+})$ by proposition \ref{IeqL(Iplus)}. As $L(I^{+})^{*} = R(I^{+})$ this proves that $I^{*} = I$.}
}

\block{\Prop{\begin{itemize}

\item If $\Acal$ is a $C^{*}$-algebra, and $\Bcal \subset \Acal$ is a sub-$C^{*}$-algebra. Then $\Bcal$ is a hereditary sub-algebra if and only if $\Bcal^{+}$ is a hereditary cone of $\Acal^{+}$.

\item If $\Ago \subset \Cgo$ is a sub-$C^{*}$-category containing all objects, then $\Ago$ is a hereditary subcategory if and only for all pair of elements $X,Y \in \Cgo$, $M_\Ago(X,Y)^{+}$ is a hereditary cone of $M_\Cgo(X,Y)^{+}$. 

\end{itemize}
}
\Dem{The second point follows from the first: indeed assume that $M_\Ago(X,Y)$ is a hereditary sub-algebra of $M_\Cgo(X,Y)$ for each $X,Y \in |\Cgo|$. Let $f \in \Ago(W,X)$, $h \in \Cgo(X,Y)$ and $g \in \Ago(Y,Z)$ for some $W,X,Y,Z \in |\Cgo|$. Then by (the proof of) lemma \ref{Lem_hered}(\ref{Lem_hered3}):

\[ \lim_{n \rightarrow \infty} g (g^{*}g)^{1/n} h (ff^{*})^{1/n} f = g \circ  h \circ f \]

But because $M_\Ago(X,Y)$ is hereditary in $M_\Cgo(X,Y)$ one deduces that:
\[(g^{*}g)^{1/n} h (ff^{*})^{1/n} \in \Ago(X,Y)\]
 and hence this concludes the proof of the second point (assuming the first).

\bigskip

For the first point, it is already proved (proposition \ref{hered_Imp_cone}) that if $\Bcal$ is hereditary then $\Bcal^{+}$ is a hereditary cone. We need to prove the converse. Let $\Bcal \subset \Acal$ a sub-$C^{*}$-algebra such that $\Bcal^{+}$ is a hereditary cone of $\Acal^{+}$.

Remark first that if $f \in \Acal^{+}$ then for any $h \in \Bcal$ one has $h^{*} f h \leqslant K h^{*} h \in \Bcal$, where $K$ is a continuous real number such that $\Vert f \Vert \leqslant K$, hence $h^{*} f h \in \Bcal$. If now $f$ is an arbitrary element of $\Acal$ then it can be decomposed into a linear combination of positive self-adjoint elements hence if $h \in \Acal$ one still has $h^{*}f h \in \Bcal$. Finally let $x,y \in \Bcal$ and $f \in \Acal$. let $b=xx^{*}+y^{*}y$ by lemma \ref{Lem_hered}(\ref{Lem_hered3}) one has:

\[ \lim_{n \rightarrow \infty } b^{1/n} x f y b^{1/n} = x f y \]

and $b^{1/n} x f y b^{1/n} \in \Bcal$ which concludes the proof.
}

Note that the second point cannot be weakened to look more like the first: one can for example consider a sub-category $\Ago$ of $\Cgo$ which contains all the endomorphisms of objects and only the zero morphisms between distinct objects. $\Ago(X)$ is hereditary in $\Cgo(X)$ for each $X$  but $\Ago$ will not be a hereditary subcategory because $ff^{*}, f^{*}f \in \Ago$ will not imply $f \in \Ago$ (unless $\Cgo$ has only zero morphisms between distinct objects).
}

\block{\label{Def_AH}\Def{Let $\Cgo$ be a $C^{*}$-category, $H$ a Hilbert $\Cgo$-module and $\Acal$ an hereditary sub-algebra of operators of $H$. One defines the sub-module $\Acal H$ of $H$ by:

\[ \Acal H (X) = \{h \in H(X) | h h^{*} \in \Acal \} \]

Or equivalently as the closure in $H$ of elements of the form $a_X h$ for $a \in \Acal$ , $X \in |\Cgo|$ and $h \in H(X)$.
}
The fact that it is a sub-module and the equivalent presentation follow directly from the proof of proposition \ref{IdealGeneratedByCone}.

}

\block{\label{AHis_small}\Prop{If $H$ is a small Hilbert $\Cgo$-module and $\Acal$ is a hereditary sub-algebra of operators on $H$ then $\Acal H$ is small. More precisely if $X$ is a set of generators for $H$ the set of $a.x$ for $a \in \Acal$ and $x \in X$ form a set of generators for $\Acal H$.}

Hence, this construction is compatible with our convention that in $\Hcal \Cgo$ every object is endowed with the choice of a set of generators.

\Dem{As $H$ is small then $\Acal$ is a set (as a subset of a set). Take $(x_i)_{i \in I}$ be a set of generators of $H$, then the $a x_i$ for $a \in \Acal$ form a set of generators of $\Acal H$. Indeed, every element of $\Acal H$ can be approximated by elements of the form $a h$ for $a \in \Acal$ and $h \in H$ and every element $h$ of $h$ can be approximated by linear combination of element of the form $x_i g$, hence any element of $\Acal H$ can be approximated by linear combination of elements of the form $a x_i g$, which proves that the $a x_i$ are generators for $\Acal H$. }

}

\block{\label{Prop_H'=AH}\Prop{Let $H$ be a small Hilbert $\Cgo$-module and $H'$ a sub-module of $H$. Let $\Acal$ be the sub-algebra of operators of $H$ of $a$ such that both $a$ and $a^{*}$ takes their values in $H'$. Then $\Acal$ is a hereditary sub-algebra of $\Hcal\Cgo(H)$ and $H' = \Acal H$.}

\Dem{$\Acal$ is clearly a sub-$C^{*}$-algebra (it is a set because $H$ is small) and one has $\Acal.\Hcal \Cgo(H) \Acal \subset \Acal$. It is also clear that $\Acal H \subset H'$ as any element of the form $a h$ for $a \in \Acal$ is in $H'$, and conversely, if $h \in H'$ then one has $hh^{*} \in \Acal$ and hence $h \in \Acal H$.
}
}

\block{\Cor{A sub-module of a small Hilbert $\Cgo$-module is small.}

\Dem{If $H' \subset H$ is a sub-module then $H' = \Acal H$ for some $\Acal$ by proposition \ref{Prop_H'=AH} and we pointed out in proposition \ref{AHis_small} that a module of this form is small.
}
}

\block{
\Prop{Let $\Cgo$ be a $C^{*}$-category, and $X$ an object of $\Cgo$. The construction $\Acal \mapsto \Acal \Yon_X $ induces a bijection between the set of hereditary sub-algebras of $\Cgo(X)$ and Hilbert sub-$\Cgo$-modules of $\Yon(X)$.
}

\Dem{This is a reformulation of proposition \ref{IeqL(Iplus)} through the (tautological) identification of sub-module of $\Yon(X)$ with right ideals on $X$ and the identification of hereditary cones of $\Cgo(X)^{+}$ with hereditary sub-algebras of $\Cgo(X)$ constructed in corollary \ref{CorBijHered}.
}
}

\block{\Def{Let $\Cgo$ be a $C^{*}$-category, $H$ and $H'$ two Hilbert $\Cgo$-modules, an operator $f$ between $H$ and $H'$ is said to be compact or $\Cgo$-compact if for all $\epsilon>0$, $f$ can be approximated up to $\epsilon$ by a linear combination of operators which can be factored into $h c g$ with $c \in \Cgo$. The (non-unital) $C^{*}$-category whose objects are objects of $\Hcal \Cgo$ and morphism are compact operators is denoted $\Kcal\Cgo$}

$\Kcal\Cgo$ is a (non unital) $C^{*}$-category and a two sided ideal of $\Hcal\Cgo$ which contains $\Cgo$, It is the smallest such two sided ideal. The set of compact operators between two Hilbert $\Cgo$-modules will occasionally be denoted $\Kcal(H,H')$ instead of $\Kcal \Cgo(H,H')$. The finite linear combinations of operators which factor into an operator in $\Cgo$ will be called the finite rank operators.

}

\block{The terms $\Cgo$-compact is here to avoid confusion with absolutely compact operators that will be defined later, and as well as because sometimes a same $C^{*}$-category can be both of the form $\Hcal \Cgo$ and a full subcategory of $\Hcal \Dgo$ for some $\Dgo$, in which case it can have two different notion of compactness: $\Dgo$-compactness and $\Cgo$-compactness. In fact one can define more generally:

\Def{For class $C$ of map in a $C^{*}$-category $\Cgo$, the elements of the total two sided ideals of $\Cgo$ generated by $C$ will be said to be $C$-compact. If $C = \{ Id_X \}$ or equivalently $C = \Cgo(X)$ we will talk about $X$-compact operators.}

}

\block{\label{densityforHmodule}\Prop{Let $H$ be a Hilbert $\Cgo$-module, and $\Acal$ a hereditary sub-algebra of operators of $H$. The following conditions are equivalent:

\begin{itemize}

\item $\Acal$ contains the algebra $\Kcal(H)$ of ($\Cgo$)-compact operators of $H$.

\item $\Acal$ is dense for the topology of point-wise convergence\footnote{We mean that for any $f$ operators on $H$, for any $h_1, \dots h_n$ elements of $H$, for any $\epsilon >0$ there is a $a \in A$ such that $\Vert f(h_i) - a(h_i) \Vert<\epsilon$. }.

\item $\Acal H = H$.

\end{itemize}

}
\Dem{By definition, $\Acal H = H$ if and only if for all $h \in H$, $hh^{*} \in \Acal$ and the algebra of compact operator is generated (as a hereditary sub-algebra) by the $hh^{*}$, hence the first and the third propositions are equivalent. The first proposition implies the second because the compact operators are dense for the topology of simple convergence and the second implies the third because $\Acal H(X)$ can be described as the closure of the set of $a .h$ for $a \in \Acal$ and $h \in H(X)$.
}

}

\subsection{Restriction of Hilbert modules}
\label{Restricton_of_Hilbert_modules}

\blockn{With ordinary category, if $F: C \rightarrow C'$ is a functor and $P$ is a pre-sheaf on $C'$ then $F^{*} P = P \circ F$ is a pre-sheaf on $C$. This construction does not make sense with $C^{*}$-categories because there is in general no way to pullback the $C'$ valued scalar product into a $C$ valued one. Their is a solution to this problem when $F$ is the inclusion of a partial hereditary subcategory.}

\block{\Def{Let $\Cgo$ be a $C^{*}$-category and $\Ago$ be a partial hereditary sub-category, and $X \in |\Cgo|$ then we define the restriction of $X$ to $\Ago$ as the following Hilbert $\Ago$-module:

\[ R_{\Ago}(X)(a) = \{ f:a \rightarrow X | f^{*}f \in \Ago \} \]

}

This is a Hilbert module because if $f,g \in R_{\Ago}(X)(a)$ then $f+g$ as well because of \ref{Lem_hered}(\ref{Lem_hered1}), and $\scal{f}{g} := f^{*}g$ satisfy:

\[ (f^{*}g)^{*} (f^{*}g) =g^{*}ff^{*}g \leqslant K g^{*}g \in \Ago \]
\[ (f^{*}g) (f^{*}g)^{*} = f^{*}gg^{*}f \leqslant K' f^{*}f \in \Ago \]

hence is in $\Ago$ by proposition \ref{hered_Imp_cone}.

This is functorial in both $a \in |\Ago|$ and in $X \in |\Cgo|$ because the set of $f$ such that $f^{*}f \in \Ago$ is a left ideal (proposition \ref{IdealGeneratedByCone}). It is complete because it is closed in $\Cgo(a,X)$.}

\blockn{There is two obstructions into making $R_{\Ago}$ into a functor from $\Cgo$ to $\Hcal \Ago$: one need to prove that for each $c \in \Cgo$ $R_{\Ago}(c)$ is small and one need for each $c$ to actually fix a set of generators (because of our definition of $\Hcal \Ago$ as detailed just below \ref{DefnonsmallHC}). The first problem will be solved in the next paragraph using the collection axiom \footnote{axiom $(A7)$ in the appendix \ref{Def_class_of_smallmaps}, see also \ref{Description_Collection}} while the second just force us to only defines $R$ as an anafunctor (as mentioned in the introduction, see \cite{makkai1996avoiding} for details about anafunctors). }

\block{\label{Prop_restrictionaresmall}\Prop{For any $C^{*}$-category $\Cgo$,  any partial hereditary subcategory $\Ago \subset \Cgo$ and any $X \in |\Cgo|$ the Hilbert $\Ago$-module $R_{\Ago}(X)$ is small.}

\Dem{Let $V$ be the class of all pair $(A,h)$ for $A \in |\Ago|$ and $h \in R_{\Ago}(X)$, i.e. $h \in \Cgo(A,X)$ such that $h^{*}h \in \Ago(A)$. There is map $h \mapsto hh^{*}$ from $V$ to $\Cgo(X)$, let $V'$ its image which is small because $\Cgo(X)$ is. One can then apply the collection axiom (see \ref{Description_Collection}) and one get a set $V''$ with a map to $V$ such that the composite to $V'$ is surjective. We claim that $V''$ is a set of generators for $R_{\Ago}(X)$. Indeed, let  $H \subset R_{\Ago}(X)$ be the sub-Hilbert module of $R_{\Ago}(X)$ generated by $V''$, by construction the $C^{*}$ algebra $K'$ generated by the $hh^{*}$ for $h \in H$ is equal to $K$ the algebra of compact operators of the Hilbert $\Ago$-module $R_{\Ago} (X)$.
If $h$ is an an arbitrary element of $R_{\Ago}(X)$ then $hh^{*} \in K$ and by lemma \ref{Lem_hered}(\ref{Lem_hered3}), $h = \lim (hh^{*})^{\frac{1}{n}} h $. But as the $(hh^{*})^{\frac{1}{n}}$ are in $K$ they can be written as combination of composition of $aa^{*}$ for $a \in H$ and this implies that all the $(hh^{*})^{\frac{1}{n}} h $ are in $H$ and hence that $h \in H$ which proves that $H = R_{\Ago}(X)$ and concludes the prof.
}

}

\block{\label{Cor_restrictionfctexist}\Cor{For any $C^{*}$-category $\Cgo$ and any partial hereditary sub-category $\Ago \subset \Cgo$ there exists a restriction (ana)functor $R_{\Ago}: \Cgo \rightarrow \Hcal \Ago$.} 

This is typical one place where one really need to use anafunctors: indeed, the previous proposition show that $R_{\Ago}(X)$ is small for all $X$ but does not provide a canonical choice of a set generators hence does not allow to construct a functor from $\Cgo$ to $\Hcal \Ago$ (we recall from \ref{DefnonsmallHC} that any object of $\Hcal \Ago$ is endowed with a fixed of generator).

\bigskip

If $\Ago \subset \Cgo$ is a partial hereditary subcategory we will also denote by $R_{\Ago}$ the restriction functor from $\Hcal \Cgo$ to $\Hcal \Ago$. It is a special case of the previous construction by considering $\Ago$ as a subcategory of $\Hcal \Cgo$ through the Yoneda embeddings of $\Cgo$, but it can also be seen as taking the elements of a Hilbert $\Cgo$-module whose scalar product with themselves lies in the subcategory $\Ago$.

\bigskip

The reader should also note that in a lot of concrete situation, $R_{\Ago}$ can indeed be made into a functor and we do not need the collection axiom to construct it: For example if $\Ago$ is small then one can choose the set of all element of $R_{\Ago}(X)$ and it gives us a canonical set of generators. Also\footnote{For this remark, we need to assume that a set indexed union of sets is a set ``in a canonical way'', and not just a small class, i.e. that the union can be represented by a map from the class of family of set to the class of sets. This is not the case in the basic framework described in appendix $B$, but it is for example the case in models of NBG, and one can always make it true by changing the universe.} for any $C^{*}$-category $\Cgo$, the functor:

\[R_{\Cgo} : \Hcal \Hcal \Cgo \rightarrow \Hcal \Cgo \]

always exists ($\Cgo$ is considered as a partial hereditary subcategory of $\Hcal \Cgo$ under the Yoneda embeddings). Indeed if one has $V \in |\Hcal \Hcal \Cgo |$ then $V$ has a (canonical) set of generators as a Hilbert $\Hcal \Cgo$-module, each generator $x$ is a an element of $V(s(x))$ for a well defined object $s(x) \in |\Hcal \Cgo|$ which also has a canonical set of generators and each of this generator $w$ is an arrow $w:s(w) \rightarrow s(x)$ by composing one obtains an arrow $w:s(w) \rightarrow V$ which is an element of $R_{\Cgo}(V)$ and one easily see that taking all those elements (for all $x$ and all $w$) is a set of generators for $V$.

}

\block{\label{Lem_boundedlinearmapIneq}\Lem{Let $\Cgo$ be a $C^{*}$-category, $X$ and $Y$ two Hilbert $\Cgo$-modules and $f: X\rightarrow Y$ a bounded linear transformation such that $\Vert f \Vert \leqslant K$. Then for any $c \in |\Cgo|,a \in X(c)$, one has:

\[ \scal{f(a)}{f(a)} \leqslant K^{2} \scal{a}{a} \]

}

\Dem{If $f$ has an adjoint, then this is classical: $\scal{f(a)}{f(a)}$ can be written as $a^{*}f^{*}fa$ in the category of Hilbert modules and by the usual properties of positivity $a^{*}f^{*}fa \leqslant a^{*}K^{2}a = K^{2}\scal{a}{a}$. Now let $k$ be a finite rank endomorphism of $X$, i.e. a finite linear combination of operators of the form $xx'^{*}$ of $x,x'$ elements of $X$, then $fk$ is an operator from $X$ to $Y$ (indeed $f\circ x \circ x'^{*} = f(x) \circ x'^{*}$ hence has an adjoint). Now for any $a \in X$ one can find such finite rank operators $k$ such that $\Vert k \Vert \leqslant 1$ and $\Vert a - ka \Vert \leqslant \epsilon$ for example by approximating $(aa^{*})^{1/n}$ by a finite rank operators which are dense among compact operators, for such a $k$ one will obtain by applying the previous equality to $fk$ that $\scal{f(a)}{f(a)} \leqslant K^{2} \scal{a}{a} - 2\epsilon \Vert a \Vert \Vert f \Vert^{2}$, this being true for any $\epsilon$ it concludes the proof.
}

}

\block{\Cor{Let $\Ago \subset \Cgo$ be a partial hereditary sub-category of a $C^{*}$-category $\Cgo$, then the restriction (ana)functor $R_{\Ago} : \Hcal \Cgo \rightarrow \Hcal \Ago$ can be extended into a restriction (ana)functor: $R'_{\Ago}: \Hcal' \Cgo \rightarrow \Hcal' \Ago$. }

\Dem{This follows immediately from lemma \ref{Lem_boundedlinearmapIneq}: if $f :X \rightarrow Y$ is a bounded linear map between two Hilbert $\Cgo$-modules then any element $x \in R_{\Ago}(X)$ is just an element $x \in X$ such that $\scal{x}{x} \in \Ago$ hence by the lemma $\scal{f(x)}{f(x)} \in \Ago$ and hence $f(x) \in R_{\Ago}(Y)$ and this defines a bounded (by the norm of $f$) linear map in $\Hcal' \Ago$. }
}

\block{\Prop{Let $\Ago \subset \Cgo$ be a partial hereditary subcategory, with $i : \Ago \hookrightarrow \Cgo$ the inclusion $*$-functor. Then $R'_{\Ago}$ is a right adjoint to $\Hcal' i$. }

Note that this is almost never true if one only look at $R_{\Ago}$ and $\Hcal i$: the co-unit of the adjunction will not be an operator but an isometric inclusion (see the next proposition) hence there is no adjunction if one only look at operators.

\Dem{Let $H$ be a Hilbert $\Ago$-module and $E$ be a Hilbert $\Cgo$-module. $\Hcal i H$ is generated by the element of $H$ (or a set of generators for $H$), hence, using the remark made in \ref{Rq_boundedlinemapongen}, a bounded linear map from $\Hcal' i H$ to $E$ is given by taking one element in $E$ for each element in $H$ satisfying the matrix inequality mentioned in \ref{Rq_boundedlinemapongen}. Because of that inequality, those elements in $E$ are automatically in $R'_{\Ago}(E)$, hence a bounded linear map from $\Hcal' i H$ to $E$ is exactly the same as a bounded linear map from $ H$ to $R_{\Ago}E$.}

}

\block{\label{PropHiwheniishereditary} With the same notation as the previous proposition, and let $\Igo$ be the (total) two sided ideal of $\Cgo$ generated by $\Ago$.

\Prop{The functor $(\Hcal' i) : \Hcal' \Ago \rightarrow \Hcal' \Cgo$ is fully faithful and the unit of the adjunction of the previous proposition is an isometric isomorphism $R'_{\Ago} \circ (\Hcal' i) \simeq Id_{\Hcal' \Ago}$. For $H \in |\Hcal \Cgo|$ a small Hilbert $\Cgo$-module the following conditions are equivalents:

\begin{enumerate}

\item $H$ is in the (essential) image of $\Hcal i$.

\item The scalar product of $H$ takes values in $\Igo$

\item $H$ is generated by elements $h$ such that $\scal{h}{h} \in \Ago$.

\end{enumerate}

And moreover, the co-unit of the adjunction $(\Hcal' i) \circ R'_{\Ago}(H) \hookrightarrow H$ corresponds to the largest sub-module of $H$ satisfying these conditions.

}

The reader should be careful with the fact that the functor $(\Hcal i)$ is not compatible with the Yoneda embeddings (it does not send a representable in $\Hcal \Ago$ to its image by $i$, see \ref{warningNonunital}), this distinction actually play a very important role in this proposition. For example, take $\Cgo$ to be the $C^{*}$-algebra $B(H)$ of bounded operator on some (infinite dimensional) Hilbert space $H$ and $\Ago = \mathbb{C} = p B(H) p$ for some rank one projection $p \in B(H)$. Then $\Igo$ is the algebra $\Kcal(H)$ of compact operators of $H$. A Hilbert $\C$-module is the same as a Hilbert $\Kcal(H)$-module and they form a full subcategory of Hilbert $B(H)$-modules (they are exactly the reflexive modules). Finally the representable Hilbert $\C$-module $\C$ is not identified with the representable Hilbert $B(H)$-module $B(H)$ but with the (non-representable) Hilbert $B(H)$-module $H$.

\Dem{We will ignore the distinction between $\Hcal$ and $\Hcal'$ as well as between $R$ and $R'$ in order to simplify notations.

For any element $h \in H(X)$ of a Hilbert module we will denote $|h| = \scal{h}{h}^{\frac{1}{2}} \in \Cgo(X)$.

We will first prove that the unit of the previous adjunction induce an isometry $R_{\Ago} \circ (\Hcal i) \simeq Id_{\Hcal \Ago}$. Let $H$ be a small Hilbert $\Ago$-module, the unit of the adjunction is the map which send any element of $H$ to the corresponding element of $\Hcal_i H$, which happen to be in $R_{\Ago} \circ (\Hcal i ) H$ this clearly preserve the scalar product, hence we just need to show that this is a surjection. For any $X \in |\Cgo|$, the space $(\Hcal i)(H)(X)$ is spammed by elements of the form $h.f$ where $h$ is in $H(Y)$ for some $Y \in |\Ago|$ and $f \in \Cgo(X,Y)$. Hence $R_{\Ago} \circ (\Hcal i)(H)(X)$ when $X \in |\Ago|$ is spammed by elements of the form $h.f.a$ with $h$ and $f$ as before and $a$ an element of $\Ago(X)$. Now, by lemma \ref{Lem_hered}(\ref{Lem_hered3}) applied to the category of modules, $h = \lim h |h|^{1/n}$, hence $R_{\Ago} \circ (\Hcal i)(H)(X)$ is spammed by elements of the form $h | h |^{1/n}.f.a$ but $|h|^{1/n}.f.a$ is in $\Ago$ because $\Ago$ is hereditary hence $h | h |^{1/n}.f.a$ is in $H(X)$, which concludes the first part of the proof.

\bigskip

We will now prove the equivalence of the three conditions. If $H$ is isomorphic to $(\Hcal i)(H')$ then for each $X \in |\Cgo|$, $H(X)$ is generated by elements of the form $hf$ for $Y \in \Cgo$, $h \in H'(Y)$ and $f : a \rightarrow X$. The scalar product of two such elements is $\scal{hf}{h'f'} = f^{*}\scal{h}{h'}f'$ which is an element of $\Igo$. By linearity and continuity this prove that $H$ satisfies the condition $(2)$.

Assume now that $H$ satisfies the condition $(2)$, then for any $X \in |\Cgo|$ and any $x \in H(X)$ one has $\scal{x}{x} \in \Igo$ hence by lemma \ref{Lem_hered}(\ref{Lem_hered3}) any $x$ in $H(X)$ can be approximated by $x|x|^{1/n}$ . By functional calculus, $|x|^{1/n} \in \Igo$ hence can be approximated by linear combination of element of the form $fag$ for $a \in \Ago$ and $f,g$ arbitrary arrows in $\Cgo$ such that this composition makes sense, in particular, $H$ is generated by element of the form $xfag$ for $a \in \Ago$. Rewriting this by forgetting the $g$ on the right and renaming $xf$, one obtains that $H$ is generated by elements of the form $va$ for $a \in \Ago$ and $v \in H(X)$ where $X$ is the target object of $a$.

Finally let $H$ be any Hilbert $\Cgo$-module and let $H'$ be its restriction to $\Ago$. Any element of $H$ such that $\scal{h}{h} \in \Ago$ corresponds by definition to an element of $H'$ and their scalar product in $H'$ and in $H$ are the same (up to identifying elements of $\Ago$ with their image in $\Cgo$). Hence $(\Hcal i)(H')$ identifies with the sub-module of $H$ generated by elements $h$ such that $\scal{h}{h} \in \Ago$. In particular, if $H$ satisfies $(3)$ then $H \simeq (\Hcal i) H'$ is in the essential image of $\Hcal i$ and $(\Hcal i)(H')$ is the largest sub-module of $H$ satisfying $(3)$.

This concludes the proof of the equivalence of the three conditions on $H$ and of the final claim.

}

}

\block{\label{Cor_eqHA_HI}\Cor{If $\Ago \subset \Cgo$ is a hereditary subcategory which generates $\Cgo$ as a two sided ideal then one has an equivalence of category:

\[ \Hcal' \Cgo \simeq \Hcal' \Ago \]

which preserve the adjunction relation and is induced by the restriction to $\Ago$.

In particular for any $C^{*}$-category $\Cgo$, one has an equivalence: 

\[ \Hcal' \Cgo \simeq \Hcal' \Kcal \Cgo.\]
}

\Dem{If $\Cgo = \Igo$ then any Hilbert $\Cgo$-module satisfies the condition $(2)$ of proposition \ref{PropHiwheniishereditary} and hence the two functor $R_{\Ago}$ and $\Hcal i $ are inverse of each other.

Finally $\Kcal \Cgo$ is exactly the two sided generated by $\Cgo$ in $\Hcal \Cgo$ and $\Cgo$ is a full subcategory of $\Kcal \Cgo$ hence is in particular a partial hereditary subcategory so this is a special case of the previous situation.}
}

\section{Complete $C^{*}$-categories}
\label{section_completeness}

\blockn{In this section we are interested in finding a correct analogue of the notion of categorical ``completeness'' for $C^{*}$-category, in the sense of existence of limits/co-limits. For example categories like $\Hcal \Cgo$ are expected to be complete and tensorisation by bi-module are expected to be (co)limits preserving functors. The notion we are after in the present paper appears to actually be more related to ``co-completeness'' (i.e. existence of co-limits) than completeness; but, because of the self duality of $C^{*}$-categories induced by the $*$ operation, there is, at least at the time this paper is being written, no reasons to think that in the world of $C^{*}$-categories limits and co-limits (and hence complete and co-complete categories) are going to be different. For this reason we will not bother to call such categories co-complete and we will always talk about complete categories and limits. If it appears in the future that there is indeed two distinct notions of completeness and co-completeness then the terminology of the present paper should be changed. For example, the asymmetry of the relation between the functors $\Hcal i$ and $R_{\Ago}$ in proposition \ref{PropHiwheniishereditary} might be a sign that there is a distinction between left adjoint functors and right adjoint functors to be made in the world of $C^{*}$-category and hence maybe as well between completeness and co-completeness, but we have not been able to make this precise at the present time.}

\subsection{Pre-complete $C^{*}$-categories}

\blockn{The general idea of our definition is that we want to force the category $\Hcal \Cgo$ of small Hilbert modules to be the ``free'' complete $C^{*}$-category on $\Cgo$ in the same way that for ordinary category the category of small pre-sheaves is the free co-complete category generated by a given category. In particular, it should exist a canonical functor $\Hcal \Cgo \rightarrow \Cgo$ for every complete $C^{*}$-category which corresponds to the natural extension of the identity functor from $\Cgo$ to $\Cgo$, and then any (co)limits in $\Cgo$ can be computed by first computing them in $\Hcal \Cgo$ and then applying this functor to $\Cgo$. This motivate the following definition:}

\block{\label{Defpre-comp}\Def{A pre-complete $C^{*}$-category is a couple $(\Cgo,L)$ where $\Cgo$ is a unital $C^{*}$-category and $L : \Hcal \Cgo \rightarrow \Cgo$ is a $C^{*}$-functor endowed with a functorial isometric isomorphism $Id_{\Cgo} \simeq L \circ \Yon$.

If $(\Cgo,L)$ and $(\Dgo,L')$ are two pre-complete $C^{*}$-categories a (unital) $C^{*}$-functor $f : \Cgo \rightarrow \Dgo$ is said to be continuous if the diagram:

\[
\begin{tikzcd}[ampersand replacement=\&]
\Hcal \Cgo \arrow{r}{\Hcal f} \arrow{d}{L} \& \Hcal \Dgo \arrow{d}{L'} \\
\Cgo \arrow{r}{f} \&  \Dgo \\
\end{tikzcd}
\]

commutes up to an isometric isomorphisms which is the canonical\footnote{The one coming from the fact that $f$ is unital.} isomorphism when restricted to the image of the Yoneda embeddings of $\Cgo$ in $\Hcal \Cgo$.

}
It is of course not expected that this notion of pre-complete category is going to be well behaved on itself. The idea is that if one take the same definition for ordinary category it suffices then to add one condition on the couple $(\Cgo,L)$ to indeed obtain a co-complete category with the correct notion of continuous functor. For example, for ordinary categories, it suffices to ask that $L$ is right adjoint to the Yoneda embeddings or that $L$ is itself a continuous functor for the natural pre-complete structure on the category of pre-sheaves. We will do the same for complete $C^{*}$-category: we start from this definition of pre-completeness and we will study the various conditions that we can add.}

\block{One has an analogue of limit in pre-complete $C^{*}$-category: let $\Cgo$ be a pre-complete $C^{*}$-category, let $\Dgo$ be a $C^{*}$-category with $f : \Dgo \rightarrow \Cgo$ a $C^{*}$-functor and $H \in \Hcal \Dgo$ a small Hilbert module. $(\Hcal f) H$ is then an object of $\Hcal \Cgo$,

\Def{In this situation We denote by $f_* H$ the object:

\[ L((\Hcal f) H ) \in \Cgo \]

}

These correspond to the usual concept of a weighted co-limit and will play the role of limit in pre-complete $C^{*}$-categories. $L$ itself is of course the special case which corresponds to $f=id_{\Cgo}$.

As we have seen in \ref{warningNonunital}, if $\Dgo$ is unital and $f$ preserve unit then $f_*(\Yon D) = f(D)$ but for a general $\Dgo$ this does not have to be true.
}

\block{Instead of specifying $L : \Hcal \Cgo \rightarrow \Cgo$ in order to define a pre-complete structure on a $C^{*}$-category $\Cgo$ one can instead specify all the functors $f_{*}: \Hcal \Dgo \rightarrow \Cgo$ when $\Dgo$ is a small $C^{*}$-category and $f:\Dgo \rightarrow \Cgo$ is a functor, one can even restrict to the case where $\Dgo$ and $f$ are unital.

The only requirement (when we restrict to unital categories and unital functors) on such choices of $f_*$ are the functoriality in $\Dgo$ and the fact that $f_*$ is an extension of $f$.

The reason why we can restrict to small categories is because any small Hilbert module $H \in |\Hcal \Cgo|$ can be written as $(\Hcal i) H'$ for $i:\Dgo \rightarrow \Cgo$ the inclusion of a small full subcategory and $H$ a $\Dgo$ Hilbert module. Indeed one simply take $\Dgo$ containing all the generators of $H$ and $H'$ is the restriction of $H$ to $\Dgo$.
}

\block{\label{Consequenceofprecompleteness}Despite what we said above about pre-completeness not being expected to be an interesting notion by itself, it already has some consequences in the $C^{*}$-categorical framework:

\Prop{Let $(\Cgo,L)$ be a pre-complete $C^{*}$-category. Then for any $H \in |\Hcal \Cgo|$ there is a canonical embeddings (i.e. an isometric linear map) of $H$ in $\Yon_{LH}$ and this embeddings is functorial in the sense that for any morphism $h$ in $\Hcal \Cgo$ the following diagram of linear map commutes:

\[
\begin{tikzcd}[ampersand replacement=\&]
H \arrow{r}{h} \arrow[hookrightarrow]{d} \& H' \arrow[hookrightarrow]{d} \\
LH \arrow{r}{Lh} \& LH' \\
\end{tikzcd}
\]

 }

In some sense, the functor $L$ construct a ``completion'' of any $\Cgo$-module into a representable $\Cgo$-modules and extend any operator between two $\Cgo$-modules to an operator between the completion in a canonical and functorial way. 

It also proves that $L$ is injective on morphisms, because an operator $h$ can be recover from $L(h)$ by restricting to the sub-module $H \subset \Yon(LH)$.

\Dem{Let $H \in |\Hcal \Cgo|$. Let $A \in |\Cgo|$ and $a \in H(A)$, the corresponding map $a: \Yon_A \rightarrow H$ give a map $La : A \rightarrow LH$ which in turn corresponds to an element of $\Yon(LH)(A)$ which we denotes by $i(a)$. This construction is linear and functorial from $H$ to $\Yon(LH)$, and it is isometric because if $a \in H(A)$ and $b \in H(B)$ then $\Yon(\scal{a}{b})$ is $a^{*}b$ hence when we apply $L$ one gets that $\scal{a}{b} = L(a^{*})L(b)= \scal{i(a)}{i(b)}$.}

}

\block{\Prop{A pre-complete $C^{*}$-category has splitting of symmetric projections and all bi-products, moreover the inclusion of $A$ into $A \oplus B$ is adjoint to the projection. }

\Dem{This is true in the category $\Hcal \Cgo$ for any $C^{*}$-category $\Cgo$ and these constructions are preserved by any $C^{*}$-functors in particular by $L$, hence they can be computed in $\Hcal \Cgo$ and then pushed to $\Cgo$ by $L$.}
}

\blockn{We will now mention a few key examples of pre-complete $C^{*}$-categories. One should add to this list the category of Hilbert module over a topos that we will present in section \ref{section_Hilbert_bundle}.}

\block{\label{Prop_HC_precomp}\Prop{Let $\Cgo$ be any $C^{*}$-category, then the category $\Hcal \Cgo$ of small Hilbert $\Cgo$-modules is pre-complete for the functor $R_{\Cgo} :\Hcal \Hcal \Cgo \rightarrow \Hcal \Cgo$ which send any $\Hcal \Cgo$ small Hilbert modules to its restriction to the partial two sided ideal $\Cgo \subset_{\Yon} \Hcal \Cgo$. }

\Dem{The only thing we need to prove is that $R_{\Cgo} \circ \Yon$ is naturally (and isometrically) isomorphic to the identity of $\Cgo$. This is essentially the Yoneda lemma: if one considers any Hilbert $\Cgo$-module $H$, then, because $\Hcal \Cgo$ is unital, for any $c \in \Cgo$  morphisms from $\Yon(\Yon(c))$ to $\Yon(H)$ in $\Hcal \Hcal \Cgo$ are the same as morphisms from $\Yon(c)$ to $H$ in $\Hcal(\Cgo)$ and those morphisms $f$ such that $f^{*}f \in \Cgo$ are exactly the compact operators from $\Yon(c)$ to $H$, i.e. the elements of $H(c)$. Hence the restriction of $\Yon(H)$ to $\Cgo$ is indeed (naturally and isometrically) isomorphic to $H$.  }

One could instead specify the $f_*H$ in a very familiar way: indeed if one has $f: \Dgo \rightarrow \Hcal \Cgo$ and $H$ a Hilbert $\Dgo$-module one can define $f_* H$ as the Hilbert $\Cgo$-module $H'$ generated by elements of the form $ h \otimes_d x \in H'(X)$ for $d \in \Dgo$, $h \in H(d)$ and $c \in f(d)(X)$ with the scalar product:

 \[ \scal{ h \otimes_d x}{ h' \otimes_{d'} x'} = \scal{x}{f(\scal{h}{h'}).x' } \]

The positivity relation is easy to check using matrix algebras and sum of Hilbert module, the matrix corresponding to a finite family $(h_i \otimes_{d_i} x_i)$ of such generators can be written as $x^{*} f(h^{*} h) x$ for $h : \oplus d_i \rightarrow H$ and $x: \oplus X_i \rightarrow \oplus f(d_i)$ whose components are given by the $h_i$ and $x_i$ and hence is obviously positive. One can also deduce the positivity from the other description of the pre-complete structure.

This is a (small) generalization of the tensor product by a bimodule: a functor from $\Dgo$ to $\Hcal \Cgo$ can be thought of as a $\Dgo-\Cgo$-bimodule (in the sense of a right Hilbert $\Cgo$-module with a left $*$-action of $\Dgo$).

}

\block{Similarly if $\Ccal$ is a pro-$C^{*}$-algebra as in \cite{phillips1988inverse}, then one can define a $C^{*}$-category of Hilbert modules over $\Ccal$, also denoted $\Hcal \Ccal$: objects are Hilbert $\Ccal$-module as defined in section 4 of \cite{phillips1988inverse} and morphisms are the adjointable operators $f$ (also as in \cite{phillips1988inverse} ) which are bounded in the sense that there exists a constant $K$ such that for each of the natural semi-norm $p$ on the space of adjointable operator, $p(f) \leqslant K$, the norm of $f$ is the smallest such constant $K$. It is proved in \cite{phillips1988inverse} that bounded adjointable operators are dense amongst general adjointable operators.

It has a natural pre-complete structure, indeed if one has $f: \Dgo \rightarrow \Hcal \Ccal$ and $H$ a Hilbert $\Dgo$-module one can define $f_* H$ exactly as for ordinary Hilbert module (with the same generators and scalar product)  with the only difference that we now need to complete with respect to a family of semi-norm instead of just a norm.
}

\blockn{Note that the category of Hilbert modules over a pro-$C^{*}$-algebra $\Ccal$ appears to be full subcategory of the category of Hilbert module over the $C^{*}$-algebra $\Ccal^{b}$ of bounded element of $\Ccal$ and the limits in the two categories are related. This is a special case of the following construction which provide a lot of examples of pre-complete $C^{*}$-category.}

\block{\Def{Let $\Cgo$ be a $C^{*}$-category. 
\begin{itemize}

\item a pre-closure operation $P$ (or $(P,\mu)$) on $\Hcal \Cgo$ is a $C^{*}$-functor $P : \Hcal \Cgo \rightarrow \Hcal \Cgo$ endowed with an isometric inclusion $\mu_X$ of $X$ into $P(X)$ for each $X \in \Hcal \Cgo$ which is functorial in $X$.

\item An object of $\Hcal \Cgo$ is said to be $P$-closed if $\mu_X$ is an isomorphism.

\item $P$ is said to be a closure operation if $P(X)$ is $P$-closed for each $X$. 
\end{itemize}
}
For example, the main result of \cite{paschke1974double} is that if $\Ccal$ is a $C^{*}$-algebra then taking the bi-dual is a closure operation on the category of Hilbert $\Ccal$-modules, closed objects being exactly the reflexive Hilbert modules. It has not been investigated yet whether or not a form of this results hold constructively or not.

If $(\Cgo,L)$ is a pre-complete $C^{*}$-category then $\Yon^{*} \circ L $ is a closure operation on $\Hcal \Cgo$ and closed objects are exactly the representable objects.
}

\block{\label{Prop_preclosedOp}\Prop{Let $\Cgo$ be a pre-complete $C^{*}$-category, $i: \Ago \hookrightarrow \Cgo$ a hereditary subcategory. $R_{\Ago} \circ i_*$ is a pre-closure operation on $\Hcal \Ago$.}

Proposition \ref{PropGeneratorsAndClosOp} and proposition \ref{Prop_quasicontainement}(\ref{Prop_quasicontainement_closureOp}) below give sufficient conditions for this being a closure operation.

\Dem{This is clearly a $C^{*}$-functor. If $X$ is a Hilbert $\Ago$-module and $x \in X(a)$ for some $a \in |\Ago| = |\Cgo|$. There is a map $\tilde{x}$ from $a$ to $i_*(X)$ corresponding to $X$, as $\tilde{x}^{*}x = \scal{x}{x}$, the map $\tilde{x}$ is an element of $R_{\Ago} \circ i_* (X) (a)$, 
This defines a map from $X$ to $R_{\Ago} \circ i_*$ which is an isometric inclusion and functorial in $X$, hence concludes the proof.
}

}

\block{\Prop{ If $P$ is a closure operation of $\Hcal \Cgo$ then the full subcategory $\Dgo$ of $P$-closed object is pre-complete for the structure:

\[ \Hcal \Dgo \overset{\Hcal i}{\rightarrow} \Hcal \Hcal \Cgo \overset{R_{\Cgo}}{\rightarrow} \Hcal \Cgo \overset{P}{\rightarrow} \Dgo \] 
}

Another way to formulate this pre-complete structure is that the $f_*$ are computed in $\Hcal \Cgo$ first and then ``completed'' into objects of $\Dgo$ by applying $P$.

\Dem{The only thing to check is the action of this functor on representable objects. the functor $i: \Dgo \rightarrow \Hcal \Cgo$ is unital so $\Hcal i$ act as $i$ on representable, so a representable in $\Dgo$ is send to itself in $\Hcal \Cgo$ and because it is already closed applying $P$ does not change anything, and this identification is functorial. }
}

\blockn{\Cor{(in classical mathematics) The category of reflexive Hilbert modules over a $C^{*}$-algebra (or over a von Neumann algebra) is a pre-complete $C^{*}$-category. The $f_*H$ are computed by computing them first in the category of Hilbert module and then taking the bi-dual.}

The only reason we need to assume classical logic and axiom of choice here is that we do not know if Paschke result of \cite{paschke1974double} on double dual modules can be formulated and proved in a constructive way.

\Dem{This follows from the fact already mentioned above (see \cite{paschke1974double}) that taking the bidual is a closure operation and reflexive Hilbert modules are exactly the closed objects for this closure operation. }

}

\blockn{The case of module over pro-$C^{*}$-algebra is also a special case of this for the closure operation on the category of Hilbert $\Ccal^{b}$-modules defined by sending any Hilbert $\Ccal^{b}$-module $H$ to the module of bounded vector\footnote{We mean by that element  $h$ of the module such that $\scal{h}{h} \in \Ccal^{b}$, i.e. is a bounded elements in the sense of \cite{phillips1988inverse}.} in the completion of $H$ as a $\Ccal$-module. Closed modules are the exactly those which are the modules of bounded vectors of a Hilbert $\Ccal$-module and morphisms between them corresponds exactly to bounded morphisms between the corresponding Hilbert $\Ccal$-modules.
This has only been proved in classical mathematics but there is apparently no obstruction into making the results constructive without any change, we will not do it here as we do not need this result.
}

\block{\label{Prop_Fctcatarecomplete}\Prop{Let $\Cgo$ be a small $C^{*}$-category, and let $\Dgo$ be a pre-complete $C^{*}$-category. Then the category $\Dgo^{\Cgo}$ of $*$-functors from $\Cgo$ to $\Dgo$ and natural transformations is a $C^{*}$-category endowed with a pre-complete structure which makes the evaluation functor at any $c \in |\Cgo|$ a continuous $*$-functor.}

\Dem{The category of functors $\Dgo^{\Cgo}$ is a $C^{*}$-category for the norm on natural transformation $\Vert \eta \Vert =\sup_{c \in \Cgo} \Vert \eta_C \Vert$ and the $*$-operation $(\eta^{*})_c = (\eta_c) ^{*}$. For each $c \in |\Cgo|$ the evaluation at $c$ induce a $*$-functor $ev_c : \Dgo^{\Cgo} \rightarrow \Dgo$ and for each $f: c \rightarrow c'$ one has a natural transformation $ev_f$ from $ev_c$ to $ev_{c'}$.

We want to construct a functor $L : \Hcal(\Dgo^{\Cgo}) \rightarrow \Dgo^{\Cgo}$, and the general idea is that ``limit should be computed objectwise''. So for each $c \in |\Cgo|$ one has the evaluation functor $ev_c: \Dgo^{\Cgo} \rightarrow \Dgo$ which induces a functor $\Hcal(ev_c) : \Hcal(\Dgo^{\Cgo}) \rightarrow \Hcal \Dgo$ which can be in turn composed with the complete structure of $\Dgo$ to obtain a functor: $L_c : \Hcal(\Dgo^{\Cgo}) \rightarrow \Dgo$, one can then see that a morphism $f: c\rightarrow c'$ induces a natural transformation from $L_c$ to $Lc'$ and that this turn the family $L_c$ into a functor $L : \Hcal(\Dgo^{\Cgo}) \rightarrow \Dgo^{\Cgo}$.

One then easily check that $L$ is the identity on representable (because the $ev_c$ are unital, $\Hcal(ev_c)$ act as $ev_c$ on the representable) and that the evaluation functor are continuous essentially by construction.
}
}

\block{We will now try\footnote{Not completely successfully.} to formulate a universal property that the $f_* H$ should satisfies (as we mentioned earlier they are the construction that we want to think of as weighted (co)limits). Here is the corresponding notion of cone:

\Def{Let $f:\Dgo \rightarrow \Cgo$ be a $C^{*}$-functor and $H$ a Hilbert $\Dgo$-module, then a cone for $(f,H)$ is an object $X$ of $\Cgo$ with for each $D \in |\Dgo|$ and $h \in H(D)$ a morphism $X_h: f(D) \rightarrow X $ such that:

 \[ X_h^{*} X_{h'} = f(\scal{h}{h'}) \]
 
\[ X_h \circ f(d) = X_{hd} \]

\[ X_h + X_{h'}= X_{h+h'} \]
}

The typical example is $f_*(H)$ itself, which comes with such a structure of cone for $(f,H)$.
Unfortunately we have not been able to find any purely $C^{*}$-categorical formulation of the fact that $f_* H$ should be the ``universal cone'' for $(f,H)$, proposition \ref{limitanddensity} below is the closest approximation to such a result we have been able to found.
}

\block{\label{PropCompactActingOnCones}\Prop{Let $f : \Dgo \rightarrow \Cgo$ be a $C^{*}$-functor, $H$ a Hilbert $\Dgo$-module and $X$ be a cone for $(f,H)$, with $X_h:f(D) \rightarrow X$ the structural map for $h \in H(D)$. Then:

\begin{enumerate}
\item The sub-$C^{*}$-algebra of $\Cgo(X)$ generated by the endomorphisms of the form $(X_h)(X_{h'})^{*}$ is (naturally) a quotient of the algebra of compact operators of the Hilbert $\Dgo$-module $H$. If $f$ is faithful then this quotient map is an isomorphism

\item The sub-$C^{*}$-algebra of $\Cgo(X)$ generated by the endomorphisms of the form $(X_h) g (X_{h'})^{*}$ for $h \in H(D), h' \in H(D')$ and  $g \in \Cgo(f(D'),f(D))$ is naturally isomorphic to the algebra of compact operators of the Hilbert $\Cgo$-module $(\Hcal f)(H)$. 

\item The second algebra is also the smallest hereditary sub-algebra of $\Cgo(X)$ which contains the first algebra.

\end{enumerate}
}

\Dem{\begin{enumerate}

\item The algebra of compact operators of the Hilbert $\Dgo$-module $H$ is generated by operators of the form $hh'^{*}$, moreover the operators of the form $(X_h) (X_h')^{*}$ compose in the following way:

\[\left(  X_{h_1} X_{h'_1}^{*} \right) \left(  X_{h_2} X_{h'_2}^{*} \right)  = X_{h_1} f(\scal{h'_1}{h_2}) X_{h'_2}^{*} = X_{h_2 (h'_1)^{*} h_2} X_{h'_2}^{*} \]

which is exactly the way the compact operators $hh'^{*}$ compose. In order to conclude we need to show that for any finite family $h_i,h'_i \in H(D_i)$ (with $i=1 ,\dots, n$) one has :

\[ \left\Vert \sum_{i=1}^{n} X_{h_i} X_{h'_i}^{*} \right\Vert \leqslant \left\Vert \sum_{i=1}^{n} h_i (h'_i)^{*} \right\Vert \]

with equality when $f$ is faithful. Let $A$ and $B$ be those two operators, and assume first that $A$ and $B$ are both self-adjoint. Then let $\tau$ be the morphism: $(X_{h_1},\dots, X_{h_n}) : \bigoplus f(D_i) \rightarrow X$, (if $\Cgo$ do not have direct sum, this exist as a morphism in $\Hcal \Cgo$), let $\tau'$ be the similar morphism with the $X_{h'_i}$ and let $\eta$ and $\eta'$ be the similar morphisms in $\Hcal \Dgo$ with the $h_i$ and $h'_i$ instead. One has, using the properties of the spectral radius (see \ref{Prop_specctral_radius}):

\begin{multline*} \Vert A \Vert= \rho(A) = \rho (\tau \tau'^{*}) = \rho(\tau'^{*} \tau) = \rho(f(\scal{h'_j}{h_i})_{i,j} ) = \rho(f(\eta'^{*} \eta)) \\ \leqslant \rho(\eta'^{*}\eta) = \rho(\eta \eta'^{*})= \rho(B)=\Vert B \Vert \end{multline*}

And there is equality if $f$ is faithful.

In the general case, when $A$ and $B$ are not self-adjoint, it suffices to use that $\Vert A \Vert =\Vert A^{*} A \Vert ^{1/2}$ and $\Vert B \Vert = \Vert B^{*} B \Vert ^{1/2}$, and using the computation rules shared by the $X_{h_i} X_{h'^{*}_i}$ and the $h_i h'^{*}_i$ one can write $A^{*}A$ and $B^{*}B$ in a form similar to the first case and concludes the proof of the first point.

\item For any $Y$, the space $(\Hcal f)(H)(Y)$ is generated by elements of the form $hu$ for $h \in H(D)$ and $u \in \Cgo(Y,f(D))$. Hence the algebra of compact operators of $(\Hcal f) (H)$ is generated by elements of the form $huu'^{*}h'^{*} = h g h^{*}$ where $g=uu'^{*}$ is an element of $\Cgo(f(D'),f(D))$. The rest of the proof is then essentially the same as the previous point.

\item If $A$ is an hereditary sub-algebra which contains all the $X_h X_{h'}^{*}$ then it also contains the $X_h g X_{h'}^{*}$ because (by lemma \ref{Lem_hered}(\ref{Lem_hered3})):
\[ X_h g X_{h'}^{*} = \lim (X_h X_h^{*})^{1/n} X_h g X_{h'}^{*} (X_{h'} X_{h'}^{*})^{1/n} \in A \] 

Conversely, the sub-algebra of $\Cgo(X)$ generated by the $X_h g X_{h'}^{*}$ is hereditary because for any endomorphism $f \in \Cgo(X)$ one has:

\[  \left( X_{h_1} g_1 X_{h_1'}^{*} \right) f \left( X_{h_2} g_2 X_{h_2'}^{*} \right) = X_{h_1} \left( g_1 X_{h_1'}^{*} f X_{h_2} g_2 \right) X_{h_2'}^{*}\]

and $g_1 X_{h_1'}^{*} f X_{h_2} g_2$ is an element of $\Cgo(f(D'_2),f(D_1))$ hence this concludes the proof.

\end{enumerate}
}

}

\block{\label{def_weakly_dense}\Def{Let $\Cgo$ be a pre-complete $C^{*}$-category, $X \in |\Cgo|$ any object and $\Acal \subset \Cgo(X)$ a hereditary sub-algebra, one says that $\Acal$ is weakly dense if $L(\Acal X)$ is isomorphic to $X$ with an isomorphism compatible with the inclusion of $\Acal X$ in $L(\Acal X)$ and in $X$.
}

$\Acal X$ denotes the submodule of $\Ygo_X$ which is the image of the action of $\Acal$ as it is defined in \ref{Def_AH}.

Without additional assumptions on $\Cgo$ it does not seem possible to prove that the isomorphism between $L(\Acal X)$ and $X$ is unique or canonical. This will be the role of the axiom $(C1)$ introduced in the next sub-section.
}

\block{What weak density means in the various examples of pre-complete $C^{*}$-categories can be found in:

\begin{itemize}
\item Proposition \ref{Prop_weakdensityinHC} for categories of Hilbert modules and this essentially corresponds to the conditions of proposition \ref{densityforHmodule}.
\item Proposition \ref{Prop_functcat_C1} for categoties of functors from a small $C^{*}$-category to a nice $C^{*}$-category.

\item For Hilbert modules over a Pro-$C^{*}$-algebra and reflexive Hilbert module over a $C^{*}$-algebra we did not treat it explicitly but one can use their description as closed module for a closure operations. In this case $\Acal \subset \Cgo(H)$ a hereditary subalgebra where $H$ is some closed module is weakly dense if and only the module $H' = \Acal H$ has $H$ for closure (with the inclusion of $H'$ into $H$ being the natural one). So for example in the case of a pro-$C^{*}$-algebra weak density corresponds to density for pointwise convergence in the pro-topology.

\end{itemize}

}

\block{\label{limitanddensity}\Prop{Let $f : \Dgo \rightarrow \Cgo$ and $H$ as before. Let $X$ be a cone for $(f,H)$,  then the following conditions are equivalent:
\begin{itemize}

\item There is an isomorphism between $X$ and $f_* H$ compatible with the cone structure. 

\item The image in $\Cgo(X)$ of the algebra of compact operators $\Kcal \Cgo ((\Hcal f)(H))$ is a weakly dense hereditary sub-algebra.

\item The hereditary sub-algebra spammed by the $X_h X_h'^{*}$ in $\Cgo(X)$ is weakly dense.

\end{itemize}
}
The first condition somehow mean that $X$ is ``the limit'' of $(f,H)$ with the slightly annoying detail that the isomorphism between $X$ and $f_* H$ is in general non unique. This will be fixed by the addition of the condition $(C1)$ of the next subsection. Of course as the notion of weak density is not purely $C^{*}$-categorical but relies on the pre-complete structure this does not provide an abstract characterization of limits in $C^{*}$-categories. But in a lot of examples, weak density has a simple description and hence this result do act as a universal properties in these cases.

\Dem{ The equivalence of the second and the third point is given by the proposition \ref{PropCompactActingOnCones} above. Let $\Kcal$ denote the algebra $\Kcal \Cgo ((\Hcal f)(H))$, it acts on $X$ by proposition \ref{PropCompactActingOnCones}. We will first construct an isomorphism of Hilbert $\Cgo$-modules between $\Kcal .X$ and $(\Hcal f)(H)$. For any generator of $\Kcal$ of the form $X_h g X_h'$ and any element $x$ of $X(c)$, $X_h g X_h' x$ is of the form $X_h u$ with $h \in H(D)$ and $u \in \Cgo(f(D),C)$ such an element can be identified with the generators $h.u$ of $(\Hcal f)(H)$ and the scalar products computed between element of the form $h.u$ or the corresponding elements of the form $X_h u$ yields the same results. Hence, as such elements are generators of $\Kcal X$ this induces an isomorphism between $\Kcal X$ and $(\Hcal f ) (H)$.

Now the weak density of $\Kcal$ in $\Cgo(X)$ mean that $L(\Kcal X)$ is isomorphic to $X$ in a way compatible to the inclusion of $\Kcal X$, which through the previously constructed isomorphism is exactly the same as saying that $f_* H = L((\Hcal f) H)$ is isomorphic to $X$ in a way compatible to the cone structure, which proves the equivalence of the first point with the second.

}
}

\subsection{The separation condition $(C1)$}

\block{\Def{A pre-complete category $(\Cgo,L)$ is said to satisfy condition $(C1)$ or the separation condition if:
\begin{itemize}
\item[(C1)] For all $H \in |\Hcal \Cgo|$ and $\lambda: LH \rightarrow A$ a map in $\Cgo$, if the restriction of $\lambda$ to $H \subset \Yon(LH)$ is zero then $\lambda = 0$.
\end{itemize}

}
This express the idea that as $LH$ is in some sense the completion of $H$, or the colimit of a diagram corresponding to $H$, then functions on $LH$ should be described by their values on $H$.

Under this condition $(C1)$ the isomorphisms of definition \ref{def_weakly_dense} and of proposition \ref{limitanddensity} are unique.

By proposition \ref{limitanddensity}, the couples $(H,LH)$ for $H \in \Hcal \Cgo$ are exactly the couple $(\Acal.X,X)$ for $X \in |\Cgo|$ and $\Acal \subset \Cgo(X)$ a weakly dense hereditary sub-algebra of $\Cgo(X)$. Hence one can immediately deduce that: 

}

\block{\Prop{Let $\Cgo$ be a pre-complete category, then $\Cgo$ satisfies $(C1)$ if and only if for any object $X \in |\Cgo|$, any weakly dense hereditary sub-algebra $\Acal \subset \Cgo(X)$ is essential, i.e. if $f\in \Cgo(X)$ satisfies $fa=0$ for all $a \in \Acal$ then $f=0$.}

\Dem{By propositions \ref{Prop_H'=AH} and \ref{limitanddensity} Any $H \in \Hcal\Cgo$ can be written in the form $\Acal.X$ with $X = LH$ and $\Acal$ weakly dense in $\Cgo(X)$. And for any map $\lambda$ from $X$ to another object $Y \in |\Cgo|$ saying that $\lambda|_H =0$ is the same as saying that for all $a \in \Acal$, $\lambda^{*} \lambda a =0$.}

}

\block{\label{Prop_weakdensityinHC}\Prop{Let $\Cgo$ be a $C^{*}$-category, the category $\Hcal \Cgo$ is endowed with its pre-complete structure of \ref{Prop_HC_precomp}, a hereditary sub-algebra $\Acal \subset \Hcal\Cgo(X)$ is weakly dense if and only if it satisfies the equivalent conditions of proposition \ref{densityforHmodule}, for example if it contains the compact operators of $X$. In particular, they are essential and $\Hcal \Cgo$ satisfies $(C1)$.}

\Dem{Let $X$ be a Hilbert $\Cgo$-module and let $\Acal \subset \Hcal\Cgo(X)$ be a hereditary sub-algebra the module $L(\Acal.X)$ is obtained by computing $\Acal X$ as a $\Hcal\Cgo$-Hilbert module and then restricting to $\Cgo$, hence it is exactly the module $\{a \in X | aa^{*} \in \Acal \}$, hence it is the natural notion of $\Acal X$. It is isomorphic to $X$ with the correct cone structure if and only the natural inclusion map $\Acal .X \rightarrow X$ is an isomorphism which was one of the equivalent condition of proposition \ref{densityforHmodule}. In particular an operator $f$ such that $f a = 0$ for all $a \in \Acal$ is automatically equal to $0$ on $\Acal.X =X$ hence $\Acal$ is essential which proves that $\Cgo$ satisfies $(C1)$.}

}

\block{\label{Prop_PC1}\Prop{Let $\Cgo$ be a $C^{*}$-category and $P$ a closure operation on $\Hcal \Cgo$. Then the pre-complete $C^{*}$-category $\Dgo$ of $P$-closed Hilbert $\Cgo$-modules as defined in \ref{Prop_HC_precomp} satisfies $(C1)$ if and only $P$ satisfies the following condition:

\begin{itemize}
\item [(PC1)] For all $X \in |\Hcal\Cgo|$, $X \overset{\mu_X}{\hookrightarrow} P(X)$ is an essential sub-module.
\end{itemize}

}

We recall that a sub-module $X \subset Y$ is essential if for map $h: Y \rightarrow Z$  if $h$ is zero on $X$ then $h=0$. By considering $h^{*}h$ one can see that it is enough to test it for $h$ an endomorphism of $Y$.

\Dem{ Let $\Dgo$ be the category of $P$-closed object and $i:\Dgo \rightarrow \Hcal \Cgo$ the canonical inclusion. For any Hilbert $\Dgo$-module $H$, $LH$ is computed as $P (i_*H)$ if $f: LH \rightarrow W$ is any map, to say that for all $h \in H$ one has $f(h)=0$ is equivalent to say that $f$ restricted to $i_* H \subset P(i_* H)$ is $0$, hence if for all $H$, $i_* H $ is essential in $P(i_* H)$ this proves that $\Dgo$ satisfies $(C1)$.

Conversely assume that the category $\Dgo$ of $P$-closed objects satisfies $(C1)$. Let $X \in |\Hcal \Cgo|$, $X$ is a sub-module of $PX$, hence by proposition \ref{Prop_H'=AH} $X$ is of the form $\Acal.(PX)$ for some hereditary sub-algebra $\Acal \subset \Cgo(X)$. By definition of the pre-complete structure on $\Dgo$, $\Acal$ is weakly dense (indeed $\Acal.P(X)$ in $\Dgo$ is $P(\Acal.P(X)) = P(X)$ in a way compatible to the cone structure), hence as $\Dgo$ satisfies $(C1)$ this implies that $\Acal$ is an essential hereditary sub-algebra of $\Dgo(PX) = \Cgo(PX)$. Hence of any $h \in \Cgo(PX)$ such that for all $x \in X$, $h(x)=0$ one has $\forall a \in \Acal, ha =0$ hence $h=0$ which concludes the proof.
}
}

\block{\label{continuous=preservedensity}\Prop{Let $\Cgo$ and $\Dgo$ be two pre-complete $C^{*}$-categories and let $f :\Cgo \rightarrow \Dgo$ be a $C^{*}$-functor and consider the following two conditions:

\begin{enumerate}

\item $f$ is continuous functor (i.e. there exists an isometric isomorphism as in definition \ref{Defpre-comp})

\item For every $X \in |\Cgo|$ and every $\Acal \subset \Cgo(X)$ a weakly dense hereditary sub-algebra, the hereditary sub-algebra spammed by $f(\Acal)$ is weakly dense.

\end{enumerate}

Then $(1) \Rightarrow (2)$ and if $\Dgo$ satisfies $(C1)$ they are equivalent.
}

\Dem{The key point is the observation that if $\Acal \subset \Cgo(X)$ is an hereditary sub-algebra and if $\Bcal$ is the hereditary sub-algebra spammed by $f(\Acal)$ then there is a canonical isomorphism between $\Hcal(f)(\Acal X)$ and $\Bcal f(X)$ in $\Hcal \Dgo$: indeed they are both the module $T$ generated element $a \in T(f(X))$ for each $a \in \Acal$ such that $\scal{a}{a'} = f(a^{*}a')$.

We denote by $L$ and $L'$ the structural functor $\Hcal \Cgo \rightarrow \Cgo$ and $\Hcal \Dgo \rightarrow \Dgo$.

Assume the first condition, then if $\Acal$ is dense the isomorphism between $L(\Acal X)$ and $X$ inducing the identity on $\Acal X$ is sent by $\Hcal f$ on an isomorphism between $L'(\Bcal f(X))$ and $f(X)$ inducing the identity on $\Hcal f(\Acal X) = \Bcal f(X)$, which proves that $\Bcal$ is weakly dense in $\Dgo(f(X))$.

Assume now the second condition and let $H \in |\Hcal \Cgo|$ a Hilbert module. As a sub-module of $\Yon LH$, there exists a hereditary sub-algebra $\Acal_H \subset \Cgo(LH)$ such that $\Acal_H LH = H$. In particular, $\Acal_H$ is weakly dense, and hence by assumption its image by $f$ is weakly dense, i.e. there exists an isometric isomorphism $\mu_H$ between $L'((\Hcal f) (H))$ and $f(LH)$ which induces the identity on $(\Hcal f ) H$.
We need to prove the functoriality of these isomorphisms, i.e. for any $: H \rightarrow H'$ an arrow in $\Hcal \Cgo$, one need to prove the commutativity of the diagram:

\[
\begin{tikzcd}[ampersand replacement=\&]
L'((\Hcal f)(H)) \arrow{r}{\mu_H} \arrow{d}{L'((\Hcal f)(h))} \& f(LH) \arrow{d}{f(Lh)} \\
L'((\Hcal f)(H')) \arrow{r}{\mu_{H'}} \&  f(LH') \\
\end{tikzcd}
\]

But the commutativity is clear when we restrict to the sub-module $(\Hcal f)(H) \subset L'((\Hcal f)(H))$, hence follows from $(C1)$. It is also clear that $\mu_H$ is the canonical isomorphism when $H$ is representable.
}
}

\block{\label{Prop_functcat_C1}\Prop{Let $\Dgo$ be a pre-complete $C^{*}$-category satisfying $(C1)$, and $\Cgo$ a small $C^{*}$-category. Then:

\begin{itemize}
\item If $F : \Cgo \rightarrow \Dgo$ is a $*$-functor and $\Acal$ is a hereditary subalgebra of endomorphisms of $F$ then $\Acal$ is weakly dense if and only if for all $c \in \Cgo$ the action of $\Acal$ on $F(c)$ generate a weakly dense hereditary sub-algebra in $\Dgo(F(c))$.

\item In particular the category $\Dgo^{\Cgo}$ of $*$-functor also satisfies $(C1)$.
\end{itemize}
}

\Dem{By proposition \ref{Prop_Fctcatarecomplete} the functors $ev_c: F \mapsto F(c)$ are continuous, and by proposition \ref{continuous=preservedensity} this implies that if $\Acal$ is weakly dense on $F$ then its image on each $F(c)$ generates a weakly dense hereditary sub-algebra. The continuity of $ev_c$ also implies that $(\Acal F)(c) = \Acal_c F(c)$ where $\Acal_c$ denotes the hereditary sub-algebra of $F(c)$ generated by $\Acal$. If each $\Acal_c$ is weakly dense one has that $(\Acal F)(c) \simeq F(c)$ but $(C1)$ in $\Dgo$ easily implies that those isomorphisms are functorial and induce an isomorphism $\Acal F \simeq F$ which preserves the cone structure which proves the first point. In particular if $\Acal$ is weakly dense and $\mu :F \rightarrow G$ is a natural transformation such that $fa = 0$ for all $a \in \Acal$ then $\mu_c a = 0$ for each $a \in \Acal_c$ hence $\mu_c=0$ as $\Acal_c$ is weakly dense hence $\mu=0$.}

}

\block{\label{Prop_weaklyNDg}\Prop{Let $\Dgo$ be a (non unital) $C^{*}$-category and $\Cgo$ a pre-complete $C^{*}$-category. Let $f : \Dgo \rightarrow \Cgo$, consider the following conditions:

\begin{enumerate}

\item There is an isometric isomorphism of functors $\mu$ between $f$ and $f_* \circ \Yon : \Dgo \rightarrow \Hcal \Dgo \Rightarrow \Cgo$, such that for any $h \in \Dgo(d',d)$ the map: 

\[f(d') \overset{f(h)}{\rightarrow} f(d) \overset{\mu_d}{\rightarrow} f_*(\Yon(d)) \]

Is the map $f(d') \rightarrow f_*(\Yon(d))$ corresponding to $h$ in the cone structure on $f_*(\Yon(d))$.

\item For all $d \in |\Dgo|$, $f(\Dgo(d))$ generates a weakly dense hereditary sub-algebra of $\Cgo(f(d))$.

\end{enumerate}

Then $(1) \Rightarrow (2)$ and if $\Cgo$-satisfies $(C1)$ they are equivalent.

}

Assuming $\Cgo$ satisfies $(C1)$ a functor satisfying those conditions will be called \emph{weakly non-degenerate}. One can see that a functor $f : \Dgo \rightarrow \Dgo'$ between two $C^{*}$-categories is non-degenerate in the sense of the remark made in \ref{warningNonunital} if and only if it is weakly non-degenerate when seen as a functor with values in $\Hcal \Dgo'$.

\Dem{By proposition \ref{limitanddensity} the second condition is equivalent to the fact that for each $d \in |\Dgo|$ there is an isomorphism $\mu_d$ between $f(d)$ and $f_*(\Yon(d))$ which is compatible with ``the cone structures'', i.e. such that for any element in $\Dgo(d',d)$ the natural morphisms from $f(d')$ to $f(d)$ and from $f(d')$ to $f_*(\Yon(d))$ commute with the isomorphism $\mu_d$. The first condition clearly gives a family of such $\mu_d$. Conversely, assume that one has such a family of $\mu_d$, and we will show that assuming $\Cgo$ satisfies $(C1)$ they form a natural transformation satisfying the condition of the first point, i.e. that for any $h \in \Dgo(d',d')$ the following diagram commutes:

\[
\begin{tikzcd}[ampersand replacement=\&]
f(d') \arrow{r}{f(h)} \arrow{d}{\mu_{d'}} \arrow{dr} \& f(d) \arrow{d}{\mu_d} \\
f_* \Yon(d') \arrow{r}{f_* \Yon(h)} \& f_* \Yon(d)\\
\end{tikzcd}
\]

Where the diagonal map is the cone structure. Now the upper triangle commute by definition of the $\mu_d$, and the lower triangle commute as soon as we pre-compose the two sides by any arrow in $f(\Dgo(d'))$, as by assumption they spam a weakly dense hereditary sub-algebra this implies by $(C1)$ that the lower triangle also commutes and concludes the proof.

}

}

\block{\label{Prop_restrictionoftransnat1}\Prop{Let $\Dgo$ be a (preferably small) $C^{*}$-category and $\Cgo$ a pre-complete $C^{*}$-category satisfying $(C1)$. Then for any two weakly non-degenerate functors\footnote{i.e. satisfying the equivalent conditions of the previous proposition} $i,j : \Dgo \rightarrow \Cgo$ there is a bijection $Hom(i,j) \simeq Hom(i_*,j_*)$ given by the functoriality of $i \rightarrow i_*$ in one direction and the restriction to $\Dgo$ in the other.
}

If we do not assume that $\Dgo$ is small, the proposition is still true but no longer makes sense in our framework: it has to be interpreted as a meta-theorem or require stronger assumption on the framework, like a class of function between class or a notion of $2$-class (with a $2$-class of all class).

\Dem{By the previous proposition the restrictions of $i_*$ and $j_*$ to $\Dgo$ are naturally isomorphic to $i$ and $j$ hence one can indeed restrict a natural transformation between $i_*$ and $j_*$ to one between $i$ and $j$. Conversely if $\alpha: i \rightarrow j$ is a natural transformation then there is a natural transformation $\alpha' :i_* \rightarrow j_*$ such that for any Hilbert module $H \in |\Hcal \Dgo|$ and any $h \in H(d), d \in |\Dgo|$ the following diagram commutes, the vertical arrows being the cone structures:

\[
\begin{tikzcd}[ampersand replacement=\&]
i(d) \arrow{r}{\alpha_d} \arrow{d}{h} \& j(d) \arrow{d}{h} \\
i_*H \arrow{r}{\alpha'_H} \& j_* H\\
\end{tikzcd}
\]

By $(C1)$, $\alpha'$ is uniquely determined by these diagrams (indeed, one know its value on every elements of the cone structure of $i_* H$), and by the previous proposition, $\alpha_d : i(d) \simeq i_* \Yon d \rightarrow j(d) \simeq j_* \Yon d$ also satisfies this condition hence up to those natural isomorphisms is equal to $\alpha'_{\Yon d}$ which concludes the proof.
 }

}

\subsection{The completeness conditions $(C2)$,$(C3)$ and $(C4)$ }

\block{\label{def_extension_prop}
\Def{Let $\Cgo$ be a $C^{*}$-category and $P$ be a closure operation on $\Hcal \Cgo$. A Pre-adjoint pair for $P$ in $\Cgo$ is the data of:
\begin{itemize}
\item Two Hilbert $\Cgo$-modules $A,B \in |\Hcal\Cgo|$.

\item Two bounded linear maps $f: A \rightarrow PB$ and $g: B \rightarrow PA$.

\item For all $X,Y \in |\Cgo|$ and for all $a \in A(X)$, $b \in B(Y)$ one has:

\[ \scal{f(a)}{\mu_B(b)} = \scal{\mu_A(a)}{g(b)} \]

\end{itemize}

An extension of a pre-adjoint pair is then an operator $\tilde{f}$ from $PA$ to $PB$ which is an extension of $f$, i.e. for all $X\in \Cgo$, $a \in A(X)$ $f(a)=\tilde{f}(a)$ in $PB(X)$.

}
}

\block{\label{Prop_uniqunessinextprop}\Prop{Let $\Cgo$ be a $C^{*}$-category and $P$ a closure operation on $\Hcal \Cgo$ satisfying\footnote{Or equivalently, such that the category of $P$-closed object satisfies $(C1)$} $(PC1)$ then the extension $\tilde{f}$ of a pre-adjoint pair is unique and $\tilde{f}^{*}$ is an extension of $g$.}

\Dem{The uniqueness follows directly from $(PC1)$: if $f'$ is any other extension of $f$, the operator $(\tilde{f}-f')$ is zero on $A$ hence also zero on $PA$.

Let $b \in B(X)$ for any $X \in \Cgo$, $\tilde{f}^{*}(b)$ and $g(b)$ are both elements of $(PA)(X)$ i.e. maps from $X$ to $PA$ such that for all $a \in A(Y)$,

\[\scal{a}{g(b)} = \scal{a}{\tilde{f}^{*}(b)} = \scal{f(a)}{b} \in \Cgo(X,Y) \]

i.e. $(g(b))^{*}$ and $(\tilde{f}^{*}(b))^{*}$ have the same restriction to $A$, hence by $(C1)$ are equals on $PA$ and hence $g(b) = \tilde{f}^{*}(b)$ which concludes the proof.

}

}

\block{\Def{Let $(\Cgo,L)$ be a pre-complete $C^{*}$-category. One consider the following three conditions on $\Cgo$:

\begin{itemize}

\item[(C2)] The functor $L : \Hcal \Cgo \rightarrow \Cgo$ is continuous in the sense of \ref{Defpre-comp}.

\item[(C3)] For any $X \in |\Cgo|$ and any pair of hereditary sub-algebras $\Acal \subset \Bcal \subset \Cgo(X)$ if $\Acal$ is weakly dense then $\Bcal$ is weakly dense.

\item[(C4)] Every pre-adjoint pair for the closure operation $L\circ \Yon$ on $\Hcal \Cgo$ has an extension.

\end{itemize}

}}

\block{\label{Prop_impbetween_Ci}\Prop{Let $(\Cgo,L)$ be a pre-complete $C^{*}$-category then:

\begin{itemize}

\item If $\Cgo$ satisfies $(C3)$ and $(C1)$ then $\Cgo$ satisfies $(C2)$.

\item If $\Cgo$ satisfies $(C4)$ and $(C1)$ then $\Cgo$ satisfies $(C3)$ and hence also $(C2)$.

\end{itemize}}

We have not been able to find counterexamples to the other implications between the conditions $(Ci)$, nor to find examples of pre-complete $C^{*}$-categories that does satisfies those conditions and we would be very interested by such examples.

\Dem{

Assume that $(\Cgo,L)$ satisfies $(C1)$ and $(C3)$. By proposition \ref{continuous=preservedensity}, in order to prove that $\Cgo$ satisfies $(C2)$, i.e. that $L:\Hcal \Cgo \rightarrow \Cgo$ is continuous, it suffice, by proposition \ref{continuous=preservedensity}, to show that it send weakly dense hereditary sub-algebras to weakly dense hereditary sub-algebras. So let $x \in |\Hcal \Cgo|$ and $\Acal \subset \Hcal\Cgo(X)$ an hereditary subalgebra which is weakly dense for the natural pre-complete structure of $\Hcal \Cgo$. We proved in \ref{Prop_weakdensityinHC} that this mean that $\Acal$ contains the algebra of compact operators of $X$, hence its image in $LX$ contain the image under $L$ of the compact operators of $X$ which spam a weakly dense hereditary algebra because of proposition \ref{limitanddensity} hence the image of $\Acal$ under $L$ also spam a weakly dense hereditary sub-algebra because of $(C3)$.

\bigskip

Let now $(\Cgo,L)$ be a pre-complete $C^{*}$-category satisfying $(C1)$ and $(C4)$, let $X \in |\Cgo|$ be an object and let $\Acal \subset \Bcal \subset \Cgo(X)$ be two hereditary sub-algebras such that $\Acal$ is weakly dense, one need to prove that $\Bcal$ is also weakly dense.

One has, in $\Hcal\Cgo$, $\Acal X \subset \Bcal X \subset \Yon X$ and there exists an isomorphism between $L(\Acal X)$ and $X$ which induces the identity on $\Acal X$.

Hence, there exists a linear map $f$ from $\Acal X$ to $L(\Bcal X)$ given by $\Acal X \hookrightarrow \Bcal X \hookrightarrow L(\Bcal X)$ and a map $g$ from $\Bcal X$ to $L(\Acal X)$ given by $\Bcal X \hookrightarrow X \simeq L(\Acal X)$. For any $a \in \Acal X$ , $b \in \Bcal X$ one has:

\[ \scal{f(a)}{b} = \scal{a}{g(b)} \]

and it is the scalar product in $X$. Hence there exists a map $u$ from $L(\Bcal X)$ to $L(\Acal X) \simeq X$ which extend the inclusion of $\Bcal X$ into $X$ and whose adjoint induces the map $\Acal X \hookrightarrow \Bcal X \hookrightarrow L(\Bcal X)$. In particular $uu^{*}$ induces the identity on $\Acal X$ and hence on $X$ by $(C1)$ and $u^{*}u$ induces the identity on $\Bcal X$ and hence on $L(\Bcal X)$ by $(C1)$.

}

}

\block{\label{Prop_HC_Ci}\Prop{Let $\Cgo$ be a $C^{*}$-category then $\Hcal \Cgo$ satisfies $(C1)$,$(C2)$,$(C3)$ and $(C4)$.}

\Dem{$(C1)$ has been proved in \ref{Prop_weakdensityinHC} hence by \ref{Prop_impbetween_Ci} it is enough to prove that $\Hcal \Cgo$ satisfies $(C4)$. This could be seen as a direct corollary of the implication $(PC4) \Rightarrow (C4)$ proved in the next propostion with $P$ the identity functor on $\Hcal \Cgo$, or as well as a corollary of proposition \ref{Prop_PC5impPCi}, but we prefer the following direct proof:

Let $A$,$B$ and $f:A \rightarrow L(B)$ $g: B \rightarrow L(A)$ be a pre-adjoint pair in $\Hcal \Hcal \Cgo$. As $L$ is just the restriction functor $R_{\Cgo}: \Hcal \Hcal \Cgo \rightarrow \Hcal \Cgo$, $f$ and $g$ after restriction just defines a pair of adjoint linear map between the restriction of $A$ and the restriction of $B$, hence an operator $\tilde{f}:LA \rightarrow LB$ which extend the restriction of $f$ to a weakly dense sub-module (the restriction of $A$ to $\Cgo$, hence extend $f$ (because $\Hcal \Cgo$ satisfies $(C1)$).
}
}

\block{\Prop{Let $\Dgo$ be a pre-complete category satisfying $(C1)$ and $(C4)$ and $\Cgo$ be a small $C^{*}$-category, then $\Dgo^{\Cgo}$ satisfies $(C1)$ and $(C4)$ (and in particular $(C2)$ and $(C3)$).}

\Dem{The fact that $\Dgo^{\Cgo}$ satisfies $(C1)$ has been proved in \ref{Prop_functcat_C1}, so we only have to prove $(C4)$.

Let $f:A \rightarrow LB$ and $g:B \rightarrow LA$ be a pre-adjoint pair in $\Hcal(\Dgo^{\Cgo})$. For each $c \in \Cgo$, the functor $\Hcal ev_c : \Hcal(\Dgo^{\Cgo}) \rightarrow \Hcal \Dgo$ acts on bounded linear map and hence one can transport this pair into a similar pair in $\Hcal \Dgo$, as $\Dgo$ satisfies $(C4)$ there will be an extension of $ev_c(f)$ in $\Dgo$ for each $c$, and using $(C1)$ one can prove that all those maps provides a natural transformation and hence are a morphism $\tilde{f}$ in $\Dgo^{\Cgo}$ that extend $f$.} }

\blockn{We will now give a formulation of each of the conditions $(C2),(C3),(C4)$ in terms of a closure operation defining the corresponding category. The corresponding statement for $(C1)$ is \ref{Prop_PC1} }

\block{\label{Prop_RelCiPCi}\Prop{Let $\Cgo$ be a $C^{*}$-category, $P$ be a closure operation on $\Hcal \Cgo$ and $\Dgo$ the pre-complete $C^{*}$-category of $P$-closed Hilbert $\Cgo$-module. Consider the following conditions on $P$:

\begin{itemize}

\item[(PC2)] The functor $P : \Hcal \Cgo \rightarrow \Dgo$ is continuous.
\item[(PC3)] If one has two Hilbert $\Cgo$-modules $X \subset Y \subset PX$ with $X \subset PX$ being the canonical inclusion then there is an isomorphism between $PY$ and $PX$ compatible with the inclusions of $Y$.
\item[(PC4)] Every pre-adjoint pair for $P$ has an extension (as in \ref{def_extension_prop}).
\end{itemize}

Then $(PC2) \Rightarrow (C2)$; $(PC3) \Leftrightarrow (C3)$ and $(PC4) \Leftrightarrow (C4)$. If moreover $\Cgo \subset \Dgo$ as subcategory of $\Hcal \Cgo$ then $(PC2) \Leftrightarrow (C2)$ as well.

}

\Dem{

$(PC2) \Rightarrow (C2) :$ The functor $L: \Hcal \Dgo \rightarrow \Dgo $ is defined as $P \circ i_*$ where $i : \Dgo \rightarrow \Hcal \Cgo$ is the canonical inclusion. As $\Hcal \Cgo$ satisfies $(C2)$, the functor $i_*:\Hcal \Dgo \rightarrow \Hcal \Cgo$ is continuous hence if $P$ is continuous this implies that $L$ is continuous. 

\bigskip

$(C2) \Rightarrow (PC2)$ when $\Cgo \subset \Dgo$: We will show that the functor $P: \Hcal \Cgo \rightarrow \Dgo$ is  $j_* = L\circ (\Hcal j)$ where $j$ is the inclusion of $\Cgo$ into $\Dgo$. As $\Hcal j$ is continuous this is enough to prove that the continuity of $L$ implies the continuity of $P$. As $L = P \circ i_*$, one has $L\circ (\Hcal j) = P \circ i_* \circ (\Hcal j)$ but $i_* \Hcal j = (i\circ j)_{*}$ and $i\circ j$ is the Yoneda embeddings of $\Cgo$ in $\Hcal \Cgo$ hence $ (i \circ j )_{*}$ is the identity functor of $\Hcal \Cgo$ which concludes the proof. 

\bigskip

$(PC3) \Rightarrow (C3) :$ Let $\Acal \subset \Bcal \subset \Dgo(X)$ two hereditary subalgebra with $\Acal$ weakly dense. One obviously has an inclusion $\Acal.X \subset \Bcal.X \subset X$ and $P(\Acal.X)$ is isomorphic to $X$ in a way compatible to the inclusion of $\Acal.X$. Hence $(PC3)$ implies that $P(\Bcal.X)$ is isomorphic to $P(\Acal.X) \simeq X$ in a way compatible to the inclusion of $\Bcal.X$, as the embedding of $\Bcal X$ in $P(\Acal.X)$ has been defined such that the isomorphism $P(\Acal.X) \simeq X$ is also compatible to the embeddings of $\Bcal .X$ this proves that one has an isomorphism between $P(\Bcal.X)$ and $X$ compatible to the inclusion of $\Bcal.X$ and hence that $\Dgo$ satisfies $(C3)$.

\bigskip

$(C3) \Rightarrow (PC3) :$ If one has $X \subset Y \subset PX$ in $\Hcal\Cgo$, they correspond by proposition \ref{Prop_H'=AH} to two hereditary sub-algebras $\Acal \subset \Bcal \subset \Dgo(PX)$ such that $X = \Acal.PX$ and $Y=\Bcal.PY$ with the correct inclusion. In particular $\Acal$ is weakly dense hence $\Bcal$ is weakly dense by $(C3)$ hence there is an isomorphism of $PY$ and $PX$ compatible with the inclusion of $Y$ and this concludes the proof.

\bigskip

$(C4)\Rightarrow (PC4):$ Let $\Dgo_{\Kcal}$ be the category of $P$-closed Hilbert $\Cgo$-modules and $\Cgo$-compact operators between them, i.e. $\Dgo_{\Kcal} = \Kcal \Cgo \cap \Dgo$. By definition $\Kcal \Cgo$ is the two sided ideal of $\Hcal \Cgo$ generated by $\Cgo$, but one can also easily see that it is the two sided ideal of $\Hcal \Cgo$ generated by $\Dgo_{\Kcal}$: indeed as any object of $\Cgo$ can be identified with a sub-object of a module in $\Dgo$, any arrow in $\Cgo$ can be factored (at least approximatively) into a compact arrow between objects of $\Dgo$. Hence by corollary \ref{Cor_eqHA_HI} the categories $\Hcal\Dgo_{\Kcal}$, $\Hcal \Cgo$ and $\Hcal \Kcal \Cgo$ are all equivalent and these equivalences are also defined on linear bounded maps and preserve the adjunction relation. Hence instead of trying to prove that the closure operation on $\Hcal\Cgo$ satisfies $(PC4)$ one can prove that the corresponding closure operation on $\Dgo_{\Kcal}$ satisfies $(PC4)$.

Let $A \rightarrow PB$ and $B \rightarrow PA$ be a pre-adjoint pair in $\Hcal\Dgo_{\Kcal}$. As $\Dgo_{\Kcal}$ is a two sided ideal of $\Dgo$, a Hilbert $\Dgo_{\Kcal}$-module can be seen as a Hilbert $\Dgo$-module whose scalar product takes values in $\Dgo_{\Kcal}$, moreover, by definition of the pre-complete structure on $\Dgo$, $L$ is the functor that take a Hilbert $\Dgo$-module, remember only the element whose self-scalar product is in $\Dgo_{\Kcal}$ and then apply $P$ to it which yields an object of $\Dgo$. Hence by seeing $A$ and $B$ as $\Hcal\Dgo$-modules one obtain exactly an extension pair in $\Hcal \Dgo$, and hence one would obtain an extension because $\Dgo$ satisfies $(C4)$.

\bigskip

$(PC4)\Rightarrow (C4) :$ Let $A$ and $B$ be two Hilbert $\Dgo$-modules and let $f:A \rightarrow LB$ and $:B \rightarrow LA$ as in the extension property (see definition \ref{def_extension_prop}). Let $ i: \Dgo \rightarrow \Hcal \Cgo$ the natural inclusion, the Hilbert $\Cgo$-module $i_* A$ (and similarly $i_* B$) is generated by elements of the form $a \otimes_d c \in i_*A(C)$ for $d \in |\Dgo|, a \in A(d), c: c \rightarrow d$ with the scalar product described in \ref{Prop_HC_precomp}. $L A$ and $LB$ are by definition $Pi_*A$ and $Pi_*B$, and one can see that the maps $f$ and $g$ can be pushed to maps from $i_* A \rightarrow Pi_* B$ and $i_* B \rightarrow Pi_* A$ by defining $f'(a \otimes_d c)= f(a).c$ one can then easily check that $f'$ and $g'$ defined this way still satisfies the condition of the extension property hence as $P$ satisfies $(PC4)$ defines a morphism $\tilde{f}$ from $LA$ to $LB$ which can be checked to be an extension of $f$.

}

}

\block{There is one additional natural and convenient conditions that one can consider on a closure operation, unfortunately it has no clear analogue in terms of properties of pre-complete $C^{*}$-category:

\Def{Let $P : \Hcal \Cgo \rightarrow \Hcal \Cgo$ a closure operation. $P$ is said to be a reflection if it satisfies:
\begin{itemize}

\item[(PC5)] Any bounded linear map from a Hilbert $\Cgo$-module $A$ to a closed Hilbert $\Cgo$-module $B$ extend uniquely into a bounded linear map from $P(A)$ to $B$.

\end{itemize}
}

This condition is not as well behaved as the four first conditions in the sense that it is not the case that if one closure operation defining a pre-complete category $\Dgo$ satisfies it then all closure operations will (as it is shown for the others in the previous proposition). This should probably be attributed to the fact that this condition is not really a $C^{*}$-categorical property. For example, if $\Cgo$ is a $C^{*}$-category then the identity of $\Hcal \Cgo$ is a closure operation satisfying $(PC5)$ and define the category of Hilbert $\Cgo$-modules, but the closure operation $\Yon \circ R_{\Cgo} : \Hcal \Hcal \Cgo \rightarrow \Hcal \Cgo \rightarrow \Hcal \Hcal \Cgo$ also has $\Hcal \Cgo$ as category of closed object but almost never satisfies $(PC5)$, indeed assume it satisfies $(PC5)$ then take $f: H \rightarrow H'$ any bounded linear map, one can define two Hilbert $\Hcal\Cgo$-modules $\Kcal H$ and $\Kcal H'$ defined by $\Kcal H (Y)= \Kcal(Y,H)$ and one easily see that $f$ defines a bounded linear map from $\Kcal H$ to $\Kcal H'$, but one easily see $P(\Kcal H) \simeq \Yon(H)$ hence $P'(f)$ should in particular be a bounded linear map from $\Yon(H)$ to $\Yon(H')$ but because of the Yoneda lemma this is automatically an operator and one easily see that it has to be equal to $f$ (because it is equal to $f$ when restricted to $\Kcal H$). Hence the existence of such an operator would imply that any bounded linear map in $\Hcal \Cgo$ is in fact an operator, which, except for some exceptional situation like $\Cgo = \C$ assuming the law of excluded middle, is not the case.

Also not that the uniqueness in the condition is (at least apparently) stronger than the condition $(PC1)$ because it allows to show that a bounded linear map is zero on $P(A)$ if and only if it is zero on $A$ while $(PC1)$ only show this for operators and it does not seems possible to deduce this uniqueness from $(PC1)$.
}

\block{\label{Prop_PC5impPCi}\Prop{$(PC5)$ implies both $(PC1)$ and $(PC4)$.}

\Dem{As we mentioned above, the uniqueness property in $(PC5)$ is a strong form of $(PC1)$: Let $f:P(A) \rightarrow B$ an arrow in $\Dgo$ which is zero on $A$, then it corresponds under the adjunction to the zero map from $A$ to $B$ hence is zero, which proves that $P$ satisfies $(PC1)$.

Assume now that one has $f:A \rightarrow P(B)$ and $g:B \rightarrow P(A)$ a pre-adjoint pair. Then one has two extension as bounded linear map given by $(PC5)$: $f':P(A) \rightarrow P(B)$ and $g':P(B) \rightarrow P(A)$, we need to prove that they are adjoint in order to conclude that they defines an operator from $P(A)$ to $P(B)$.

Fix a $a \in A(X)$, then for any $b \in B(Y)$, where $X$ and $Y$ are two objects of $\Cgo$, one has $\scal{a}{g'(b)} = \scal{f'(a)}{b} \in \Cgo(Y,X)$, hence as bounded linear maps from $B$ to $\Yon(X)$, $\scal{a}{g'( \_)}$ and $\scal{f'(a)}{\_}$ agree, but they are also defined as maps from $P(B)$ to $\Yon(X)$, hence their extension to $P(B)$ agree (if $\Yon(X)$ is not it self closed one can just compose them with the embeddings of $\Yon(X)$ in $P(\Yon(X))$), hence one has for all $a \in A$ for all $b \in P'(B)$ $\scal{a}{g'(b)} = \scal{f'(a)}{b}$ applying the same argument a second time reversing the role of $a$ and $b$ gives that $f'$ and $g'$ are adjoint.

}

}

\blockn{One can also observe that $(PC5)$ turn the category of $P$-closed Hilbert modules and bounded linear map into a reflexive subcategory of the category of Hilbert modules and bounded linear maps. $P$ is the reflection, and the natural inclusion $\mu_X : X \hookrightarrow P(X)$ is the unit of adjunction. This is why we call such closure operator a reflection.}

\block{\Prop{If $\Ccal$ is a pro-$C^{*}$-algebra then the closure operation on $\Ccal^{b}$-modules defined by taking the bounded element of the completion for the pro-topology is a reflection on $\Hcal \Ccal^{b}$.

Assuming classical mathematics, for any $C^{*}$-algebra $C$ taking the bidual is a reflection on $\Hcal C$.

In particular both the category of Hilbert modules over a pro-$C^{*}$-algebra and (assuming classical mathematics) of reflexive Hilbert modules over a $C^{*}$-algebra are pre-complete categories satisfying all the conditions $(C1),(C2),(C3)$ and $(C4)$.  }

\Dem{For the case of prop-$C^{*}$-algebra, any bounded linear transformation is uniform for the pro-topology so the existence and uniqueness of the extension follow easily from usual result on completion of uniform space, the fact that the extension indeed takes value in bounded elements is a corollary of lemma \ref{Lem_boundedlinearmapIneq}.

For the bidual construction, if $f$ is a bounded linear map from a Hilbert $C$-module $A$ to a reflexive Hilbert $C$-module $B$ then one obtains a map $f':B' \rightarrow A'$ by precomposition and a map $f'':A'' \rightarrow B''=B$ again by pre-composition. This is an extension of $f$ and it is proved in \cite{paschke1974double} that the norm of a linear map on $A''$ is the same as the norm of its restriction to $A$, hence this extension is unique.

By the previous proposition this implies that the corresponding closure operation satisfies the conditions $(PC1)$ and $(PC4)$, which in turn by \ref{Prop_PC1} and \ref{Prop_RelCiPCi} that the corresponding category satisfies $(C1)$ and $(C4)$ and hence it also satisfies $(C2)$ and $(C3)$ by \ref{Prop_impbetween_Ci}.
}
}

\subsection{Generators and comparison theory}
\label{subsection_Generators_comparison}

\block{\label{DefGenerator}\Def{Let $\Cgo$ be a pre-complete $C^{*}$-category, $\Ago$ a hereditary sub-category with $i: \Ago \rightarrow \Cgo$ the inclusion $C^{*}$-functor, $\Ago$ is said to be generating if there is an isomorphism of $C^{*}$-functor: $\mu : Id_{\Cgo} \simeq i_* R_{\Ago}$ such that for any arrow $f \in \Cgo(X,Y)$ with $f^{*}f \in \Ago(X)$ the following diagram commute:

\[
\begin{tikzcd}[ampersand replacement=\&]
X \arrow{rd}{L(\tilde{f})} \arrow{d}{f} \& \\
Y \arrow{r}{\mu} \& i_* R_{\Ago}(Y)\\
\end{tikzcd}
\]

Where $\tilde{f}$ is the map  $X \rightarrow \Hcal i R_{\Ago}(Y)$ induced by $f$.

}

This definition is not very convenient to check in practice, but proposition \ref{Prop_GeneratingHered=generatingIdeal} and \ref{Prop_generating=dense} below gives an easy equivalent condition to test whether a hereditary subcategory is generating or not, in the case of a category satisfying $(C1)$.
}

\blockn{In order to obtain a weaker condition for a hereditary subcategory to be generating one can start with the following observation:}

\block{\label{Prop_GeneratingHered=generatingIdeal}\Prop{Let $\Cgo$ be a pre-complete $C^{*}$-category and $\Ago$ a hereditary subcategory then $\Ago$ is generating if and only if the two sided ideal $\Igo$ spammed by $\Ago$ is generating.}

\Dem{Let $i : \Ago \rightarrow \Cgo$ , $j : \Ago \rightarrow \Igo$ and $k : \Igo \rightarrow \Cgo$ the natural inclusion. One has $i = k \circ j$ and hence $i_* = k_* \circ \Hcal i$.

We denote $R_{\Ago}$ and $R_{\Igo}$ the restriction functor from $\Cgo$ to $\Hcal\Ago$ and $\Hcal \Igo$, and by $F$ the restriction functor $\Hcal \Igo \rightarrow \Hcal \Ago$, one has $R_{\Ago}= F \circ R_{\Igo}$. Hence:

\[ i_* R_{\Ago} \simeq k_* \circ  \Hcal j \circ F \circ R_{\Igo}\]

But by corollary \ref{Cor_eqHA_HI}, the functors $\Hcal j$ and $F$ are two equivalences inverse of each other. Hence:

\[ i_* R_{\Ago} \simeq k_* R_{\Igo} \]

and hence one can pass from an isomorphism $i_* R_{\Ago} \simeq Id_{\Cgo}$ into an isomorphism $k_* R_{\Igo}  \simeq Id_{\Cgo}$. It remains to check that this preserves the compatibility condition. For any $f \in \Cgo(X,Y)$ such that $f^{*}f \in \Ago(X) \subset \Igo(X)$, the isomorphism $i_* R_{\Ago}(Y) \simeq k_* R_{\Igo}(Y)$ is compatible to the maps $L(\tilde{f})$ hence if $\Igo$ is generating then $\Ago$ is generating. Conversely if the one has an isomorphism which satisfies the compatibility condition for $\Ago$, then for any $f \in \Igo$ of the form $g a h$ with $a \in \Ago$ one can see that $\tilde{f} = (\Hcal i R_i)(g) \circ \tilde{a} \circ h$ hence one can deduce the compatibility condition for $f$ from the condition for $a$ and the functoriality of the isomorphism $\mu$, as $\Igo$ is spammed by such maps this concludes the proof.
}

}

\block{\label{Prop_generating=dense}\Prop{Let $\Cgo$ be a pre-complete $C^{*}$-category and $\Igo$ be a (total) two sided ideal. If $\Igo$ is generating then for any object $X \in |\Cgo|$, $\Igo(X)$ is a weakly dense hereditary sub-algebra of $\Cgo(X)$. If $\Cgo$ satisfies $(C1)$ then the converse is true: $\Igo$ is generating if and only if $\Igo(X)$ is weakly dense in $\Cgo(X)$ for all $X$.
}

Note that we in fact only need that $\Cgo$ satisfies $(C1)$ for two sided ideal (i.e. that weakly dense two sided ideal are essential)

\Dem{Let $L : \Hcal \Cgo \rightarrow \Cgo$ be the structural functor, $i : \Igo \hookrightarrow \Cgo$ the inclusion. Let $X \in \Cgo$ we denote by $\Igo X$ the Hilbert $\Cgo$-module $\Igo (X) \Yon(X)$. One easily see that $(\Hcal i) (R_{\Igo}(X))$ is naturally isomorphic to $\Igo X$. Hence $i_* \circ (R_{\Igo})$ is identified with $L(\Igo X)$. In particular if $\Igo$ is generating one gets an isomorphism between $L(\Igo X)$ and $X$, and the compatibility condition on this isomorphism gives the compatible to the inclusion of $\Igo X$ which proves the density of $\Igo(X)$ in $\Cgo(X)$.

Conversely, if $\Igo(X)$ is dense in $\Cgo(X)$ for all $X$, then one has an isomorphism between $i_* \circ R_{\Igo}(X)$ and $X$ for all $X$. But using $(C1)$ one easily shows that this isomorphism is unique and functorial. The compatibility condition for this isomorphism follows directly from the condition satisfied by the ``density'' isomorphism.
}
}

\blockn{We will now investigate the consequence of having a generating subcategories.}

\block{\label{PropGeneratorsAndClosOp}\Prop{Let $\Cgo$ be a pre-complete $C^{*}$-category and let $\Ago$ be a generating hereditary subcategory. Then $P= R_{\Ago} \circ i_*$ is a closure operation on $\Hcal \Ago$ and the category of $P$-closed Hilbert $\Ago$-module is equivalent to $\Cgo$ as a $C^{*}$-category.

Moreover, if $\Cgo$ satisfies $(C2)$ or if $i_*$ is continuous then this is an equivalence of pre-complete $C^{*}$-category.

}

\Dem{We will first show that for each $X \in \Cgo$, $R_{\Ago}(X)$ is $P$-closed. Because $\Ago$ is generating, one has for each $X \in \Cgo$ an isomorphism $\mu_X$ between  $i_* R_{\Ago} (X)$ and $X$, and the compatibility condition of definition \ref{DefGenerator} exactly assert that this isomorphism $\mu_X$ induce the canonical inclusion of $R_{\Ago}(X)$ into $R_{\Ago} i_* R_{\Ago}$ and hence that $R_{\Ago}(X)$ is $P$-closed.  

This already proves that $R_{\Ago} i_*$ is a closure operation on $\Hcal \Ago$ : it is a pre-closure by \ref{Prop_preclosedOp} and we just proved that $R_{\Ago} X$ is always closed.

Let $\Dgo$ be the category of $P$-closed Hilbert $\Ago$-module. One has two functors: $R_{\Ago} : \Cgo \rightarrow \Dgo$ and $i_* : \Dgo \rightarrow \Cgo$ which are inverse of each other: $R_{\Ago} i_*$ is the identity on $\Dgo$ exactly because $\Dgo$ is the category of $P$-closed object and $i_* R_{\Ago}$ is the isomorphic to the identity of $\Cgo$ because $\Ago$ is generating. 

If $\Cgo$ satisfies $(C2)$ then $i_*$ is continuous, hence preserve the $f_*$ construction, and hence the $f_*$ construction computed one both side of this equivalence are the same which proves that one has an isomorphism of pre-complete $C^{*}$-categories.
}
}

\blockn{We can also extend proposition \ref{Prop_restrictionoftransnat1} to any generating subcategory:}

\block{\Prop{Let $\Cgo$ and $\Dgo$ be two pre-complete $C^{*}$-categories, $\Ago \subset \Cgo$ a generating partial two sided ideal \footnote{We mean that the corresponding hereditary subcategory obtained by adding all the zero homomorphism is generating} and assume that $\Dgo$ satisfies $(C1)$.

Then the restriction functor from the category of continuous $C^{*}$-functors from $\Cgo$ to $\Dgo$ to the category of $C^{*}$-functors from $\Ago$ to $\Dgo$ is fully faithful.}

Note that, strictly speaking, this statement does not makes sense in our framework as it quantify over functors between large category which is not allowed. But it can be interpreted as a meta-theorem, or by strengthening our foundations either with a notion of 2-class or by allowing the class of function between two class.

\Dem{Let $\Igo$ be the total two sided ideal generated by $\Ago$. By propositions \ref{Prop_GeneratingHered=generatingIdeal} and \ref{Prop_generating=dense} for each object $X \in \Cgo$, $\Igo(X)$ is weakly dense. Moreover as $\Ago$ is already a partial two sided ideal $\Ago(X) = \Igo(X)$ for any $X \in |\Ago|$.

Let $i \Ago \rightarrow \Cgo$ be the inclusion, as pointed out above, $i$ satisfies the weak density condition of proposition \ref{Prop_weaklyNDg}, hence as $\Dgo$ satisfies $(C1)$ its image in $\Dgo$ under any continuous functor is weakly non-degenerated. Moreover, As $\Ago$ is generating one has an isomorphism $\mu: Id_{\Cgo} \simeq i_* R_{\Ago}$.

Let $F : \Cgo \rightarrow \Dgo$ be a continuous $C^{*}$-functor, $F \simeq F \circ i_* \circ R_{\Ago}$ as $F$ is continuous, $F \circ i_* \simeq (F \circ i)_*$ and $F \circ i$ is weakly non-degenerated. Hence, by proposition \ref{Prop_restrictionoftransnat1} for any two such continuous functors $F,G : \rightrightarrows \Dgo$, natural transformations from $F \circ i$ to $G \circ i$ are the same as natural transformations from $(F \circ i)_*$ to $(G \circ i)_*$, hence gives rise to a natural transformation from $F \simeq F \circ i_* \circ R_{\Ago}$ to $G \simeq G \circ i_* \circ R_{\Ago}$.

Moreover as $\Ago$ is a partial two sided ideal, one can check that $R_{\Ago}$ send the embeddings of $\Ago$ in $\Cgo$ to the Yoneda embeddings of $\Ago$ into $\Hcal \Ago$, hence this natural transformation from $F$ to $G$ is an extension of the natural transformation from $F \circ i$ to $G \circ i$ it comes from.

Finally such an extension is unique by the exact same argument as in the proof of proposition \ref{Prop_restrictionoftransnat1}, or more explicitly: If $\mu:F \rightarrow G$ is a natrual transformation whose restriction $\mu'$ is equal to zero, then for $X \in |\Cgo|$ any object let $h=f \circ a \circ g$ be an endomorphism of $X$ which factor into an arrow $a \in \Ago(U,V)$. Because $\mu$ is zero on all objects of $\Ago$ this immediately show that $\mu_X F(h) = G(f) \mu_V F(a \circ g) = 0$. Hence for any $i \in \Igo(X)$ one has $\mu_X F(i)= 0$. As $F$ is continuous $F(\Igo(X))$ is weakly dense and hence as $\Dgo$ satisfies $(C1)$ this shows that $\mu_X = 0$ and hence that the restriction functor is faithful.  

}

}

\block{\Prop{If $\Cgo$ is a pre-complete $C^{*}$-category and $\Ago$ is a generating hereditary subcategory then the functor $R_{\Ago}$ is fully faithful.}

\Dem{This follows directly from the fact that $i_* R_{\Ago} \simeq Id$ and the fact that $i_* = L \circ \Hcal i$ is faithful as a composite of two faithful functors.}
}

\blockn{As for $W^{*}$-category theory, the notion of generators is closely related to a notion of quasi-containement:}

\block{\label{Def_quasicontain}\Def{Let $\Cgo$ be a pre-complete $C^{*}$-category and let $X,Y \in |\Cgo|$ be any object. We say that $Y$ is quasi-contained in $X$, and write $Y \prec X$ if the total two sided ideal of $\Cgo$ generated by $\Cgo(X)$ is weakly dense in $\Cgo(Y)$.

We will also note $Y \prec \Ago$ and say that $Y$ is weakly contained in $\Ago$ if $\Ago$ is any subcategory of $\Cgo$ and the two sided ideal ideal generated by $\Ago$ is weakly dense in $\Cgo(Y)$.}}

\block{\Def{Let $\Cgo$ be a pre-complete $C^{*}$-category. One says that $\Dgo$ is a pre-complete full subcategory of $\Cgo$ if:

\begin{itemize}
\item $\Dgo$ is a full subcategory of $\Cgo$.
\item If $c \in |\Cgo|$ is isometrically isomorphic to $d \in |\Dgo|$ then $c \in |\Dgo|$.
\item if $i : \Ago \rightarrow \Dgo$ is a $*$-functor and $H$ is a Hilbert $\Ago$-module, then $i_* H \in \Dgo$.

\end{itemize}
}

Of course, such a subcategory $\Dgo$ is automatically endowed with an induced pre-complete structure such that the $i_*$ for the pre-complete structures of $\Cgo$ and of $\Dgo$ are the same.
}

\block{\label{Prop_quasicontainement}\Prop{Let $\Cgo$ be a pre-complete $C^{*}$-category satisfying $(C1)$ and $(C3)$ and let $X \in |\Cgo|$ then:
\begin{enumerate}

\item An object $Y$ is quasi-contained in $X$ if and only if it is isomorphic to an object of the form $i_* H$ for $i$ the inclusion of $\Cgo(X)$ into $\Cgo$ and $H$ a $\Cgo(X)$-module. More generally, if $\Ago \subset \Cgo$ is a partial hereditary subcategory then $Y$ is weakly contained in $\Ago$ if and only if it is of the form $i_*H$ for $H$ a Hilbert $\Ago$-module.

\item Quasi-containement is a transitive and reflexive relation on objects of $\Cgo$.

\item \label{Prop_quasicontainement_closureOp} For any $\Ago \subset \Cgo$ a hereditary subcategory, the pre-closure operation on $\Hcal \Ago$ induced by $\Cgo$ as in \ref{Prop_preclosedOp} is a closure operation and the category of closed $\Ago$-module is equivalent to the full subcategory of $\Cgo$ of object weakly contained in $\Ago$.

\item A hereditary subcategory $\Ago$ is generating $\Cgo$ if and only if any object is quasi-contained in $\Ago$.

\item The full sub-category $\Cgo_X$ of $\Cgo$ of objects weakly contained in $X$ is the smallest pre-complete full subcategory of $\Cgo$ which contains $X$.

\item $\Cgo_X$ is the unique pre-complete full sub-category of $\Cgo$ generated by $X$.

\end{enumerate}

}

For this proposition we did not try to optimise the assumption on $\Cgo$: some of these points can clearly be proven under weaker assumption than $(C1)$ and $(C3)$

The last two points can also be formulated for a subcategory $\Ago$ with the subtleties that without additional assumptions on $\Ago$ there is no reason for the objects of $\Ago$ to be themselves weakly contained in $\Ago$, we need for example the assumption that $\Ago$ is non-degenerated in the sense that for every $a \in |\Ago|$, $\Ago(a)$ is weakly dense in $\Cgo(a)$. Once an assumption of this kind is added then the exact same proof caries over to this case.

\Dem{

\begin{enumerate}

\item We only need to prove the more general form (the other case follow by taking $\Ago = \Cgo(X)$).

Let $\Kcal_{\Ago}$ be the two sided ideal generated by $\Ago$ (the class of $\Ago$-compact operators) and let $i : \Ago \hookrightarrow \Cgo$ the inclusion functor.
Any object $Y$ has naturally the structure of a cone for $(i,R_{\Ago}(Y))$, indeed any element of $R_{\Ago}(Y)(a)$ actually is a map from $a$ to $Y$. Moreover, the sub-algebra generated by this cone structure as in proposition \ref{PropCompactActingOnCones} (any of them) is the algebra generated by arrows factoring into $\Ago$ hence it is exatly $\Kcal_{\Ago}(X)$.

By proposition \ref{limitanddensity}, this gives that $Y$ is weakly contained in $\Ago$ if and only if $Y=i_* R_{\Ago} Y$, hence one already have one implication: if $Y$ is weakly contained in $\Ago$ then it is of the form $i_* H$ for some $\Ago$-module $H$.

Conversely, if $Y=i_* H$ for $H$ a Hilbert $\Ago$-module, then the algebra of operators induced by the cone structure on $Y$ is made of $\Ago$-compact operator and hence is included in $\Kcal_{\Ago} \wedge \Cgo(Y)$, hence as $\Cgo$ satisfies $(C3)$ the weak density of the algebra generated by the cone structure implies that $Y$ is quasi-contained in $\Ago$.

\item For the second point, the reflexivity of $\prec$ is obvious. We assume that $X \prec Y$ and $Y \prec Z$. By the previous observation, $Y \simeq i_* \Cgo(X,Y)$ and $Z \simeq j_* H$ where $i$ and $j$ denotes respectively the inclusion of $\Cgo(X)$ and $\Cgo(Y)$ into $\Cgo$ and $H$ is a Hilbert $\Cgo(Y)$-module. $\Cgo(Y)$ acts on $\Cgo(X,Y)$ (as endomorphisms of a Hilbert $\Cgo(X)$-module.) Hence the functor $j$ can be factored into $i_*\circ j'$ where $j'$ in the functor corresponding to the action of $\Cgo(Y)$ on $\Cgo(X,Y)$ in the category of Hilbert $\Cgo(X)$-modules. As $\Cgo$ satisfies $(C1)$ and $(C3)$ it also satisfies $(C2)$ by proposition \ref{Prop_impbetween_Ci}, hence $i_*$ is a continuous $C^{*}$-functor, and hence $Z \simeq j_*(H)= i_*(j'_*(H))$can be written as $i_*H'$ with $H' = j'_*(H)$ a Hilbert$\Cgo(X)$ modules which proves that $X \prec Z$

\item We still denote by $i$ the inclusion of $\Ago$ in $\Cgo$. Let $H$ be any Hilbert $\Ago$-module. By the point $(1)$, $i_* H$ is weakly contained in $\Ago$. Hence (as shown in the proof of point $(1)$) $i_* R_{\Ago} i_*H \simeq i_* H$ and finally there is a natural isomorphism between $R_{\Ago} i_* R_{\Ago} i_* H$ and $R_{\Ago} i_*H$ which shows that $R_{\Ago} i_*$ is indeed a closure operation. If $X$ is an object of $\Cgo$ weakly contained in $\Ago$, then $X=i_*H$ and hence $R_{\Ago} X$ is a closed $\Ago$-module, hence $R_{\Ago}$ and $i_*$ form a pair of functors between the category of closed $\Ago$-module and object of $\Cgo$-weakly contained in $\Ago$ and restricted to these subcategory they are inverse of each other because of the isomorphism mentioned above.

\item As $\Cgo$ satisfies $(C1)$, propositions \ref{Prop_GeneratingHered=generatingIdeal} and \ref{Prop_generating=dense} implies that $\Ago$ is generating if and only if the two sided ideal generated by $X$ is weakly dense on every object which is exactly the statement that any object is quasi-contained in $X$.

\item Let $\Cgo_X$ be the full subcategory of objects quasi-contained in $X$, as any object of $\Cgo_X$ can be written $i_* H$ for $H$ a Hilbert $\Cgo(X)$-module any full subcategory of $\Cgo$ pre-complete for the induced structure contains $\Cgo_X$. We just need to show that $\Cgo_X$ itself is pre-complete for the induced structure. One has a functor $\Cgo(X,\_)$  from $\Cgo_X$ to $\Hcal \Cgo(X)$ (it is in fact $R_{\Cgo(X)}$) and one has observed above that for any object $Y \in |\Cgo_X|$ $Y \simeq i_* \Cgo(X,Y)$, moreover (by $(C1)$) this can be made into an isomorphism of functors between the identity of $\Cgo_X$ and $i_* R_{\Cgo}$. Hence for any $*$-functor $j : \Dgo \rightarrow \Cgo_X$, $j$ can be factored into $j \simeq i_* \circ j'$ where $j'$ is $R_{\Cgo} \circ j$, and because $\Cgo$ satisfies $(C2)$, $i_*$ is a continuous functor hence for any Hilbert $\Dgo$-module $H$, one has $j_* H \simeq i_* (j'_* H) \in |\Cgo_X|$ which concludes the proof.

\item In a full subcategory $\Dgo$ of $\Cgo$ which is pre-complete for the pre-complete structure of $\Cgo$, being weakly dense for the pre-complete structure of $\Dgo$ or of $\Cgo$ is the same, hence for $Y,Z \in \Dgo$ , $Y \prec Z$ in $\Cgo$ if and only if $Y \prec Z$ in $\Dgo$. Hence any such full subcategory in which $X$ is generating has to be the category of $Y$ such that $Y \prec X$.

\end{enumerate}

}

}

\subsection{Compactness and categories of Hilbert modules}
\label{subsection_compact}
\blockn{Finally we would like to obtain a characterization of category of the form $\Hcal \Cgo$. It is known (see proposition \ref{Cor_eqHA_HI}) that such a category is written under the form $\Hcal \Dgo$ in a canonical way by taking $\Dgo = \Kcal \Cgo$ the sub-category of compact operators. Hence it seems natural that one should seek for an abstract characterization of this class of morphisms. This is done by mimicking to role of compact operator in proposition \ref{densityforHmodule}, which has we have seen in proposition \ref{Prop_weakdensityinHC} characterize weak density in the pre-complete category $\Hcal \Cgo$.
}

\block{\label{Def_Abscpt}\Def{Let $\Cgo$ be a pre-complete category. A morphism $k \in \Cgo(X,Y)$ is said to be absolutely compact if for all object $Z$ of $\Cgo$, for all morphisms $h_1:Z \rightarrow X$ and $h_2:Z \rightarrow Y$ and for all $A \subset \Cgo(Z)$ a weakly dense hereditary sub-algebra $h_2^{*} k h_1 \in A$.

Absolutely compact operators form a total two sided ideal of $\Cgo$ denoted $\Cgo_{\Kcal}$.
}

One can see that in $\Hcal \Cgo$ absolutely compact operators are exactly the $\Cgo$-compact operator, and it will become clear in the rest of this section.
}

\block{The ideal $\Cgo_{\Kcal}$ of absolutely compact operators satisfies a kind of universal property:

\Prop{Let $\Cgo$ be a complete $C^{*}$-category and $\Ago \subset \Cgo$ be a hereditary sub-category then the functor $ R_{\Ago} : \Cgo \rightarrow \Hcal \Ago$ is continuous if and only if $\Ago \subset \Cgo_{\Kcal}$.

}

In concrete cases where the continuity of a functor is related to the preservation of some sort of weak convergence, like for $W^{*}$-category or Hilbert modules over a pro-$C^{*}$-algebra, this proposition means that absolutely compact operators are those such that composition with them turn weak convergence into norm convergence. 

\Dem{Assume first that $\Ago \subset \Cgo_{\Kcal}$. Because $\Hcal \Ago$ satisfies $(C1)$ it suffice to check that $R_{\Ago}$ preserves weakly dense sub-algebras to check that it is a continuous functor.

Let $Z \in |\Cgo|$ and $\Acal$ be weakly dense hereditary sub-algebra of $\Cgo(Z)$. Let $a \in |\Ago|$ and $ h \in R_{\Ago}(Z)(a)$ i.e. $h$ is a map in $\Cgo(a,Z)$ from $a$ to $Z$ such that $h^{*}h \in \Ago(a)$, in particular $h^{*}h$ is absolutely compact and as absolutely compact operators form a two sided ideal this also proves that $h$ is absolutely compact and hence that $hh^{*} \in \Acal$ (as $\Acal$ is weakly dense) which proves that in $\Hcal\Acal$ $hh^{*}$ is in $R_{\Ago}(\Acal)$. This being true for any $h \in R_{\Ago}(Z)$ this proves that $R_{\Ago}(\Acal)$ contains all the $\Ago$-compact operator of $R_{\Ago}(Z)$ and hence is weakly dense.

\bigskip

Conversely assume that $R_{\Ago}$ is continuous, let $X,Y \in |\Ago|$ , $a \in \Ago(X,Y)$ one need to prove that $a$ is absolutely compact. Let $Z \in |\Cgo|$ be any object, $\Acal \subset \Cgo(Z)$ be a weakly dense hereditary sub-algebra and let $h_1 : Z \rightarrow X$ and $h_2 : Z \rightarrow Y$, we need to prove that $h_2^{*} a h_1 \in \Acal$. One can always assume that $a$ can be written $a=a_2 a_1$ with both $a_1$ and $a_2$ in $\Ago$, (for example, approximately using lemma \ref{Lem_hered}(\ref{Lem_hered3}) ) one hence has $(h_2)^{*}(a_2)$ and $ (h_1)^{*}(a_1)^{*}$ which are two element of $R_{\Ago}(Z)$. Because $R_{\Ago}$ is continuous, the hereditary sub-algebra generated by $R_{\Ago}(\Acal)$ is weakly dense, hence for any $\lambda_i$ an approximate unit of $\Acal$, one has $\lambda_i (h_2)^{*}(a_2)$ and $\lambda_i  (h_1)^{*}(a_1)^{*}$ which converge respectively to  $(h_2)^{*}(a_2)$ and $ (h_1)^{*}(a_1)^{*}$ (in $R_{\Ago}(Z)$ but this is the same as the convergence in $\Cgo$). In particular, $\lambda_i h_2^{*}a_2 a_1 h_1 \lambda_i$ converge in $\Cgo$ to $h_2^{*} a h_1$ hence $ h_2^{*} a h_1 \in \Acal$ which concludes the proof.
}

}

\block{\label{Th_caracofHmod}\Th{A pre-complete $C^{*}$-category is equivalent to $\Hcal \Dgo$ if and only if $\Dgo$ can be embedded in $\Cgo$ as a generating hereditary subcategory of absolutely compact operator. Moreover the equivalence between $\Cgo$ and $\Hcal \Dgo$ is induced by $R_{\Dgo}$.}

Note that given such an embeddings of $\Dgo$ into $\Cgo$, the identification of $\Cgo$ with $\Hcal \Dgo$ does not necessarily turn this embeddings into the Yoneda embeddings. One can check that this will be the case if and only if $\Dgo$ is chosen as a two sided ideal of a full subcategory of $\Cgo$.

\Dem{If $\Cgo$ is $\Hcal \Dgo$ then the image of the Yoneda embeddings of $\Dgo$ is generating and absolutely compact.

Conversely, let $\Dgo \subset \Cgo$ be a generating subcategory of absolutely compact operators. Let $\Igo$ be the two sided ideal generated by $\Dgo$. By proposition \ref{Prop_GeneratingHered=generatingIdeal} $\Igo$ is still generating and because absolutely compact form a two sided ideal $\Igo$ is still composed\footnote{One can actually easily see that $\Igo$ is the ideal of absolutely compact operator} of absolutely compact operator. Moreover as pointed out in \ref{Cor_eqHA_HI} the category of Hilbert $\Dgo$-modules and Hilbert $\Igo$-modules are equivalents, hence it is enough to check the result for $\Igo$. Let $i$ denote the inclusion of $\Igo$ into $\Cgo$. One has two functors $i_*: \Hcal \Dgo \rightarrow \Cgo$ and $R_{\Igo} : \Cgo \rightarrow \Hcal \Dgo$.

As $\Igo$ is generating one has $i_* R_{\Igo} \simeq Id$.

As $\Igo$ is composed of compact operators, the functor $R_{\Igo}$ is continuous, hence $R_{\Igo}(i_*) \simeq (R_{\Igo} \circ i)_*$ but because $\Igo$ is a two sided ideal one easily check that $R_{\Igo}\circ i$ is the Yoneda embeddings of $\Igo$ into $\Hcal \Igo$ and hence $R_{\Igo}(i_*)$ is equivalent to the identity of $\Hcal \Igo$.

Finally as the equivalence between $\Hcal \Igo$ and $\Hcal \Dgo$ is given by the functor $R_{\Dgo}$ this proves that in general the equivalence between $\Cgo$ and $\Hcal \Dgo$ is also given by $R_{\Dgo}$.
}
}

\section{Hilbert modules over toposes}
\label{section_Hilbert_bundle}

\subsection{The $C^{*}$-category of Hilbert modules over a topos}
\label{subsection_HB1}

\blockn{In this subsection, we will consider a Grothendieck topos $\Tcal$ and $\Ccal$ a $C^{*}$-algebra object of $\Tcal$, or even $\Cgo$ a small $C^{*}$-category objects of $\Tcal$. More precisely, a $C^{*}$-algebra $\Ccal$ of $\Tcal$ is an object of $\Tcal$ endowed with structures making it internally into a $C^{*}$-algebra. In this situation we will say that $\Ccal$ is a  $C^{*}$-algebra over $\Tcal$. More generally, a $C^{*}$-category object $\Cgo$ of $\Tcal$ is the data of an object $|\Cgo|$ of $\Tcal$ together with an object $\Cgo$ of $\Tcal/(|\Cgo| \times |\Cgo|)$ and all the additional structure such that internally $\Cgo$ satisfies the axioms for being a $C^{*}$-category with $|\Cgo|$ its set of objects }

\blockn{Our goal is to introduce the category of Hilbert $\Cgo$-modules or of Hilbert $\Ccal$-modules over $\Tcal$ and prove that it is a $C^{*}$-category. This category will be denoted by $\Hcal_{\Tcal} \Cgo$ (or $\Hcal_{\Tcal} \Ccal$). This fact is not new at all, but we could not found a proof of it in the literature so we decided to included it here. }

\blockn{It is known that when $\Tcal$ is the topos of sheaves over a reasonable\footnote{How reasonable the space of $X$ has to be depends on the definition of semi-continuous fields you are using. If they are defined as in \cite{burden1979banach} or \cite{mulvey1980banach} then no hypothesis on $X$ are required.} topological space $X$ then a $C^{*}$-algebra $\Ccal$ over $X$ (over the topos $\sh(X)$) is automatically the sheaf of locally bounded continuous sections of a uniquely defined semi-continuous fields of $C^{*}$-algebras over $X$ (see \cite{burden1979banach}), and that this induces an equivalence between the category of $C^{*}$-algebras over $\sh(X)$ and the category of semi-continuous fields of $C^{*}$-algebras over $X$. This has been proved in the case of Banach space in \cite{mulvey1980banach}, see also \cite{burden1979banach}, and apply all the same to all kind of ``Banach space with structure'', like $C^{*}$-algebra, Hilbert spaces, Hilbert module over a $C^{*}$-algebra etc... Hence this is a very reasonable notion of ``(semi)-continuous fields'' of $C^{*}$-algebras over $\Tcal$.
Asking the continuity of the field is the same as asking that (internally in $\Tcal$) the norm of every element of $\Ccal$ is a continuous real number (see subsection \ref{subsection_genAnaprelim}). We will not explicitly use any of these remark in the present paper, but this is what motivates the study of $C^{*}$-algebras over toposes and hence is the main reason why we have developed everything in the present paper in a constructive framework. Also in order to help the intuition of the reader familiar with continuous fields we will for most of our constructions explain what they mean in terms of semi-continuous fields in the case where $\Tcal$ is a topos of sheaves over a reasonable topological space. These facts will not be proved nor use anywhere in the article (but most of them are direct consequence of the result on semi-continuous fields of Banach spaces).
For example a $C^{*}$-category over a topological space $X$ would be a sheaf $|\Cgo|$ over $X$ and a semi-continuous fields of Banach spaces $\Cgo$ over the etale space of $|\Cgo|\times_X |\Cgo|$ such that for any $x \in X$ the fiber $\Cgo$ and $|\Cgo|$ over $x$ is a $C^{*}$-category and all the operations (composition, $*$, identity) are continuous on the total space.
}

\block{In this situation, a Hilbert $\Ccal$-module or Hilbert $\Ccal$-module over $\Tcal$ is just an object of $\Tcal$ endowed with some structure making it into (internally) a right Hilbert $\Ccal$-module. More generally, a Hilbert $\Cgo$-module over $\Tcal$ is an object of $\Tcal_{/|\Cgo|}$ endowed with a structure making it internally into a Hilbert $\Cgo$-module.
If $H$ and and $H'$ are two Hilbert $\Cgo$-module (or $\Ccal$-module) over $\Tcal$, a globally bounded operator from $H$ to $H'$ is the a map $f:H \rightarrow H'$ such that:

\begin{itemize}

\item $f$ is internally a bounded operator between the Hilbert $\Cgo$-module $H$ to $H'$. In particular $f$ has an adjoint.

\item There exists (externally) a constant $K$ such that internally $\Vert f \Vert <K$.

\end{itemize}

The infimum of all such constants $K$ (as a upper semi-continuous real number) is denoted by $\Vert f \Vert_{\infty}$. Is is described by: $\Vert f \Vert_{\infty} < q$ if and only if there exists $q' <q$ such that internally in $\Tcal$ one has $\Vert f \Vert <q'$.

In the case where $\Tcal$ is a topos of sheaves over a topological space the first condition assert that $f$ is a locally bounded continuous family of operators, and the second that it is in fact globally bounded, with $\Vert f \Vert_{\infty}$ the supremum of the norm at every point (while $\Vert f \Vert$ corresponds to the semi-continuous function on $X$ giving the norm at every point).
}

\blockn{Before going further we will need the following Lemma:}

\block{\label{Lem_sectionarecomplete}\Lem{Let $\Tcal$ be a topos, $B$ be a Banach space object in $\Tcal$. We denote by $\Gamma_b(B)$ the set of global sections of $B$ which are externally bounded endowed with the norm $\Vert x \Vert_{\infty}$ given by the infimum of the external bound of the norm of $x$ as above. Then $\Gamma_b(B)$ is a Banach space.}

\Dem{The algebraic properties and the axiom of the norm pass to sections without difficulties, so the only non trivial thing is that it is complete. We will use the ``Cauchy approximation'' based notion of completeness. Let $A_n$ be a Cauchy approximation on $\Gamma_b(B)$, let $i_n: A_n \hookrightarrow \Gamma_b(B) \hookrightarrow \Gamma(B) = p_* (B)$ be the natural inclusion, where $p$ is the the canonical geometric morphism $\Tcal \rightarrow \{ *\}$. By adjunction one has a morphism $j_n: p^{*} A_n \rightarrow B$ we denote by $B_n$ its image. Because the object of natural number in a Grothendieck topos is an infinite coproduct of copies of the terminal object, an external indexing a family $B_n$ of object is the same as an internal indexing and proving something for all $n$ internally is the same as proving it for each external $n$, hence internally in $\Tcal$ one has a family $B_n$ of subobject of $B$, which is going to be a Cauchy filter: each $B_n$ is inhabited because each $A_n$ was, the inclusion $B_{n+1} \subset B_n$ follow from the same inclusion for the $A_n$ and as for all $x,y \in A_n$ one has $\Vert x - y \Vert_{\infty}<\frac{1}{n}$ one has internally that for all $x,y \in p^{*}A_n, \Vert i(x)-i(y) \Vert < 1/n$ and hence $B_n$ has the same properties because it is the image of $A_n$ in $B$.

This proves that internally $B_n$ converges to a limit $b \in B$ because the limit of a Cauchy approximation is unique $b$ is given by a global section also denoted $b$. If we pick an element $a \in A_n$ one has internally, because $i(a) \in B_n$ and $B_n$ converge to $x$ that $\Vert b -a \Vert \leqslant \frac{1}{n}$ hence $b$ is in $\Gamma_b(B)$ and for each $n$ and each $a \in A_n$, one has $\Vert b- a \Vert_{\infty} \leqslant \frac{1}{n}$, which proves that $A_n$ converge to $b$ and concludes the proof.

}

}

\block{Note that applying this result internally show that if $f$ is a geometric morphism between two topos $f:\Tcal \Rightarrow \Ecal$ then one has a functor $f_{*,b}$ from Banach space over $\Tcal$ to Banach space over $\Ecal$ defined by taking internally in $\Ecal$ the bounded global sections over $\Tcal$. This functor is easily seen to be right adjoint to the pullback of Banach spaces $f^{\sharp}$ which can be defined either as the pullback of the localic completion, or the Hausdorff completion of the pullback as semi-normed space as mentioned in section \ref{subsection_genAnaprelim}.}

\block{\Prop{The category $\Hcal_{\Tcal} \Cgo$ of Hilbert $\Ccal$-modules over $\Tcal$ and globally bounded operators between them is a $C^{*}$-category with:

\begin{itemize}

\item Addition, complex multiplication and adjunction are given by the corresponding internal operations.

\item the norm is the norm $\Vert \_ \Vert_{\infty}$.

\end{itemize}

}

The first point mean that (in the case of a topological space) these operations are computed pointwise. Of course this also applies to $\Hcal_{\Tcal} \Ccal$.

\Dem{The fact that it is a $\mathbb{C}$-linear category follow easily from the analogous internal statement, moreover $\Vert \_ \Vert_{\infty}$ is easily checked to be a norm compatible to addition and composition from the fact $\Vert \_ \Vert$ is such a norm internally. The strong $C^{*}$-inequality of proposition \ref{Prop_alternatedDefCstarcat} is also easily verified simply because it is also true internally, hence (because of proposition \ref{Prop_alternatedDefCstarcat}) the only thing that remains to be proved is the completeness of the morphism space. But the set of globally bounded operators from $H$ to $H'$ can be written as $\Gamma_b(Op(H,H'))$ where $Op(H,H')$ is the (internal) Banach space of operators from $H$ to $H'$ endowed with the operator norm, and hence completeness follow from lemma \ref{Lem_sectionarecomplete}.
}

The exact same argument can be extended to all kind of structures that defines an internal $C^{*}$-category: like an internal pro-$C^{*}$-algebra for example.

}

\block{Finally, it is frequent that, given a $C^{*}$-category over a topos $\Tcal$ there is a $C^{*}$-algebra $\Ccal$ over $\Tcal$ such that $\Hcal_{\Tcal} \Cgo$ is equivalent to $\Hcal_{\Tcal} \Ccal$. For example, if $\Cgo$ has a small decidable\footnote{An object $X$ is said to be decidable if internally for $x,y \in X$, $x=y$ or $x \neq y$.} set of objects then one can construct a Hilbert $\Cgo$-module $\Hcal = \bigoplus_{c \in |\Cgo|} \Yon_c$, more precisely, $\Hcal$ has one generator $e_c$ for each $c \in |\Cgo|$ and $\scal{e_c}{e_{c'}}=0$ if $c \neq c'$ and $1_c$ if $c=c'$ (if $\Cgo$ was non unital one can adapt using generators of the $\Yon_c$ instead). Then it is easy to see internally that restricting to the algebra $\Ccal$ of compact operators of $\Hcal$ induce an equivalence between the category of Hilbert $\Cgo$-module and Hilbert $\Ccal$-module, which clearly induce an (external) equivalence between the category $\Hcal_{\Tcal} \Cgo$ and $\Hcal_{\Tcal} \Ccal$.

More generally, if $|\Cgo|$ admit a covering by a decidable object $X$ one can perform a similar construction with $\Hcal=\bigoplus_{x \in X} \Yon_x$ and obtains the same result, and for a large number of example of topos (for examples, all topological spaces) any sheaf admit a covering by a decidable sheaf, hence over those toposes this construction can always be performed and there is very little gain to consider $C^{*}$-categories instead of $C^{*}$-algebras. But in the more general situation $C^{*}$-category are more flexible, for example if $\Tcal$ is any topos and $v$ is an element of $H^{2}(\Tcal, \U)$ where $\U$ is the commutative group object of continuous complex number of module $1$ then there is a notion of $v$-twisted Hilbert space\footnote{It is a generalization of $v$-twisted unitary representations of a group $G$ for $v$ a $2$-cocycle on $G$ with value in complex number of module $1$.} over $\Tcal$ which corresponds to Hilbert modules over a $C^{*}$-category over $\Tcal$ naturally attached to $v$ (this $C^{*}$-category is well defined up to unique external Morita equivalence and is internally Morita equivalent to $\C$). If we have additional assumption on $\Tcal$ like the fact that any object can be covered by a decidable object then this $C^{*}$-category over $\Tcal$ can be replaced by a $C^{*}$-algebra over $\Tcal$ but this not always possible.  }

\block{Note that in this section we have assumed that $\Cgo$ was a small $C^{*}$-category over $\Tcal$. This was only to keep the exposition simple but we can get rid of this assumption and we will. A general $C^{*}$-category over $\Tcal$ will be exactly the same except that the object $|\Cgo|$ can be a large sheaf (see appendix \ref{appendiw_foundation}), and $\Hcal_{\Tcal} \Cgo$ is now the category of Hilbert $\Cgo$-module over $\Tcal$ which are small in the sense that (externally) there exists a small object $X$ of $\Tcal$ with a map from $X \times X$ to $\Cgo$ which is internally a set of generator of our Hilbert module. The proof given above that $\Hcal_{\Tcal} \Cgo$ is a $C^{*}$-category only rely on the fact that for each pair of objects one has a small sheaf over $\Tcal$ of operators between them, that this form internally a $C^{*}$-category and that morphisms in $\Hcal_{\Tcal} \Cgo$ corresponds to bounded global sections of this sheaf, hence it works exactly the same when $\Cgo$ is no longer small. 

It is important that the set of generators is given externally so that we can actually talk about the class of small Hilbert module (it is the class of sets of generator) while the 'class of Hilbert modules for which there exists internally a set of generators' does not clearly exists in our framework. Moreover it is easy to pass from a Hilbert module which is internally small to a Hilbert module with a externally given by a set of generators: internally small mean that one has a set of generators for $S \times \Hcal$ in $\Tcal_{/S}$ for some inhabited $S$ then forgetting the map to $S$ and applying the projection $\Hcal \times S \rightarrow S$ gives us a set of generators for $\Hcal$ in $\Tcal$.
 }

\subsection{Completeness for the $C^{*}$-category of Hilbert modules over a topos}

\block{Let $\Tcal$ be a topos, $\Cgo$ be a $C^{*}$-algebra over $\Tcal$. In the previous sub-section we defined the $C^{*}$-category $\Hcal_{\Tcal} \Cgo$ of Hilbert $\Cgo$-modules over $\Tcal$ and globally bounded operators between them. In this sub-section we will study its properties from the point of view of pre-complete $C^{*}$-categories. Although we are mainly interested in the case where $\Cgo$ is just a $C^{*}$-algebra, it is actually helpful to consider the case of a $C^{*}$-category.
}

\block{Let $F:\Dgo \rightarrow \Hcal_{\Tcal} \Cgo$ be a $C^{*}$-functor from a small $C^{*}$-category, and let $H$ be a Hilbert $\Dgo$-module. We want to define $F_* H \in \Hcal_{\Tcal} \Ccal|$. Let $p : \Tcal \rightarrow *$ the canonical geometric morphism to the point (the topos of sets), one can pullback $\Dgo$ and $H$ to a $C^{*}$-category $p^{\sharp}\Dgo$ and a Hilbert module $p^{\sharp} H$ over it in $\Tcal$: one pullback the set of objects ($p^{*}$) and then pullback all the Banach spaces involved in the $C^{*}$-category and the Hilbert module ($p^{\sharp}$). Using the adjunction between $p^{\sharp}$ and $\Gamma_b$ one obtains a $*$-functor $p^{\sharp}\Dgo \rightarrow \Hcal \Cgo$ internally in $\Tcal$ and one can use the pre-complete structure of $\Hcal \Cgo$ internally in $\Tcal$ to construct a well defined $\Cgo$-module $F_* H$ which is internally small and hence can be seen as an element of $\Hcal_{\Tcal} \Cgo$. 

This object $F_* H$ can be described more explicitly. Let $F:\Dgo \rightarrow \Hcal_{\Tcal} \Cgo$ be a $*$-functor and $H$ a Hilbert $\Dgo$-module. Then the Hilbert $\Cgo$-module $H' = F_* H$ admit a cone structure given by morphisms $H'_x$ from $F(d)$ to $H'$ for each $d \in |\Dgo|, x \in H(d)$ they satisfies $(H'_x)^{*} H'_{x'} = F(\scal{x}{x'})$ and $H'$ is internally the closure of the image of all those $H_x$ hence it can be written as the closure of a pre-Hilbert module generated by elements of the form $H'_x(v)$ for $d \in \Dgo x \in H(d)$ and $v \in F(d)$.

It is clear that this is functorial in $\Dgo$ (because the pullback to $\Tcal$ and the application of internal pre-complete structure are) and that if $F$ is unital and $H$ representable then $F_*H$ is just $F(H)$ (either using the explicit description or because those properties are also preserved by pullback and that this is true internally) hence this indeed defines a pre-complete structure on $\Hcal_{\Tcal} \Cgo$.

\Prop{$\Hcal_{\Tcal} \Cgo$ is a pre-complete $C^{*}$-category for the structure described above.}

In general this structure is an anafunctor.

One could also apply a similar argument to any kind of internally a pre-complete $C^{*}$-category: for example, if $\Ccal$ is a pro-$C^{*}$-algebra over $\Tcal$ then the category of Hilbert $\Ccal$-modules over $\Tcal$ is also going to be a pre-complete $C^{*}$-category for a similar structure.
}

\block{\label{Prop_weakdensityinHTC}\Prop{Let $\Acal \subset \Hcal_{\Tcal} \Cgo (X)$ a hereditary sub-algebra of endomorphisms of an object $X \in |\Hcal_{\Tcal} \Cgo|$. Then following conditions are equivalent:

\begin{itemize}

\item $\Acal$ is weakly dense for the pre-complete structure of $\Hcal_{\Tcal} \Cgo $.

\item The hereditary subalgebra $\Acal'$ of $\Hcal \Cgo (X)$ generated internally by $\Acal$ is (internally) weakly dense, i.e. satisfies internally the following equivalent conditions of proposition \ref{Prop_weakdensityinHC}.

\item If $X' \subset X$ is a sub-module such that any element of $\Acal$ can be factored (as a linear map) into $X \rightarrow X' \hookrightarrow X$ where the second arrow is the inclusion, then $X'=X$

\item There is a directed increasing net $(a_{\lambda})_{\lambda \in \Lambda}$ of positive elements of norm smaller than $1$ in $\Acal$ which internally converges to $1$ pointwise, i.e. such that internally, for all $c \in |\Cgo|$ and $x\in X(c)$ the net $a_{\lambda} x$ converges to $x$.

\end{itemize}

Moreover, these conditions implies that $\Acal$ is essential and hence $\Hcal_{\Tcal} \Ccal$ satisfies $(C1)$.

}

\Dem{ $\Acal X$, for the pre-complete structure of $\Hcal_{\Tcal} \Cgo$ is simply (by definition of the pre-complete structure) $\Acal' X$ computed internally, i.e. it can be describe internally as the closure of $\{ a.x \in X(c) | x \in X(c)  a \in p^{*}\Acal \}$ for any $c \in |\Cgo|$, and it is isomorphic to $X$ in a way compatible to the cone structure if and only if it is (internally) equal to $X$ as a subspace (the compatibility to the cone structure exactly says that the element $a.x$ of $\Acal.X$ should corresponds to the element $a.x$ of $X$). Hence one obtains the equivalence between the first and the second condition. Moreover, the action of $\Acal$ can be factored into a submodule $X'$ of $X$ if and only if $X'$ contains the submodule defined internally as the closure of $\{ a.x \in X | x \in X  a \in p^{*}\Acal \}$, hence the first and the third condition are also equivalent.

If one has a net as in the last condition, then in particular, internally, any element of $X$ can be approximated by elements of the form $ax$ hence it implies the first two conditions. Conversely, if one assume the first two conditions one can pick for $(a_{\lambda})$ an approximate unit as for example the one constructed in proposition \ref{Prop_approximate_unit}, and it will be true internally that for any $a \in p^{*} \Acal$, one has $a_{\lambda} a \rightarrow a $ and internally one has that for any $c \in |\Cgo|, x \in X(c)$ and any $\epsilon>0$  there exists $x' \in X(c)$ and $a \in p^{*}\Acal$ such that $\Vert ax'-x \Vert<\epsilon$ (this is what the second condition claims), using that $a_{\lambda}a \rightarrow a$, there exists a $\lambda_0$ such that for any $\lambda \geqslant \lambda_0$, $\Vert a_{\lambda} a -a \Vert < \epsilon$, hence for any $\lambda \geqslant \lambda_0$ one has:

\[
\begin{array}{c c c c c c c c }
 \Vert a_{\lambda}x -x \Vert & \leqslant & \Vert a_{\lambda} (x - ax') \Vert & + & \Vert a_{\lambda} a x' -ax' \Vert & + & \Vert ax' -x \Vert \\
 & \leqslant & \epsilon &+ & \Vert x \Vert \epsilon  & + & \epsilon
\end{array} 
 \]

and hence concludes the proof of the equivalence of the four conditions.

Finally, if $\Acal$ is weakly dense, and $f$ is an operator such that $fa = 0$ for all $a \in \Acal$ then one can prove internally that $f=0$ (because for any element of the form $ax, fax=0$ and such element are dense) and hence $f=0$, this proves that $\Acal$ is essential and hence that $\Hcal_{\Tcal} \Cgo$ satisfies $(C1)$.

}

}

\block{\Prop{The pre-complete $C^{*}$-category $\Hcal_{\Tcal} \Cgo$ satisfies all conditions $(C1)$,$(C2)$,$(C3)$ and $(C4)$.}

\Dem{By proposition \ref{Prop_impbetween_Ci} and the fact proved just above that this category satisfies $(C1)$, it suffices to prove that is satisfies $(C4)$. So let $f:A \rightarrow L(B)$ and $g: B \rightarrow L(A)$ a pre-adjoint pair (as in definition \ref{def_extension_prop}) in $\Hcal \Hcal_{\Tcal} \Cgo$. The two objects $L(A)$ and $L(B)$ are Hilbert $\Cgo$-module over $\Tcal$ and we need to construct an operator $\tilde{f}$ between them. 

Let $v$ be an element of $A(X)$ for some $X \in |\Hcal_{\Tcal} \Cgo|$, in particular, $v$ induces an operator from $X$ to $L(A)$ (we identifies $L(A)$ with its image by the Yoneda embeddings), for any (internal) $c\in|\Cgo|, x \in X$ it is natural to define $\tilde{f}(v(x))$ as $f(v)(x)$ ($f(v)$ is an element of $L(B)(X)$, hence an operator from $X$ to $L(B)$). As it is not clear at this point that $\tilde{f}(v(x))$ only depends on $v(x)$, we will denote\footnote{On a very formal level, this means that for each external $v$ and $X$ one has a function $\tilde{f}(v, \_)$ between $X$ and $L(A)$ as object of $\Tcal$ defined internally by the given formula.} $\tilde{f}(v,x) := f(v)(x)$ for now. Similarly, for any $w \in B(Y)$ and (internally) $ y \in Y $ one defines $\tilde{g}(w,y) := g(w)(y) \in L(A)$. One can then observe that:

\[ \scal{w(y)}{ \tilde{f}(v,x)} = \scal{w(y)}{f(v)(x)}= \scal{y}{w^{*}f(v) (x)} \]

Now, because of the pre-adjunction formula, $w^{*}f(v)$ is, as an operator from $X$ to $Y$, equal to $g(w)^{*}v$ hence:

\[ \scal{w(y)}{ \tilde{f}(v,x)} = \scal{y}{g(w)^{*} v (x)} = \scal{g(w)(y)}{v(x)} = \scal{\tilde{g}(w,y)}{v(x)}\]

By definition of the pre-complete structure on $\Hcal_{\Tcal} \Cgo$, elements of the form $w(y)$ and $v(x)$ (with $X$ and $Y$ allowed to vary) are internally dense in $LA$ and $LB$, hence the above equality proves that $\tilde{f}(v,x)$ and $\tilde{g}(w,y)$ only depends on $v(x)$ and $w(y)$ and that they are adjoint to each other when they are defined. Finally from the fact they are densely defined and bounded (by the norm of $f$ and $g$) one can deduce that they extend into operators between $L(A)$ and $L(B)$ and hence this concludes the proof.
}

}

\block{\label{Prop_QcontaininHTC}\Prop{Let $Y$ be an object of $\Hcal_{\Tcal} \Cgo$ and $\Ago$ be a hereditary subcategory of $\Hcal_{\Tcal} \Cgo$, then $Y$ is quasi-contained in $\Ago$ in the sense of definition \ref{Def_quasicontain} if and only if it satisfies the following property:

If $Y' \subset Y$ is a sub-module of $Y$ such that (externally) any operator $f$ with value in $Y$ such that $f^{*}f \in \Ago$ can be factored (as a linear map) into the inclusion of $Y'$ in $Y$, then $Y'= Y$.}

\Dem{This is essentially the same as the equivalence between the first and the third point of proposition \ref{Prop_weakdensityinHTC}: $Y$ is quasi-contained in $\Ago$ if the two sided ideal of $\Ago$-compact operator is weakly dense in $Y$ and the corresponding submodule is essentially by construction the smallest sub-module that contains the image of all the operators $f :a \rightarrow Y$ such that $f^{*}f \in \Ago$, moreover, it is isomorphic to $Y$ in a way compatible to the cone structure if and only if it is equal to $Y$ as a submodule. }
}

\subsection{The topos theoretic Green-Julg theorem}

\blockn{In this subsection, we want to prove a topos theoretic analogue of the Green-Julg theorem stating that for a certain class of toposes the category of Hilbert space over $\Tcal$ (and more generally categories of the form $\Hcal_{\Tcal} \Cgo$) is a category of Hilbert modules over a $C^{*}$-algebra attached to $\Tcal$ (more generally, a $C^{*}$-algebra $\Cgo \rtimes \Tcal$). This result will of course rely on the characterization of such categories given in proposition \ref{Th_caracofHmod}, and hence we first need to find criterion for an arrow in $\Hcal_{\Tcal} \Ccal$ to be absolutely compact in the sense of definition \ref{Def_Abscpt}, and then prove that those categories have enough such operators.}

\blockn{A section of the sheaf of upper semi-continuous real numbers of a topos $\Tcal$ is the same as a geometric morphism from the topos $\Tcal$ to the locale of upper semi-continuous numbers. In classical mathematics, this is the space of real numbers with the topology of upper semi-continuity, i.e. where the only open subsets are the $(-\infty, a)$, in constructive mathematics it is defined as the classifying space for the theory of upper semi-continuous real number which is a locale because this theory only have propositions in its language (See \cite{sketches} D4.7, especially lemma D4.7.2 and corollary D4.7.3). In particular it is the same as a morphism from the localic reflection of $\Tcal$ to the local of upper semi-continuous real number.}

\block{\Def{A global section $f$ of the sheaf of positive upper semi-continuous real numbers of a topos $\Tcal$ is said to be \emph{zero at infinity} if for any rational number $q>0$ the closed complement of the open subspace $`f < q'$ of the localic reflection of $\Tcal$ is compact. }

By the above discussion, $f$ can be seen as a positive upper semi-continuous function on the localic reflection of $\Tcal$ and hence $f <q$ indeed defines an open subspace of this locale.
}

\block{\label{Prop_compactnessHTC}\Prop{Let $\Tcal$ be a topos, $\Cgo$ a $C^{*}$-category over $\Tcal$ and $f : H \rightarrow H'$ a globally bounded operator between two Hilbert $\Cgo$-modules over $\Tcal$. Assume that:

\begin{itemize}

\item $f$ is internally in $\Tcal$ a $\Cgo$-compact operator.

\item $\Vert f \Vert$, seen as a section of the sheaf of positive upper semi-continuous real numbers is zero at infinity.

\end{itemize}

Then $f$ is absolutely compact in the pre-complete $C^{*}$-category $\Hcal_{\Tcal} \Ccal$.

}

\Dem{Let $f:H' \rightarrow H$ an operator satisfying the conditions of the proposition, and let $(a_{\lambda})_{\lambda \in \Lambda}$ be a bounded (with $K$ such that $\Vert a_{\lambda} \Vert <K$ for all $\lambda$) net of operators from $H$ to $Z$ in $\Hcal_{\Tcal} \Cgo$ which converges to an operator $a:H \rightarrow Z$ in $\Hcal_{\Tcal} \Cgo$ in the sense of internal pointwise convergence, i.e. such that internally one has : $\forall c \in |\Cgo|, \forall h \in H(c), \Vert a_{\lambda} h - a h \Vert \rightarrow 0$. We will prove that $a_{\lambda} \circ f $ tends to $a \circ f $ in norm.

We will first prove that $a_{\lambda} \circ f $ tends to $a \circ f $ internally in norm. Internally, as $f$ is compact it can be approximated by ``finite rank'' operator, let $h$ be a finite rank operator such that $\Vert f - h \Vert < \epsilon$, the pointwise convergence of $a_{\lambda}$ automatically implies the norm convergence of $a_{\lambda} \circ h$ (it suffice to apply the convergence the finite set of vector that appears in the definition of $h$), but as $(a_{\lambda})$ is bounded, one has $ \Vert a_{\lambda} \circ f - a_{\lambda} \circ h \Vert < K \epsilon$, hence this also implies that $a_{\lambda} \circ f \rightarrow a \circ f$ internally in norm.

We will now prove the external convergence of $a_{\lambda} \circ f$ to $a \circ f$. Let $\epsilon>0$ be a rational number, because $a_{\lambda} \circ f$ converge internally in norm, one has a covering $(U_i)_{i \in I}$ of the terminal object, such that for each $i \in I$ there exists a $\lambda_0 \in \Lambda$ such that for all $\lambda \geqslant \lambda_0$ $\Vert a_{\lambda} \circ f - a \circ f \Vert < \epsilon$ on $U_i$. The $U_i$ naturally arise as arbitrary object, but on can replace them by their image on the terminal object and hence assume that they are sub-object of the terminal object. Because $\Vert f \Vert $ is zero at infinity, the closed complement of the subterminal object on which $\Vert f \Vert < \frac{\epsilon}{2K}$ is compact, hence there exists a finite subset $J \subset I$ such that the terminal object is covered by the $(U_j)_{j \in J}$ and the subobject $W$ on which $\Vert f \Vert < \frac{\epsilon}{2K}$. One can then fix a $\lambda_0$ such that for each $\lambda \geqslant \lambda_0$,  $\Vert a_{\lambda} \circ f - a \circ f \Vert < \epsilon$ on $U_j$ for all $j \in J$, and this also holds on $W$ because $\Vert a \circ f - a_{\lambda} \circ f \Vert \leqslant 2K \Vert f \Vert < \epsilon$. Hence $\forall \lambda \geqslant \lambda_0$, $\Vert a_{\lambda} \circ f - a \circ f  \Vert_{\infty}<\epsilon$ which show as claimed above that $a_{\lambda} \circ f \rightarrow a \circ f$ in $\Hcal_{\Tcal} \Cgo$.

We can now conclude on the absolute compactness of $f$: Operators satisfying the conditions of the proposition clearly form a closed two sided ideal of $\Hcal_{\Tcal} \Cgo$. Moreover if $\Acal \subset \Hcal_{\Tcal} \Cgo (X)$ is a weakly dense hereditary sub-algebra, and $f:X \rightarrow X$ satisfies the condition of the proposition, then for any approximate unite $a_{\lambda}$ of $\Acal$ one has $a_{\lambda} f \rightarrow f$ hence $f$ belong to the right ideals generated by $\Acal$ and $ff^{*} \in \Acal$. As $f^{*}$ satisfies the same properties, $f^{*}$ is also in the right ideal generated by $\Acal$, hence $f^{*} f$ and $ff^{*}$ are both in $\Acal$ which implies that $f \in \Acal$ as $\Acal$ is hereditary. And this implies that those operators are absolutely compact by the definition of absolute compactness.

}
}

\blockn{We will also need the following lemma:}

\block{\label{Lem_finitegeneration}\Lem{Let $\Tcal$ be a topos which is separated and locally decidable, then it admit a generating family $(X_i)$ of objects such that for each $i$ there exists a sub-object $U_i$ of the terminal object of $\Tcal$ such that $X_i$ is over $U$ and is finite and decidable over $U$.}

\Dem{The proof of this is essentially the same as the main theorem of \cite{henry2015finitness}. Let $\Lcal$ be the localic reflection of $\Tcal$, we can see $\Tcal$ as a topos over $\Lcal$, and internally in $\Lcal$, $\Tcal$ is hyperconnected by definition, separated because of proposition $II.2.3$ of \cite{moerdijk2000proper} and locally decidable because of lemma $5.1$ of \cite{henry2015finitness}, hence one can apply theorem $4.6$ of \cite{henry2015finitness} to it and conclude that (internally) in $\Lcal$, $\Tcal$ admit a generating family of finite objects.

Let $Y' \subset Y$ be a subobject of an object of $\Tcal$, and assume that any map from an object $X$ as described in the proposition (finite decidable over a sub-terminal object) to $Y$ factor into $Y'$, in order to conclude we just have to show that $Y' =Y$. We will first show that is true internally in $\Lcal$ that any map from a finite decidable object of $\Tcal$ to $Y$ factor into $Y'$: indeed this is equivalent to the fact that for any $U$ object of $\Lcal$, for any map $f:V \rightarrow U$ in $\Tcal$ which is fiberwise finite and decidable any map from $V$ to $Y$ factor into $Y'$ and this can be proved by covering $U$ by subterminal objects.

Hence, as it is internally true in $\Lcal$ that $\Tcal$ is generated by finite decidable object, this proves that $Y' = Y$ and concludes the proof.   }

}

\blockn{We can now prove our main theorem:}

\block{\label{main_Th1}\Th{Let $\Tcal$ a topos such that:

\begin{itemize}

\item $\Tcal$ is separated.

\item $\Tcal$ is locally decidable.

\item The localic reflection of $\Tcal$ is locally compact and completely regular\footnote{Assuming the axiom of depend choice, the hypothesis completely regular can be removed as it follows from locally compact and separated by the localic version of Urysohn's lemma, see \cite{picado2012frames} sections V.5, VII.2 and XIV.6.2 or \cite{henry2014localic} section 2.6 together with \cite{sketches} C.3.2.10.}.

\end{itemize}

Let also $\Cgo$ be a $C^{*}$-category (possibly non-small) over $\Tcal$ then there exists a $C^{*}$-category $\Cgo \rtimes \Tcal$ (in the topos of sets) well defined up to Morita equivalence such that the category of Hilbert $\Cgo$-modules over $\Tcal$ is equivalent to the category of Hilbert $\Cgo \rtimes \Tcal$-modules, i.e.: 

\[ \Hcal_{\Tcal} \Cgo \simeq \Hcal (\Cgo \rtimes \Tcal) \]

Moreover, $\Cgo \rtimes \Tcal$ can be chosen canonically to the $C^{*}$-category $\Kcal_{\Tcal} \Cgo$ of Hilbert $\Cgo$-modules over $\Tcal$ and operators between them which satisfies the conditions of proposition \ref{Prop_compactnessHTC}. 
}

\Dem{By proposition \ref{Prop_compactnessHTC}, the category $\Kcal_{\Tcal} \Cgo$ as defined in the theorem form a two sided ideal of $\Hcal_{\Tcal} \Cgo$ of absolutely compact operators, hence by theorem \ref{Th_caracofHmod} it is enough to show that $\Kcal_{\Tcal} \Cgo$ is generating, and hence by proposition \ref{Prop_generating=dense} and the fact that $\Hcal_{\Tcal} \Cgo$ satisfies $(C1)$ (proposition \ref{Prop_weakdensityinHTC}) it is enough to show that for every object $H \in |\Hcal_{\Tcal} \Cgo |$ the algebra $\Kcal_{\Tcal} \Cgo(H)$ is weakly dense in $\Hcal_{\Tcal} \Cgo(H)$, and this will be done using the second criterion of proposition \ref{Prop_weakdensityinHTC}.

So, let $H$ be a Hilbert $\Cgo$-module over $\Tcal$, and let $H'$ be a sub-module of $H$ such that any element of $\Kcal_{\Tcal} \Cgo(H)$ have its image included in $H'$, and we need to prove that $H=H'$ in order to conclude. We will use the fact (lemma \ref{Lem_finitegeneration}) that $\Tcal$ has a generating family of objects $X_i$ such that for each $i$, $X_i$ is finite and decidable over a sub-terminal object. We fix one of these objects, and we assume that one has a map from $X$ to $H$ (here we see $H$ as an object over $|\Cgo|$, so it means internally that for any $x \in X$ one has an object $c_x \in |\Cgo|$ and an element $h_x \in H(c_x)$) , we only need to prove that this map factor into $H'$.

Because the localic reflection of $\Tcal$ is locally compact and completely regular, $U_i$ admit a covering by sub-object $V_j$ such that for each $j$ one has a (closed) compact sub-locale $K_j$ of the localic reflection such that $V_j \subset K_j \subset U_i$ and a function $f_j$ from the localic reflection to the locale of complex number such that $f_j$ is equal to $1$ on $V_j$ and has its support inside $K_j$.

We fix one of this $V_j$, and we will construct an endomorphism $k_j$ of $H$ defined as:

\[ k(v) = f_j \sum_{x \in X} h_x \scal{h_x}{v} \]

This is actually a definition by case: as $K_j \subset U_i$ the terminal object admit a covering by $U_i$ and $K_j^{c}$ the open complement of $K_j$, seeing those two subobjects of the terminal object as proposition, one has internally $U_i$ or $K_j^{c}$. Assuming $K_j^{c}$ then $f_j=0$ hence $k(v)=0$ is well defined, assuming $U_i$, then $X$ is a finite decidable set so it is legitimate to take a sum indexed by $X$, moreover if $v \in H(c)$ then $\scal{h_x}{v}$ is an element of $\Cgo(c,c_x)$ hence $h_x \scal{h_x}{v}$ is an element of $H(c)$ hence it is well defined, and when both $U_i$ and $K_j^{c}$ holds the two expressions agrees so this indeed defined an endomorphism of $H$. This endomorphism is obviously internally a finite rank (hence compact) bounded operator and as it is non-zero only on a compact it is obviously globally bounded and with zero norm at infinity, hence it is an element of $\Kcal_{\Tcal} \Cgo (H)$, hence it factor in $H'$.

Assume (internally) $V_j$, hence $f_j=1$, $X$ is a finite decidable object and $k(v) = \sum_{x \in X} h_x \scal{h_x}{v}$, $k$ can be written as $g.g^{*}$  where $g$ is the map from $\oplus_{x \in X} Y_{c_x}$ to $H$ which send each component to $H$ by multiplication by $h_x$, the fact that $gg^{*}$ factor into $H'$ implies that $g$ itself factor in $H'$, hence that $h_x \in H'$. As such $V_j$ cover $U_i$ it proves that for any $x \in X$, $h_x \in H'$, hence that the map from $X$ to $H$ factor into $H'$ and this concludes the proof.

}

}

\block{\label{Prop_smallnessinMainTh}\Prop{Assume that in the previous theorem, $\Cgo$ is a small $C^{*}$-category over $\Tcal$, then $\Tcal \rtimes \Cgo$ can be chosen to be a small $C^{*}$-category. }

\Dem{We will show that in this case the proof of the theorem \ref{main_Th1} only use a set of object in $\Kcal_{\Tcal} \Cgo$, and not the whole category. Fix $I$ an indexing of a set generators $(X_i)$ of $\Tcal$ as in the previous proposition (with $X_i$ finite and decidable over $U_i \subset 1$ for each $i$), for each $i \in I$ and each $c :X_i \rightarrow \Cgo$ one can construct a Hilbert $\Cgo$-module $H_{i,c}$ : if $\Cgo$ is unital it has one generator for each $x \in X$ and the scalar product is given by $\scal{e_x}{e_x'} = 0$ if $x \neq x'$ and $1_{c(x)}$ if $x=x'$ (note that $X_i$ is always decidable), if $\Cgo$ is non unital one can use instead one generator for each $x \in X$ and each endomorphism of $c(x)$.

Then (assuming $f_j$ is real and positive, which is always possible) the operator $k$ which we used in the proof of \ref{main_Th1} can be written as $gg^{*}$ where $g$ is an operator from a $H_{i,c}$ to $H$: take $g(e_x)=(f_j)^{1/2} h_x$ and the same proof by cases as above show that $g$ is a well defined bounded operator (and globally bounded because of compactness) and that $gg^{*} = k$, hence the exact same proof as above show that one can take $\Cgo \rtimes \Tcal$ to be the full subcategory of $\Kcal_{\Tcal} \Cgo$ on the object of the form $H_{i,c}$ and this is a small category.

}
}

\blockn{Finally, one can also deduces that:}

\block{\Th{Assume that the ground topos is boolean, or just that it is true (internally in the ground topos) that every object admit a covering by a decidable object. Then for any topos $\Tcal$ which is separated, locally decidable and whose localic reflection is locally compact and completely regular\footnote{Once again, this hypothesis is automatically true if one assume the axiom of dependant choice in the ground topos.}, and any $C^{*}$-algebra $\Ccal$ over $\Tcal$ there exists a $C^{*}$-algebra $\Ccal \rtimes \Tcal$ (in the ground topos) such that there is an equivalence of category:

\[ \Hcal_{\Tcal} \Ccal \simeq \Hcal (\Ccal \rtimes \Tcal) \]

Moreover $\Ccal \rtimes \Tcal$ (endowed with this equivalence) is unique up to unique\footnote{The Morita equivalence being it self unique up to unique isomorphism} Morita equivalence, and is\footnote{More precisely, any $C^{*}$-algebra giving such an equivalence will be.} isomorphic to the algebra of endormophism which are internally compact and whose norm is zero at infinity of a Hilbert module over $\Tcal$.
}

\Dem{Theorem \ref{main_Th1} and proposition \ref{Prop_smallnessinMainTh} show that one can find a small $C^{*}$-category $\Cgo$ such that $\Hcal_{\Tcal} \Ccal$ is equivalent to $\Hcal \Cgo$, hence if the set of objects of $\Cgo$ can be covered by a decidable object one can construct a Hilbert module over $\Cgo$ as $\bigoplus_{c \in \Cgo} \Yon(c)$ and show that $\Hcal \Cgo$ is equivalent to the category of Hilbert modules over the algebra of $\Cgo$-compact operators of this module, hence this prove the equivalence.

For any other $C^{*}$-algebra $\Dcal$ satisfying this property there will be a canonical equivalence of category between $\Hcal \Dcal$ and $\Hcal (\Ccal \rtimes \Tcal)$ which hence induces a unique morita equivalence between $\Dcal$ and $\Ccal \rtimes \Tcal$. Finally, as it has been remarked above, the category $\Kcal \Cgo$ of $\Cgo$-compact operators is identified in the above equivalence to the category $\Kcal_{\Tcal} \Ccal$, but any $C^{*}$-algebra such that $\Ccal$ such that $\Hcal \Ccal \simeq \Hcal \Cgo$ turns $\Ccal$ into the algebra of compact operator of a Hilbert $\Cgo$-module (the modules corresponding to $\Ccal$ under the equivalence), hence an algebra of endomorphism in $\Kcal_{\Tcal} \Cgo$, i.e. the algebra of operator that are internally compact and whose norm tends to zero at infinity over some generating Hilbert $\Cgo$-module over $\Tcal$.
}
}

\appendix
\renewcommand{\thesubsubsection}{\Alph{section}.\arabic{subsubsection}}
\section{Localic spectrum and positivity}
\label{appendix_positivity}

\blockn{The goal of this appendix is to introduce the notion of spectrum of an element in a Banach algebra or a $C^{*}$-algebra, to prove that for $C^{*}$-algebras if $x \in B \subset A$ then the spectrum of $x$ in $A$ and $B$ are the same and to deduce from this the various properties of positive elements in a $C^{*}$-algebras. For a reader not interested in constructive aspect or willing to put the strong $C^{*}$-inequality as in proposition \ref{Prop_alternatedDefCstarcat} or Corollary \ref{Cor_strongineq} instead of the $C^{*}$-equality in the definition of $C^{*}$-algebras, this appendix can be completely ignored.}

\blockn{It is very natural to expect that in constructive mathematics the correct notion of spectrum of an element $x$ in a $C^{*}$-algebra should be a sub-locale of the locale $\C$ rather than a subset of $\C$. Indeed the spectrum of a commutative algebra is already known to be a locale rather than a topological space, and moreover the locale of complex numbers is always locally compact (which will make our spectrum compact) while the topological space of complex number is locally compact only if it is homeomorphic to the locale of complex number, hence if one defines spectrum as sets they are not going to be compact for their natural topology in general. }

\blockn{Because we are going to work with localic spectrum it will be slightly more natural to work with a localic Banach algebra as in \cite{henry2014localic}. But our result apply all the same to an ordinary Banach algebra: indeed if $A$ is an ordianry banach algebra then its localic completion (see \cite[3.6]{henry2014localic}) is a localic Banach algebra whose points are exactly the elements of $A$. We will denote both the ordinary algebra and its localic completion by the same letter $A$.}

\blockn{Also, in this appendix all the algebras are \emph{assumed to be unital} and the sub-algebras to contains the unit. The case of non unital algebras will be briefly mentioned at the very of the appendix, in \ref{Rk_nonunital}.}

\block{\Prop{Let $A$ be a localic Banach algebra, let $A^{\times}$ be the locale corresponding to $\{(x,y) \in A^{2} | y x = x y =1  \}$, and let $i:A^{\times} \rightarrow A$ be the map which send $(x,y)$ to $x$. Then $i$ is an open inclusion.}

This is a geometric formulation of the fact that the set of invertible is open, but it also contain the fact that the inversion is continuous on the set of invertible elements.

\Dem{ The map from $A^{\times}$ to $A$ is a monomorphism, indeed in terms of generalised elements, if one has $(x,y)$ and $(x,y')$ in $A^{\times}$ then $y = y'xy = y'$, hence (see \cite[Corollary C3.1.12]{sketches}) it suffices to prove that this map is an open map.

Let 
\[ A_{1/2} = \{ (x,y) \in A^{2} | \Vert 1 - xy \Vert < \frac{1}{2} ; \Vert 1 - yx \Vert < \frac{1}{2} \} \]

By continuity of the multiplication it is an open sub-locale of $A^{2}$, moreover as $A$ is locally positive the map $\pi_2 : A \times A \rightarrow A$ is open hence the first component map $A_{1/2} \rightarrow A$ (send $(x,y)$ to $x$) is an open morphism.

There is also a canonical inclusion map $A^{\times} \hookrightarrow A_{1/2}$, we will construct a retraction $A_{1/2} \rightarrow A^{\times}$ of this map over $A$ which will conclude the proof, for example because once the retraction is constructed one can construct the direct image in $A$ of any open subset $U \subset A^{\times}$ by first pulling it back to $A_{1/2}$ and then taking the direct image in $A$, hence the resulting sub-locale of $A$ will be always open. 

In order to construct this retraction, we will work internally in $A_{1/2}$.

In $A_{1/2}$ one dispose of a pair of elements $(x,y)$ such that $\Vert 1 - xy \Vert <1/2$ and $\Vert 1- yx \Vert <1/2$.

Hence one can define an element $c \in A$ by:

\[ c = \sum_{k=0}^{\infty} (1-xy)^{k} \]

and one has $(1-xy)c= c - 1$, hence $xyc=1$ hence defining $u:=yc$ one has an element $ u \in A$ such that $ux=1$, symmetrically one can construct an element $v \in A$ such that $xv=1$ and hence $v=uxv=u$, hence one has constructed internally in $A_{1/2}$ an element $u \in A$ such that $xu=ux=1$, as this construction does not depends on any choice it corresponds to a morphism $A_{1/2} \rightarrow A^{\times}$ which send $(x,y)$ to $(x,u)$, and this is the retraction\footnote{Indeed, if $xy=1$ then $c=1$ and $u=y$.} that we needed.
}
}

\block{\label{defSpec}Let $A$ be a localic Banach algebra, and let $a \in A$ be a point of $A$, we define $D(a)$ as the sub-locale of $\C$ defined by $\{ \lambda \in \C, (a - \lambda 1) \in A^{\times} \}$ which is clearly an open sub-locale of $\C$. 

\Def{Let $A$ be a localic Banach algebra, we define $\spec[A] a$ to be the closed complement of $D(a)$.}

\Prop{$\spec[A] a$ is a compact sub-locale of $\C$.}

\Dem{$\spec[A] a$ is closed by construction, and it is well known that closed bounded sub-locale of $\C$ are compact so we just have to prove that $\spec a$ is bounded, and we will prove more precisely that for any rational number $q$ such that $\Vert a \Vert <q$, one has $\spec[A] a \subset \{ |\lambda| \leqslant q \}$. 

By taking the open complements, it is the same to show $\{ |\lambda|>q \} \subset D(a)$ and this can be done on generalised elements by the usual power series argument: if $|\lambda| > q$ then the series:

\[ \frac{-1}{\lambda} \sum_{k=0}^{\infty} \left(\frac{a}{\lambda}\right)^{k}\]

converges and is an inverse for $(a - \lambda.1)$ which proves the result.

}

}

\block{Note that the construction of the spectrum is clearly geometric in the sense that if $f : \Ecal \rightarrow \Tcal$ is a geometric morphism between two toposes, $A$ is a localic Banach algebra in $\Tcal$ and $a \in A$ is a point of $A$ in $\Tcal$ then $f^{\sharp}(\spec[A] a) \simeq \spec[f^{\sharp} A] f^{\sharp} a$.

This allows to explain how the notion of spectrum can be developed properly for an arbitrary localic Banach algebra: one can define (by working internally over the locale $A$) the spectrum for the ``generic'' element of a localic Banach algebra which will be a proper separated map $\spec[A] * \rightarrow A$ which factor into a closed inclusion $\spec[A] * \subset \C \times A$, and by geometricity of the spectrum, the fiber over any point $a$ of $A$ is canonically isomorphic to $\spec[A] a$. This will not be used in the present paper.

}

\block{\label{Prop_spectrum_basic}\Prop{Let $A$ be a Banach algebra, and let $a,b\in A$ then :

\[ (\spec ab) \cup \{0 \} = (\spec ba) \cup \{0 \}. \]

}

\Dem{It is the same to prove that $ D(ab) \cap \C^{*} = D(ba) \cap \C^{*} $, and by symmetry it is enough to prove that $D(ab) \cap \C^{*} \subset  D(ba) $ and this can be done on generalised element exactly as in the classical case: if $\lambda$ is invertible and such that $(\lambda - a b)$ is invertible, then let $r$ be the inverse of $(\lambda - ab)$, and a simple algebraic computation using $abr=rab = \lambda r -1$ show that $\lambda^{-1}(1+ b r a)$ is a two sided inverse of $(\lambda - ba)$.}
}

\block{\Prop{Let $A$ be a $C^{*}$-algebra. If $ u \in A$ is a unitary element ($u^{*}u=uu^{*} = 1$) then the spectrum of $u$ is included in the unit circle. If $s \in A$ is self adjoint ($s^{*}=s$) the spectrum of $s$ is included in the real line.}

\Dem{Let $u$ be a unitary element, one has$\Vert u \Vert^{2} = \Vert u^{*} u \Vert =1$  hence the spectrum of $u$ is included in the unit disk. $u^{-1}=u^{*}$ is also unitary hence the spectrum of $u^{-1}$ is also included in the unit disk. Let $ a \in \spec u$ (a generalised element), as $u$ is invertible one has $|a|>0$. Assume $|a|<1$ then $|1/a| >1$ hence $\frac{1}{a} -u^{-1}$ is invertible, but $(u-a) u^{-1} \frac{1}{a} =  \frac{1}{a} -u^{-1}$ hence $(u-a)$ is invertible which is impossible. Hence one has $a \geqslant 1$ and as we already know that $a \leqslant 1$ this proves the first part of the proposition.

Let $s$ be a self-adjoint element. Then $exp(is)$ defined using the exponential power series is a unitary element by a formal computation, hence its spectrum is included in the unit disk. Let $\lambda$ be a complex number such that $|Re(\lambda)|>0$, hence $|exp(i\lambda)| < >1$ hence $exp(is)-exp(i \lambda)$ is invertible. But a computation\footnote{Convergence questions are dealt with exactly as one will do in the classical case} on series gives:

\[exp(is)-exp(i \lambda)= exp(i\lambda) (s-\lambda) \left( \sum_{k=0}^{\infty} \frac{i^{k+1}(s-\lambda)^{k}}{(k+1)!} \right) \]

Hence $(s- \lambda)$ is invertible, i.e. $\lambda \in D(s)$ which concludes the proof.
}

}

\block{\label{propspectrumInv}\Prop{Let $B \subset A$ be a localic Banach algebra with a localic sub-algebra, and let $a \in B$ any point (element) of $B$.
\begin{itemize}
\item[(a)]  $B^{\times}$ is both open and closed in $A^{\times} \cap B$. 
\item[(b)] $D_{B}(a) \subset D_{A}(a)$ is also open and closed.
\end{itemize}

If moreover $A$ is a $C^{*}$-algebra and $B$ a sub-$C^{*}$-algebra, then:

\begin{itemize}
\item[(c)] $B^{\times} = A^{\times} \cap B$.
\item[(d)] $\spec[B] a= \spec[A] a$.
\end{itemize}
}

\Dem{\begin{itemize}
\item[(a)] $B^{\times}$ is open in $B$ hence it is also open in $A^{\times} \cap B$. But on the other hand, $B^{\times}$ is homeomorphic to the set of $(b,b')$ in $B^{2}$ such that $b'b=bb'=1$ while $A^{\times} \cap B$ is homeomorphic to the set of $(b,a)$ in $B \times A$ such that $ab=ba=1$, hence $B^{\times} =( A^{\times} \cap B) \cap (B \times B)$ is clearly a closed subset of $( A^{\times} \cap B)$.

\item[(b)] This follows directly from the previous point and the fact that $D_A(a)$ and $D_B(a)$ are pre-image of $A^{\times}$ and $B^{\times}$ by the map $\lambda \mapsto a- \lambda$ which takes values in $B$.

\item[(c)] On generalized elements: Let $a$ be an element of $B$ invertible in $A$. We want to prove that $a$ is invertible in $B$. It is enough to prove that $a^{*}a$ is invertible in $B$ and $a^{*}a$ is also invertible in $A$ hence one can freely assume that $a$ is self adjoint. As $a$ is invertible in $A$, $0 \in D_A(a)$. moreover, as the spectrum of $a$ is included in the real line one can construct a continuous path in $D_A(a)$ from $0$ to some complex number $z$ of absolute value bigger than $\Vert a \Vert$. As $|z|$ is bigger than $\Vert a \Vert$, $z \in D_B(a)$, but because $D_B(a)$ is open and closed in $D_A(a)$ is has to contain the whole connected component of $z$, in particular it contains $0$ hence $a$ is invertible in $B$.

\item[(d)] This follows directly from the previous point.
\end{itemize}

}
}

\block{From this one easily deduce the following proposition:

\Prop{Let $x$ be a normal\footnote{This means $xx^{*}=x^{*}x$.} element of a $C^{*}$-algebra $A$, then the $C^{*}$-algebra generated by $x$ is isomorphic to the algebra of complex valued continuous function on $\spec x$, with $x$ corresponding to the canonical map from $\spec x$ to $\C$. }

\Dem{Indeed, the $C^{*}$-algebra $C$ generated by $x$ is commutative because $x$ is normal, hence it is isomorphic to the algebra $\Ccal(K)$ of continuous complex valued function on some compact locale $K$ by the localic Gelfand duality. Through this identification $x$ corresponds to some function $f$ on $K$. This function $f:K \rightarrow \C$ has to be an homeomorphism on its image because any function in the algebra generated by $f$ can be seen as a function on the image of $f$ pre-composed by $f$, hence functions on $K$ and on the image of $f$ are the same and as they are both compact completely regular locale this implies that they are isomorphic. Finally, if $f$ is a function on some compact completely regular locale $K$, the spectrum of $f$ is clearly the image of $K$ by $f$ hence this concludes the proof.}

}

\block{\label{Prop_def_positivity}\Prop{Let $x \in A$ be any element of a unital $C^{*}$-algebra. The following conditions are equivalent:
\begin{itemize}
\item[(i)] $x=a^{*}a$ for some $a \in A$.
\item[(ii)] $x = u^{2}$ for some $u \in A$ such that $u^{*}=u$.
\item[(iii)]$x$ is normal and $\spec x \subset \R_+$
\item[(iv)] $x$ is self-adjoint and $\Vert K -x \Vert \leqslant K$ for any continuous number $K \geqslant \frac{1}{2} \Vert x \Vert$.
\item[(v)] $x$ is self-adjoint $\Vert K -x \Vert \leqslant K$ for some continuous number $K$.
\item[(vi)] There is a commutative sub-$C^{*}$-algebra $\Ccal(X) = C \subset A$ with $X$ a compact completely regular locale and $x$ corresponds to a function on $X$ with value in positive real numbers.
\end{itemize}

Elements of $A$ satisfying those conditions are called positive elements, the set of positive element of $A$ is denoted by $A^{+}$ and satisfies:

\begin{itemize}
\item[(a)] If $x \in A^{+}$ and $\lambda \in \R_+$ then $\lambda x \in A^{+}$
\item[(b)] If $x,y \in A^{+}$ then $x+y \in A^{+}$.
\item[(c)] If $x \in A^{+}$ and $-x \in A^{+}$ then $x=0$.
\end{itemize}

}

\Dem{ The proof will go as follow: we will prove that all the conditions from $(ii)$ to $(v)$ are each equivalent to $(vi)$, then we will prove properties $(a),(b)$ and $(c)$ and then we will prove that $(i)$ is equivalent to the other properties.

The key results here is of course proposition \ref{propspectrumInv} which show that the spectrum computed in $A$ or $C$ are the same. From this, one can see that properties $(vi)$ clearly implies all the other, because all these properties holds for a positive function in the algebra of continuous functions. Conversely, assuming $(ii)$ one can take $C$ to be the algebra  generated by $u$, and assuming any of the proposition $(iii)$ to $(v)$ one can take $C$ to be the algebra generated by $x$, and $x$ clearly has to identify with a positive function under any of these assumptions.

The properties $(a)$ and $(c)$ follow easily from characterization $(iii)$ of positivity, and $(b)$ follow from characterization $(v)$: if $\Vert K -x \Vert \leqslant K$ and $\Vert K'- y \Vert \leqslant K'$ then $\Vert K+K' - (x+y)  \Vert \leqslant \Vert K - x \Vert + \Vert K'- y \Vert  \leqslant K+K'$.

Property $(ii)$ clearly implies $(i)$ hence all we have to do is to prove that $a^{*}a$ is positive for any $a \in A$, and we have now all the element to apply the argument of \cite[1.6.4]{dixmierCStar}:

As $a^{*}a$ is self-adjoint one can use functional calculus to decompose it into $a^{*}a = u^{2} - v^{2}$ with $u$ and $v$ self-adjoint such that $uv=vu=0$. Then $(av)^{*}(av)= va^{*}av= v u^{2} v - v^{4} = -v^{4} \in -A^{+}$. Writing $av = s+it$ with $s$ and $t$ self-adjoint one has:
\[ (av)(av)^{*} = -(av)^{*}(av) +(s-it)(s+it) +(s+it)(s-it) = v^{4} + 2s^{2} +2t^{2} \]

Hence $(av)(av)^{*}$ is positive because of $(a)$ and $(b)$. Now $(av)(av)^{*}$ and $(av)^{*}(av)$ have the same spectrum (up to adding zero) by proposition \ref{Prop_spectrum_basic} hence $(av)^{*}(av)$ is also positive (by $(iii)$) and hence $-v^{4}=0$ by $(c)$ so that $v=0$ and $a^{*}a = u^{2}$ which concludes the proof.

}}

\block{\label{Cor_strongineq}All the classical properties of positivity in a $C^{*}$-algebra can then be deduced from those stated in \ref{Prop_def_positivity} and from the Gelfand duality. For example, if $A$ is a $C^{*}$-algebra then let $A^{s}$ be the subset of self-adjoint element of $A$, $A^{s}$ admit an order relation: $x \leqslant y$ if $y-x \in A^{+}$. If $x$ is a positive element, then by working in the commutative $C^{*}$-algebra generated by $x$ one can see that for any continuous real number $r$ one has $x \leqslant r$ if and only if $\Vert x \Vert \leqslant r$ and in particular if $x \leqslant y$ then $\Vert x \Vert \leqslant \Vert y\Vert$. And finally one has:

\Cor{Any $C^{*}$-algebra $A$ satisfies the strong $C^{*}$-inequality: for all $x,y \in A$ one has:

\[ \Vert x \Vert ^{2} \leqslant \Vert x^{*}x +y^{*}y \Vert \]

}
indeed $\Vert x \Vert^{2} = \Vert x^{*} x \Vert \leqslant \Vert x^{*}x +y^{*}y \Vert$. 
}

\block{\Def{Let $a \in A$ be an element of a unital Banach algebra, we denote by $\rho(a)$ or $\rho_A(a)$ the spectral radius of $a$ which is the upper semi-continuous real number defined by:

\[ \rho_A(a) < q \Leftrightarrow \spec[A] a \subset \{ | \lambda | < q \} \]

}
 
The fact that $\rho_A(a)$ defined this way is indeed a upper semi-continuous real number follows easily from the compactness of $\spec a$. Moreover the proof of \ref{defSpec} also show that $\rho(a) \leqslant \Vert a \Vert$.
}

\block{\label{Prop_specctral_radius}\Prop{If $a$ is a normal element of a $C^{*}$-algebra $A$ then $\rho(a) = \Vert a \Vert$ and for any two elements $a,b \in A$ one has $\rho(ab)=\rho(ba)$.}

In particular the norm of any element can be computed in term of spectral radius: $\Vert x \Vert =\rho(x^{*}x)^{1/2}$.

\Dem{If $a$ is normal then the sub-$C^{*}$-algebra generated by $a$ is isomorphic to the algebra of functions over the spectrum of $a$ with $a$ the canonical map to $\C$, hence $\Vert a \Vert <q$ if and only if $\rho(a)<q$. The second part is clear from proposition \ref{Prop_spectrum_basic}.}

}

\block{\Cor{Let $f :A \rightarrow B$ be any morphism of $C^{*}$-algebra then for any element $a \in A$ one has:
\[ \spec f(a) \subset \spec a \]
\[\rho(f(a)) \leqslant \rho(a) \]
\[\Vert f(a) \Vert \leqslant \Vert a \Vert \]

Moreover, when $f$ is injective, all three are equalities.}

\Dem{In terms of generalized points: if $\lambda \in \spec f(a)$ then $f(a) - \lambda$ is non-invertible, hence $a- \lambda$ is also non invertible: indeed, if it was then the image by $f$ of its inverse would be an inverse for $f(a) - \lambda$. This easily implies the relation on $\rho$ and the final relation follow from the fact that $\Vert a \Vert = \Vert a^{*}a \Vert ^{1/2} = \rho(a^{*}a)^{1/2}$. If $f$ is injective then proposition \ref{propspectrumInv} already says that the first inclusion is an equality and the two other equality follows.
}
}

\block{\label{Rk_nonunital}All the results of this section can be adapted to the non-unital case by considering the unitarisation process $A \mapsto A^{+}$ explained in \cite{henry2014nonunital}. There is only one subtleties (that is also present in classical mathematics but easier to ignore in this case):

\Prop{ If $A$ is a unital $C^{*}$-algebra, and $a \in A$, then $\spec[A^{+}] (a) = \spec[A] (a) \cup \{0 \}$. }

In classical mathematics it is not really important because a $C^{*}$-algebra is either unital or non unital and we don't need to unitarize it in the first case, so the only consequence of this is that in the case of non-unital inclusion $B \subset A$ of $C^{*}$-algebra $\spec[A]$ and $\spec[B]$ only coincident when we add zero to them. In constructive mathematics it is a bit more problematic because when we have a general $C^{*}$-algebra we cannot treat it differently depending on if it is unital or not. For all the results that does not explicitly mention the spectrum (like all the results about positivity) this is irrelevant, and for the other possible tricks would be to use the multiplier algebra instead of $A^{+}$, to use a modified construction of $A^{+}$ so that if $A$ is unital then $A=A^{+}$ or to consider two different notions of spectrum: one that always contains zero and one that is only defined for unital algebras. As me will not make any mention of the spectrum of an element in the present paper (outside this appendix) this question will be irrelevant.
}

\blockn{We will conclude this appendix by a very classical results on approximate units in $C^{*}$-algebras.}

\block{\Lem{Let $A$ be a unital $C^{*}$-algebra, $a,b \in A^{+}$ invertible, then:

\begin{itemize}
\item $a \leqslant 1$ if and only if $1 \leqslant a^{-1}$
\item If $a \leqslant b$ then $b^{-1} \leqslant a^{-1}$
\item $x \mapsto (1-x)^{-1} - 1$ induce an order preserving bijection from the set of $x \in A^{+}$ such that $\Vert x \Vert <1$ and $A^{+}$.
\end{itemize}
}

\Dem{From the first point of proposition \ref{Prop_def_positivity}, one easily see that if $x$ is positive than for any $y$, $y^{*}x y$ is positive, hence by definition of the order relation if $a \leqslant b$ then $y^{*} a y \leqslant y^{*}by$ for any $y$. This gives in particular: if $a \leqslant 1$ then $a^{\frac{-1}{2}} a a^{\frac{-1}{2}} \leqslant a^{\frac{-1}{2}}a^{\frac{-1}{2}}$ hence $1 \leqslant a^{-1}$, and similarly in the other direction, moreover if $a \leqslant b$ then $1 \leqslant a^{\frac{-1}{2}} b a^{\frac{-1}{2}}$ hence by the first point, $a^{\frac{1}{2}} b^{-1}a^{\frac{1}{2}} \leqslant 1$ and by multiplying by $a^{\frac{-1}{2}}$ on each side one more time one gets that $b^{-1} \leqslant a^{-1}$ which conclude the proof of the first two points. For the third point, if $x$ is positive and $\Vert x \Vert <1$ it means that there exists $0<q<1$ a rational number such that $0 \leqslant x \leqslant q$, hence $1-q \leqslant 1-x \leqslant 1$ and hence (working in the commutative algebra generated by $x$ and $1$) one has that $1-x$ is invertible hence $(1-x)^{-1} - 1$ indeed exists. Moreover, if $x \leqslant y$ then:

\[  1-x \geqslant 1-y\]
\[ (1-x)^{-1} \leqslant (1-y)^{-1}\]
\[ (1-x)^{-1}- 1 \leqslant (1-y)^{-1} -1\]

Hence this map is order preserving and has zero is sent to $0$ it takes values in $A^{+}$. Finally, it is a bijection because $y \mapsto 1-(y+1)^{-1}$ is an inverse (and if $y \in A^{+}$ then $y+1$ is always invertible, and if $y \leqslant K$ then $1-(y+1)^{-1} \leqslant 1-(K+1)^{-1} <1$.
}}

\block{\label{Prop_approximate_unit}\Prop{Let $A$ be a (possibly non-unital) $C^{*}$-algebra. Then the set $\Lambda_A = \{x \in A^{+}, \Vert x \Vert <1 \}$ ordered with the ordering of $A^{+}$ is directed and for all $x \in A$,

\[ \lim_{\lambda \in \Lambda_A} \lambda x = \lim_{\lambda \in \Lambda_A} x \lambda = x\]

Where the limit is taken in the sense of limit of a net.

}

Note that a net with those properties is generally called an approximate unit. Such nets does not always exist for Banach algebras.

\Dem{When $A$ is unital the fact that $\Lambda_A$ is directed follow from the last point of the previous lemma which shows that it is in order preserving bijection with $A^{+}$ which is clearly directed: for any $a,b \in A^{+}, a,b \leqslant a+b$. In the general case one can take $\tilde{A}$ be some unital algebra in which $A$ is a two-sidded ideal (so $\tilde{A}$ can be either the unitarization defined in \cite{henry2014nonunital} or the multiplier algebra defined in the present paper) and observe that this bijection preserves any two sided ideal of $A$ (as it is defined by continuous functional calculus with respect to a function that send $0$ to $0$, see lemma \ref{Lem_hered}(\ref{Lem_hered2}) ), hence the same argument applies.

For the limit, let $x \in A$, we will first prove that $x^{*} \lambda x$ converge to $x^{*}x$. let $K$ be a rational number such that $\Vert x \Vert < K$, then by lemma \ref{Lem_hered}(\ref{Lem_hered3}) one has that the sequences $x^{*}(x^{*}x/K^{2})^{1/n} x$ converge to $x^{*}x$, moreover the sequence $(x^{*}x/K^{2})^{1/n}$ is a subsequence of the net $\Lambda_A$, so we at least have converging subsequence. Moreover  if $a \leqslant b \leqslant 1$ and if $\Vert x^{*}x -  x^{*} ax \Vert < q$ then $\Vert x^{*}x - x^{*} bx \Vert < q$.

Indeed, $ 0 \leqslant x^{*}(1-b)x  \leqslant x^{*}(1-a)x$ hence $ \Vert x^{*}(1-b)x \Vert \leqslant \Vert x^{*}(1-a)x \Vert$ so this proves the convergence of $x^{*} \lambda x$.

For the case of $\lambda x$:
\[\Vert (1- \lambda) x \Vert^{2} = \Vert x^{*} (1- \lambda)^{2} x \Vert \leqslant \Vert x^{*}(1-\lambda) x \Vert \]

where the last inequality is just the fact that as $0 \leqslant 1 - \lambda \leqslant 1$ one has $(1-\lambda)^{2} \leqslant 1 - \lambda $.

}

}

\section{Foundation: Heyting pretopos and class of small maps}
\label{appendiw_foundation}
\blockn{As mentioned in the introduction the general framework of this paper is the internal logic of a Heyting pretopos with a natural number object and a class of small maps satisfying all the additional axioms $(S3),(S4)$ and $(S5)$ of \cite{joyal1995algebraic}. In this section we will recall the definition and theses axioms, detail a little what it mean in terms of the internal logic and explain what are the intended examples of such categories.}

\blockn{A Heyting pretopos with a natural number object is essentially a category in which one can interpret first order intuitionist logic and which has an object $\N$ of natural number. For the concrete definition we refer the reader either to the appendix $B$ of \cite{joyal1995algebraic} or to \cite[A1.4]{sketches} for a more systematic introduction. For the definition of a natural number object we refer to \cite[A2.5]{sketches}.}

\blockn{We recall from chapter I\S 1 of \cite{joyal1995algebraic}:}

\block{\label{Def_class_of_smallmaps}\Def{Let $\Ccal$ be a Heyting pretopos, a class $\Scal$ of map of $\Ccal$ is called a class of open maps if it satisfies the following axioms:

\begin{itemize}

\item[(A1)] Any isomorphism belongs to $\Scal$ and $\Scal$ is closed under composition.

\item[(A2)] (\emph{Stability}) $\Scal$ is stable under pullback.

\item[(A3)] (\emph{Descent}) If the pullback of $f$ along an epimorphism is in $\Scal$ then $f$ is in $\Scal$.

\item[(A4)] The maps $\emptyset \rightarrow *$ from the initial object to the terminal object and the map $* \coprod * \rightarrow *$ are in $\Scal$.

\item[(A5)](\emph{Sums}) If $f:X \rightarrow Y$ and $f' : X' \rightarrow Y'$ are both in $\Scal$ then $f \coprod g' : X \coprod X' \rightarrow Y \coprod Y'$ is in $\Scal$.

\item[(A6)] (\emph{Quotient}) If $p$ is an epimorphism and $f \circ p \in \Scal$ then $f \in \Scal$.

\item[(A7)] (\emph{Collection axiom}) For any two arrow $p: Y \twoheadrightarrow X$ and $f:X \rightarrow A$ with $p$ an epimorphism and $f \in \Scal$ then there exists a square:

\[
\begin{tikzcd}[ampersand replacement=\&]
Z \arrow{r} \arrow{d}{g} \& Y \arrow[twoheadrightarrow]{r}{p} \& X \arrow{d}{f} \\
B \arrow[twoheadrightarrow]{rr}{h} \& \& A \\
\end{tikzcd}
\]

Such that $g \in \Scal$, $h$ is an epimorphism and the canonical map $Z \rightarrow B \times_A X$ is an epimorphism.

\end{itemize}

It is said to be a class of small maps if it satisfies the following additional axioms:

\begin{itemize}

\item[(S1)] (\emph{Exponentiability Axiom}) If $f :Y \rightarrow X$ is in $\Scal$ then $Y$ is exponentiable\footnote{An object $A \in C$ is said to be exponentiable if for any $B \in |C|$ there exists an object $[A,B]$ satisfying the universal property $C(X,[A,B]) \simeq C(X \times A,B)$.} in $\Ccal/X$.

\item[(S2)] (\emph{Representability axiom}) There exists a map $\pi:E \rightarrow U$ in $\Scal$ which is universal in the following sense: for any map $f : Y \rightarrow X \in \Scal$ there exists a diagram:

\[
\begin{tikzcd}[ampersand replacement=\&]
Y \arrow{d}{f} \& \arrow{l} Y' \arrow{d}{f'} \arrow{r} \& E \arrow{d}{\pi} \\
X \& \arrow[twoheadrightarrow]{l}{p} X' \arrow{r}{c} \& U \\
\end{tikzcd}
\]

in which $p$ is an epimorphism and the two square are pullback square.

\item[(S3)](\emph{Power-Set}) For every $X \in \Ccal$ the category $\Scal/X$ of arrows $Y \rightarrow X$ which are in $\Scal$ is an elementary topos\footnote{\label{footnote_S3_Ps}It is proved in \cite[I.\S 3]{joyal1995algebraic} that if $Y \in \Ccal/X$ there is an object $P_s(Y)$ of $\Ccal/X$ which classifies sub-objects of $Y$ which are in $\Scal/X$, hence this axiom can be stated as if $Y \in \Scal/X$ then $P_s(Y) \in \Scal/X$. }.

\item[(S4)](\emph{Separation Axiom}) Every monomorphism is in $\Scal$.

\item[(S5)](\emph{Infinity}) The map $N \rightarrow * \in \Scal$ where $N$ is the natural number object of $\Ccal$.

\end{itemize}

}

Arrows in $\Scal$ are called small maps, and objects whose unique map to $*$ is in $\Scal$ are called small objects.

}

\blockn{We will now discuss these axioms from the perspective of internal logic. First of all an Heyting pretopos is essentially a category where one can interpret first order logic, but that might not have a sub-object classifier, power sets or functions sets. Axioms $(A2),(A3)$ (and partly $(S2)$) ensure that the notions of small map and small object can be discussed internally (see \cite[I.1.6]{joyal1995algebraic}), and that internally a map $f :Y \rightarrow X$ is small if and only if internally $\forall x \in X, f^{-1} \{x\}$ is small.

\bigskip

Internally we will call ``class'' the object of $\Ccal$. Axiom $(S1)$ ensure that if if $X$ is a small class and $Y$ a class then there is a class of functions from $X$ to $Y$, if $Y$ is small then (using $(S3)$) we can prove that this class of function is itself small.

\bigskip

Axiom $(S2)$ gives us a ``weak universe''. The class $U$ is ``morally'' the class of all small class, and $E$ is the class of pair $(u \in U, x \in u)$. This only works up to the point that it would not make sense to says that any small class is equal to a small class of the form $E_u =\{e \in E | \pi(e)=u \}$ for $u \in U$ because equality of small class does not make sense internally.

Internally, we will call ``\emph{set}'' the elements of $U$, and identifies them with the corresponding small class ``$E_x$'', this way $U$ is indeed the class of all sets by definition and a small class is a class which is internally in bijection with a set.

This choice of terminology has some drawback and may seem unnatural from the perspective of categorical logic, but this seems to be the only\footnote{One could also simplify this distinction by introducing a variant the univalence axiom for our universe but this would force us to use higher categories instead of ordinary categories as the universe would be a groupoid.} way to be able to talk of the class of all sets, or to quantify over sets latter.
}

\block{Axiom $(A1)$ tell us that $*$ is small and that a (disjoint) union of small sets indexed by a small set is again small, axiom $(A4)$ that the empty set and the two element set are both small (together with $(A1)$ this show that the cardinal finite class $\{1,\dots...,n\}$ are small).

Axioms $(A2),(A3),(A5)$ are tautological when stated in term of the internal logic, but $(A2)$ and $(A3)$  where of course required to show that the notion ``small'' makes sense internally, and it appears that $(A5)$ follows from $(A2),(A3)$ and $(S2)$.

Axioms $(A6)$ say that if a class can be covered by a small class then it is small hence one can for example deduce that the Kuratowski\footnote{See \cite[D5.4]{sketches}, only the equivalent definition involving the natural number object makes sense here.} finite class are small.}

\block{Axiom $(S4)$ says that a subclass of a small class is small. It corresponds to the ``separation axiom'' of set theory only if you consider that a class is exactly something that can be defined by a predicate. Axiom $(S5)$ just say that the ``class of natural number'' is small, hence corresponds to the ordinary axiom of infinity.}

\block{Axiom $(S3)$ is probably the more subtle. Given a class $X$, one can (using the universe $U$ granted to us by the representability axiom) construct the class of all small set endowed with an injection to $X$, and then construct its quotient by the relation ``equality of image'',this gives the class of all small subclass of $X$, this is the object $P_s(X)$ mentioned in footnote \ref{footnote_S3_Ps} above. Axiom $(S3)$ assert that when $X$ is small this class is small, hence for any set $X$ one has a set (a small class in fact) of subset of $X$, making the category of sets an elementary topos, hence allowing to use full impredicative Higher order logic for sets. }

\block{\label{Description_Collection}Finally, there is the collection axiom $(A7)$, internally it corresponds to: if $p:V \twoheadrightarrow X$ is a surjection with $V$ a class and $X$ a small class then there exists a set $Y$ and a map from $f:Y \rightarrow V$ such that $p \circ f : Y \rightarrow X$ is already surjective. The only places in the present paper where this axiom is used in an essential way is in proposition \ref{Prop_restrictionaresmall} and its corollary \ref{Cor_restrictionfctexist}.}

\block{We conclude this appendix by providing examples of such categories. The first example is any model of Neumann-Bernays-Gödel\footnote{ Abbreviated as $(NBG)$. It is an axiomatization of a theory of class and sets with axioms very similar to our framework and which is a conservative extension of $(ZF)$, i.e. a statement in the language of $(ZF)$ can be proved in $(ZF)$ if and only if it can be proved in $(NBG)$.} set theory. Given such a model, the category of all class, with $\Scal$ the family of map between class whose fiber are sets, satisfies all of the above axioms. This is probably the simplest way to obtain such an example as a model of $(NBG)$ can be constructed out of any model of $(ZF)$ (essentially by defining the class as being the definable predicate). In this case $\Ccal$ itself is not an elementary topos (or at least does not have to be one), nor even a cartesian closed category, so it is not possible to talk about the class of all functions between two class.

\bigskip

Similarly, from a model of $(ZFC)$ with an inaccessible cardinal $\kappa$ one can also construct such an example by taking $\Ccal$ to be the category of all sets and $\Scal$ the class of map whose fiber have cardinality smaller than $\kappa$ (and $U$ will be the set all subsets of some representative of $\kappa$). In this case $\Ccal$ itself is a topos and one can speak about the class of all subclass of a given class or the class of functions between two class internally, but it require stronger set theoretic foundation than the previous example (the existence of an inaccessible cardinal being unprovable in $(ZFC)$). A Grothendieck universe would do exactly the same.

\bigskip

Of course if these where the only models in which we are interested in then this categorical framework and avoiding the law of excluded middle and the axiom of choice every where would have been useless. What we are interested in are sheaves models:
}

\block{Assume we are working internally in such a category $\Ccal$ satisfying all the axioms (so for example either a model of $(NBG)$ or a model of $(ZFC)$ with an inaccessible cardinal). Let $(C,J)$ be a small site, and take $\Ccal'$ to be the (meta)category of class valued sheaves over $(C,J)$ and $\Scal'$ be the class of all map $f$ in $\Ccal'$ such that for each $c \in C$, $f_c$ is a small map.

Then this satisfies all the axiom $(A1-7)$ and $(S1-5)$ (one need $(C,J)$ to be small only for the power set axiom $(S3)$) this is explained in \cite[IV.3]{joyal1995algebraic} in the case where the ground model is a model of $(ZFC)$ with an inaccessible cardinal but it works all the same in the general framework.

\bigskip

The category of small class is (equivalent to) the ordinary category of set valued sheaves on our site and class are class valued sheaves. This is the key example that we have in mind throughout all the paper.
}

\block{Let us also briefly mention that there is a ``realizability'' model discussed in \cite[IV.4]{joyal1995algebraic}, showing that theorems proved in this framework have a computational content. }

\bibliography{Biblio}{}

\begin{thebibliography}{10}

\bibitem{banaschewski2006globalisation}
Bernhard Banaschewski and Christopher~J Mulvey.
\newblock A globalisation of the {G}elfand duality theorem.
\newblock {\em {Annals of Pure and Applied Logic}}, 137(1):62--103, 2006.

\bibitem{bourbaki1966elements}
Nicolas Bourbaki.
\newblock {\em Elements of {M}athematics: {G}eneral {T}opology}.
\newblock Hermann, 1966.

\bibitem{burden1979banach}
CW~Burden and CJ~Mulvey.
\newblock {Banach spaces in categories of sheaves}.
\newblock In {\em Applications of sheaves}, pages 169--196. Springer, 1979.

\bibitem{coquand2009constructive}
Thierry Coquand and Bas Spitters.
\newblock Constructive {G}elfand duality for ${C}^*$-algebras.
\newblock In {\em {Mathematical Proceedings of the Cambridge Philosophical
  Society}}, volume 147, pages 323--337. {Cambridge Univ Press}, 2009.

\bibitem{dixmierCStar}
Jacques Dixmier.
\newblock {\em Les $C^*$-alg{\`e}bres et leurs repr{\'e}sentations}.
\newblock {Editions Jacques Gabay}, 1969.

\bibitem{wstarcat}
P~Ghez, Ricardo Lima, and John~E Roberts.
\newblock ${W}^*$-categories.
\newblock {\em {Pacific J. Math}}, 120(1):79--109, 1985.

\bibitem{henry2014nonunital}
Simon Henry.
\newblock Constructive {G}elfand duality for non-unital commutative
  ${C}^{*}$-algebras.
\newblock {\em {arXiv preprint arXiv:1412.2009}}, 2014.

\bibitem{henry2014localic}
Simon Henry.
\newblock Localic metric spaces and the localic {G}elfand duality.
\newblock {\em {arXiv preprint arXiv:1411.0898}}, 2014.

\bibitem{henry2014measure}
Simon Henry.
\newblock Measure theory over boolean toposes.
\newblock {\em {arXiv preprint arXiv:1411.1605}}, 2014.

\bibitem{henry2015finitness}
Simon Henry.
\newblock On toposes generated by cardinal finite objects.
\newblock {\em arXiv preprint arXiv:1505.04987}, 2015.

\bibitem{henry2015toward}
Simon Henry.
\newblock Toward a non-commutative {G}elfand duality: {B}oolean locally
  separated toposes and monoidal monotone complete ${C}^*$-categories.
\newblock {\em {arXiv preprint arXiv:1501.07045}}, 2015.

\bibitem{sketches}
P.T. Johnstone.
\newblock {\em Sketches of an elephant: a topos theory compendium}.
\newblock Clarendon Press, 2002.

\bibitem{joyal1984extension}
A.~Joyal and M.~Tierney.
\newblock {\em An extension of the {G}alois theory of {G}rothendieck}.
\newblock {American Mathematical Society}, 1984.

\bibitem{joyal1995algebraic}
André Joyal and Ieke Moerdijk.
\newblock {\em Algebraic set theory}, volume 220.
\newblock Cambridge University Press, 1995.

\bibitem{lance1995hilbert}
E~Christopher Lance.
\newblock {\em Hilbert ${C}^*$-modules: a toolkit for operator algebraists},
  volume 210.
\newblock Cambridge University Press, 1995.

\bibitem{makkai1996avoiding}
Michael Makkai.
\newblock Avoiding the axiom of choice in general category theory.
\newblock {\em Journal of Pure and Applied Algebra}, 108(2):109--173, 1996.

\bibitem{moerdijk2000proper}
Ieke Moerdijk and Jacob Johan~Caspar Vermeulen.
\newblock {\em Proper maps of toposes}, volume 705.
\newblock AMS Bookstore, 2000.

\bibitem{mulvey1980banach}
Christopher~J Mulvey.
\newblock Banach sheaves.
\newblock {\em {Journal of Pure and Applied Algebra}}, 17(1):69--83, 1980.

\bibitem{paschke1974double}
William~L Paschke.
\newblock The double {B}-dual of an inner product module over a ${C}^*$-algebra
  {B}.
\newblock {\em Can. J. Math}, 26(5):1272--1280, 1974.

\bibitem{phillips1988inverse}
N~Christopher Phillips.
\newblock Inverse limits of ${C}^*$-algebras.
\newblock {\em J. Operator Theory}, 19(1):159--195, 1988.

\bibitem{picado2012frames}
Jorge Picado and Ale{\'e}s Pultr.
\newblock {\em Frames and {L}ocales: topology without points}.
\newblock {Springer}, 2012.

\bibitem{renault1980groupoid}
Jean Renault.
\newblock {\em A groupoid approach to $C^*$-algebras}.
\newblock Berlin; New York: Springer-Verlag, 1980.

\bibitem{vasselli2007bundles}
Ezio Vasselli.
\newblock Bundles of ${C}^*$-categories.
\newblock {\em Journal of Functional Analysis}, 247(2):351--377, 2007.

\end{thebibliography}
\bibliographystyle{plain}

\end{document}